\documentclass[openany]{memo-l}
\usepackage{amssymb}

\input{tcilatex}

\renewcommand{\theequation}{\thechapter.\arabic{equation}}

\def\iint{\mathop{\displaystyle\int\!\!\!\int}}

\def\iint{\mathop{\displaystyle\int\!\!\!\int}}
\def\text{\mbox}

\begin{document}

\begin{center}
\mathstrut

{\huge \textbf{Operator Valued Hardy Spaces \medskip }}

\medskip

{\Large T{\small AO} M}{\small EI}
\end{center}

\medskip

{\large Author address: }

{\large \medskip Department of Mathematics, Texas A\&M University, College
Station, TX, 77843, U. S. A. }

{\large \medskip Email address: tmei@math.tamu.edu } 

{\large \medskip Jan. 01, 2004} 

\newpage 

\tableofcontents
 
\footnotetext{{\small \textbf{Key words }\ Hardy space,\ BMO space,\
Hardy-Littlewood maximal function,\ von Neumann algebra, non-commutative $%
L_p $ space, interpolation, Lusin integral. }
\par
{\small \textbf{2000 MR Subject Classification.} \ 46L52, 32C05.}}

\begin{center}

\textbf{Abstract}
\end{center}

\medskip {We give a systematic study on the Hardy spaces of functions with
values in the non-commutative $L^p$-spaces associated with a semifinite von
Neumann algebra $\mathcal{M}.$ This is motivated by the works on matrix
valued Harmonic Analysis (operator weighted norm inequalities, operator
Hilbert transform), and on the other hand, by the recent development on the
non-commutative martingale inequalities. Our non-commutative Hardy spaces
are defined by the non-commutative Lusin integral function. The main results
of this paper include:}

\begin{enumerate}
\item[(i)]  The analogue in our setting of the classical Fefferman duality
theorem between $\mathcal{H}^1$ and $\mathrm{BMO}$.

\item[(ii)]  The atomic decomposition of our non-commutative $\mathcal{H}^1.$

\item[(iii)]  The equivalence between the norms of the non-commutative Hardy
spaces and of the non-commutative $L^p$-spaces $(1<p<\infty ).$

\item[(iv)]  The non-commutative Hardy-Littlewood maximal inequality.

\item[(v)]  A description of BMO as an intersection of two dyadic BMO.

\item[(vi)]  The interpolation results on these Hardy spaces.
\end{enumerate}
\newpage

\begin{center}
{\Large \textbf{Introduction}}
\end{center}

\medskip \setcounter{theorem}{0}\setcounter{equation}{0} This paper gives a
systematic study of matrix valued (and more generally, operator valued)
Hardy spaces. Our motivations come from two closely related directions. The
first one is matrix valued Harmonic Analysis. It consists in extending
results from classical Harmonic Analysis to the operator valued setting. We
should emphasize that such extensions not only are interesting in themselves
but also have applications to other domains such as prediction theory and
rational approximation. A central subject in this direction is the study of
''operator valued'' Hankel operators (i.e. Hankel matrices with operator
entries). As in the scalar case, this is intimately linked to the operator
valued weighted norm inequalities, operator valued Carleson measures,
operator valued Hardy spaces.... A lot of works have been done notably by F.
Nazarov, S. Treil and A. Volberg; see, for instance, the recent works \cite
{[GPTV]}, \cite{NPTV}, \cite{[P]}, \cite{[Pe]}, \cite{PS}).

The second direction which motivates this paper is the non-commutative
martingale theory. This theory had been initiated already in the 70's. For
example, I. Cuculescu (\cite{C}) proved a non-commutative analogue of the
classical Doob weak type (1,1) maximal inequality. This has immediate
applications to the almost sure convergence of non-commutative martingales
(see also \cite{JR}, \cite{JR2}). The new input into the theory is the
recent development on the non-commutative martingale inequalities. This has
been largely influenced and inspired by the operator space theory. Many
inequalities in the classical martingale theory have been transferred into
the non-commutative setting. These include the non-commutative
Burkholder-Gundy inequalities, the non-commutative Doob inequality, the
non-commutative Burkholder-Rosenthal inequalities and the boundedness of the
non-commutative martingale transforms (see \cite{[PX2]}, \cite{J}, \cite{JM}%
, \cite{JX}, \cite{Rm}).

One common important object in the two directions above is the
non-commutative analogue of the classical \textrm{BMO} space. Because of the
non-commutativity, there are now two non-commutative \textrm{BMO }spaces,
the column \textrm{BMO} and row \textrm{BMO. }As expected, these
non-commutative \textrm{BMO} spaces are proved to be the duals of some
non-commutative $H^1$ spaces. To be more precise and to go into some
details, we introduce these spaces in the case of matrix valued functions.
Let $\mathcal{M}_d$\ be the algebra of $d\times d$\ matrices with its usual
trace $tr.$ Then the column \textrm{BMO }space is defined by 
\[
\rm{BMO}_c(\Bbb{R},\mathcal{M}_d)=\left\{ {\varphi} :\Bbb{R}\rightarrow 
\mathcal{M}_d,\left\| {\varphi} \right\| _{\mathrm{BMO}_c}<\infty \right\} 
\]
where 
\[
\left\| {\varphi} \right\| _{\mathrm{BMO}_c}=\sup_h\left\{ \left\| {\varphi}
(\cdot )h\right\| _{\mathrm{BMO}(l_2^d)},h\in l_2^d,\left\| h\right\|
_{l_2^d}\leq 1\right\} . 
\]
Similarly, the row \textrm{BMO} space is 
\[
\mbox{BMO}_r(\Bbb{R},\Bbb{\mathcal{M}}_d)=\left\{ {\varphi} :\Bbb{%
R\rightarrow \mathcal{M}}_d,\left\| {\varphi} \right\| _{\mathrm{BMO}%
_r}=\left\| {\varphi} ^{*}\right\| _{\mathrm{BMO}_c}<\infty \right\}. 
\]
We will also need the intersection of these \textrm{BMO} spaces, which is 
\[
\mbox{BMO}_{cr}(\Bbb{R},\Bbb{\mathcal{M}}_d)=\mbox{BMO}_c(\Bbb{R},\mathcal{M}%
_d)\cap \mbox{BMO}_r(\Bbb{R},\Bbb{\mathcal{M}}_d) 
\]
equipped with the norm $\left\| {\varphi} \right\| _{\mathrm{BMO}_{cr}}=\max
\{\left\| {\varphi} \right\| _{\mathrm{BMO}_c},\left\| {\varphi} \right\| _{%
\mathrm{BMO}_r}\}.$ When $d=1,$\ all these BMO spaces coincide with the
classical BMO\ space which is well known to be the dual of the classical
Hardy space $H^1.$ This result can be extended to the case of $d<\infty $\
very easily. Let 
\[
H^1(\Bbb{R},\Bbb{\mathcal{S}}_d^1)=\left\{ f:\Bbb{R\rightarrow \mathcal{S}}%
_d^1;\int \sup_{y>0}\left\| f(x,y)\right\| _{\mathcal{S}_d^1}dx<\infty
\right\} , 
\]
where $\Bbb{\mathcal{S}}_d^1\,$is the trace class over $l_d^2,$ and $f(x,y)$%
\ denotes the Poisson integral of $f$ corresponding to the point $x+iy.\,$
Then 
\[
(H^1(\Bbb{R},\mathcal{S}_d^1))^{*}=\mbox{BMO}_{cr}(\Bbb{R},\mathcal{M}_d) 
\]
\ and 
\[
c_d^{-1}\left\| {\varphi} \right\| _{\mathrm{BMO}_{cr}(\Bbb{R},\mathcal{M}
_d)}\leq \left\| {\varphi} \right\| _{(H^1(\Bbb{R},\Bbb{\mathcal{S}}%
_d^1))^{*}}\leq c_d\left\| {\varphi} \right\| _{\mathrm{BMO}_{cr}(\Bbb{R},
\mathcal{M}_d)}. 
\]
Here the constant $c_d\rightarrow +\infty $\ as $d\rightarrow +\infty .$\
Thus this duality between $H^1(\Bbb{R},\Bbb{\mathcal{S}}_d^1)$\ and $\mathrm{%
BMO}_{cr}(\Bbb{R},\Bbb{\mathcal{M}}_d)$\ fails for the infinite dimensional
case. One of our goals is to find a natural predual space of $\mathrm{BMO}
_{cr}$\ with relevant constants independent of $d.$

In the case of non-commutative martingales, this natural dual of \textrm{BMO}%
$_{cr}$ has been already introduced by Pisier and Xu in their work on the
non-commutative Burkholder-Gundy inequality. To define the right space $
\mathcal{H}^1,$\ they considered a non-commutative analogue of the classical
square function for martingales. Motivated by their work, we will introduce
a new definition of $H^1$\ for matrix valued functions by considering a {%
non-commutative analogue of the classical Lusin integral (Recall that, in the
classical case, a scalar valued function is in }$H^1$ {if and only if its
Lusin integral is in }$L^1,$ see \cite{fs}, \cite{[St]}$).$ For matrix valued
function $f, f\in L^1((\Bbb{R},\frac {dt}{1+t^2}),\mathcal{M}_d), 1\leq p<\infty ,$ 
let 
\[
\left\| f\right\| _{\mathcal{H}_c^p(\Bbb{R},\mathcal{M}_d)}^p=tr\int_{-%
\infty }^{+\infty }(\mathop{\displaystyle\int\!\!\!\int}_\Gamma |\nabla
f(t+x,y)|^2dxdy)^{\frac p2}dt, 
\]
where $\Gamma =\{(x,y)\in \Bbb{R}:|x|<y,y>0\}$ and 
\[
|\nabla f|^2=(\frac{\partial f}{\partial x})^{*}\frac{\partial f}{\partial x}%
+(\frac{\partial f}{\partial y})^{*}\frac{\partial f}{\partial y}. 
\]
Then we define 
\[
\mathcal{H}_c^p(\Bbb{R},\mathcal{M}_d)=\left\{ f:\Bbb{R\rightarrow \mathcal{M%
}}_d;\left\| f\right\| _{\mathcal{H}_c^p(\Bbb{R},\mathcal{M}_d)}<\infty
\right\} . 
\]
Similarly, set 
\[
\mathcal{H}_r^p(\Bbb{R},\mathcal{M}_d)=\left\{ f:\Bbb{R\rightarrow \mathcal{M%
}}_d;\left\| f\right\| _{\mathcal{H}_r^p(\Bbb{R},\mathcal{M}_d)}=\left\|
f^{*}\right\| _{\mathcal{H}_c^p(\Bbb{R},\mathcal{M}_d)}<\infty \right\} . 
\]
Finally, if $1\leq p<2,$ we define 
\[
\mathcal{H}_{cr}^p(\Bbb{R},\mathcal{M}_d)=\mathcal{H}_c^p(\Bbb{R},\mathcal{M}%
_d)+\mathcal{H}_r^p(\Bbb{R},\mathcal{M}_d) 
\]
equipped with the norm 
\[
\left\| f\right\| _{\mathcal{H}_{cr}^p(\Bbb{R},\mathcal{M}_d)}=\inf
\{\left\| g\right\| _{\mathcal{H}_c^p}+\left\| h\right\| _{\mathcal{H}%
_r^p}:f=g+h,g\in \mathcal{H}_c^p(\Bbb{R},\mathcal{M}_d),h\in \mathcal{H}_r^p(%
\Bbb{R},\mathcal{M}_d)\}. 
\]
If $p\geq 2,$ let 
\[
\mathcal{H}_{cr}^p(\Bbb{R},\mathcal{M}_d)=\mathcal{H}_c^p(\Bbb{R},\mathcal{M}%
_d)\cap \mathcal{H}_r^p(\Bbb{R},\mathcal{M}_d) 
\]
equipped with the norm 
\[
\left\| f\right\| _{\mathcal{H}_{cr}^p(\Bbb{R},\mathcal{M}_d)}=\max
\{\left\| f\right\| _{\mathcal{H}_c^p(\Bbb{R},\mathcal{M}_d)},\left\|
f\right\| _{\mathcal{H}_r^p(\Bbb{R},\mathcal{M}_d)}\}. 
\]
One of our main results is the identification of $\mathrm{BMO}_c(\Bbb{R},%
\mathcal{M}_d)$ as the dual of $\mathcal{H}_c^1(\Bbb{R},\mathcal{M}_d):(%
\mathcal{H}_c^1(\Bbb{R},\mathcal{M}_d))^{*}=\mathrm{BMO}_c(\Bbb{R},\mathcal{M%
}_d)$ with equivalent norms, where the relevant equivalence constants are
universal. Similarly, $\mathrm{BMO}_r(\Bbb{R},\mathcal{M}_d)$ (resp. $%
\mathrm{BMO}_{cr}(\Bbb{R},\mathcal{M}_d))$ is the dual of $\mathcal{H}_c^1(%
\Bbb{R},\mathcal{M}_d)($resp. $\mathcal{H}_{cr}^1(\Bbb{R},\mathcal{M}_d)).$
Another result is the equality $\mathcal{H}_{cr}^p(\Bbb{R},\mathcal{M}%
_d)$ $=L^p(L^\infty (\Bbb{R})\otimes \mathcal{M}_d)$ with equivalent norms for
all $1<p<\infty .$ This is the function space analogue of the
non-commutative Burkholder-Gundy inequality in \cite{[PX2]}. It is also
closely related to the recent work (\cite{JMX}) by Junge, Le Merdy and Xu on
the Littlewood-Paley theory for semigroups on non-commutative $L^p$-spaces.

We also prove the analogue of the classical Hardy-Littlewood maximal
inequality. Our approach to this inequality for functions consists in
reducing it to the same inequality for dyadic martingales. We should
emphasize that this approach is new even in the scalar case. It is extremely
simple. The same idea allows to write \textrm{BMO }as an intersection of two
dyadic \textrm{BMO. }This latter result plays an important role in this
paper. It permits to reduce many problems involving \textrm{BMO} (or its
variant \textrm{BMO}$^q,$ which is the dual of $\mathcal{H}^p$ for $1\leq
p<2,\frac 1p+\frac 1q=1$) to dyadic \textrm{BMO, }that is, to \textrm{BMO }%
of dyadic non-commutative martingales. For instance, this is the case of the
interpolation problems on our non-commutative Hardy spaces.

All results mentioned above remain valid for a general semifinite von
Neumann algebra $\mathcal{M}$ in place of the matrix algebras.$\,$

We now explain the organization of this paper. Chapter 1 (the next one)
contains preliminaries, definitions and notations used throughout the paper.
There we define the two non-commutative square functions which are the
non-commutative analogues of the Lusin area integral and Littlewood-Paley $g$%
-function. These square functions allow to define the corresponding
non-commutative Hardy spaces $\mathcal{H}_c^p(\Bbb{R},\mathcal{M}),$ where $%
\mathcal{M}$ is a semifinite von Neumann algebra. This chapter also contains
the definition of $\mathrm{BMO}_c(\Bbb{R},\mathcal{M})$ and some elementary
properties of these spaces.

The main result of Chapter 2 is the analogue in our setting of the famous
Fefferman duality theorem between $\mathcal{H}^1$ and $\mathrm{BMO}.$ As in
the classical case, this result implies an atomic decomposition for our
Hardy spaces $\mathcal{H}_c^1(\Bbb{R},\mathcal{M})$ (as well as $\mathcal{H}%
_r^1(\Bbb{R},\mathcal{M}),\mathcal{H}_{cr}^1(\Bbb{R},\mathcal{M})).$ Another
consequence is the characterization of functions in $\mathrm{BMO}_c(\Bbb{R},%
\mathcal{M})$ (as well as $\mathrm{BMO}_r(\Bbb{R},\mathcal{M}),\mathrm{BMO}%
_{cr}(\Bbb{R},\mathcal{M}))$ via operator valued Carleson measures.

The objective of Chapter 3 is the non-commutative Hardy-Littlewood maximal
inequality. As already mentioned above, our approach to this is to reduce
this inequality to the corresponding maximal inequality for dyadic
martingales. To this end, we construct two ''separate'' increasing
filtrations $\mathcal{D=}\{\mathcal{D}_n\}_{n\in \Bbb{Z}}$ and $\mathcal{D}%
^{\prime }\mathcal{=}\{\mathcal{D}_n^{\prime }\}_{n\in \Bbb{Z}}$ of dyadic $%
\sigma $-algebras. One of them is just the usual dyadic filtration on $\Bbb{R%
};$ while the other is a kind of translation of the first. The main point is
that any interval of $\Bbb{R}$ is contained in one atom of some $\sigma $%
-algebra of them with comparable size. This approach will be repeatedly used
in the subsequent chapters. We also prove the non-commutative Poisson
maximal inequality and the non-commutative Lebesgue differentiation theorem.

In Chapter 4, we define the $L^p$-space analogues of the $\mathrm{BMO}$
spaces introduced in Chapter 1, denoted by $\mathrm{BMO}_c^q(\Bbb{R},%
\mathcal{M}),$ $\mathrm{BMO}_r^q(\Bbb{R},\mathcal{M}),$ $\mathrm{BMO}_{cr}^q(%
\Bbb{R},\mathcal{M}).$ These spaces are proved to be the duals of the
respective Hardy spaces $\mathcal{H}_c^p(\Bbb{R},\mathcal{M}),$ $\mathcal{H}%
_r^p(\Bbb{R},\mathcal{M}),$ $\mathcal{H}_{cr}^p(\Bbb{R},\mathcal{M})$ for $%
1<p<2$ ($q=\frac p{p-1}).$ The proof of this duality is also valid for $p=1.$
In that case, we recover the duality theorem in Chapter 2. However, for $%
1<p<2,$ we need, in addition, the non-commutative maximal inequality from
Chapter 3. This is one of the two reasons why we have decided to present
these two duality theorems separately. Another is that the reader may be
more familiar with the duality between $H^1$ and \textrm{BMO} and those only
interested in this duality can skip the case $1<p<2.$ It is also proved in
this chapter that $\mathrm{BMO}_c^q(\Bbb{R},\mathcal{M})=\mathcal{H}_c^q(%
\Bbb{R},\mathcal{M})$ with equivalent norms for all $2<q<\infty .$ The third
result of Chapter 4 is the following: Regarded as a subspace of $%
L^p(L^\infty (\Bbb{R)}\otimes \mathcal{M},L_c^2(\widetilde{\Gamma })),%
\mathcal{H}_c^p(\Bbb{R},\mathcal{M})$ is complemented in $L^p(L^\infty (\Bbb{%
R)}\otimes \mathcal{M},L_c^2(\widetilde{\Gamma }))$ for all $1<p<\infty .$
This result is the function space analogue of the non-commutative Stein
inequality in \cite{[PX2]}. This chapter is largely inspired by the recent
work of M. Junge and Q. Xu, where the above results for non-commutative
martingales have been obtained.

In Chapter 5, we further exploit the reduction idea introduced in Chapter 3,
in order to describe $\mathrm{BMO}_c^q(\Bbb{R},\mathcal{M})$ as $\mathrm{BMO}%
_c^{q,\mathcal{D}}(\Bbb{R},\mathcal{M})\cap $ $\mathrm{BMO}_c^{q,\mathcal{D}%
^{\prime }}$ $(\Bbb{R},\mathcal{M}).$ These two latter $\mathrm{BMO}$ spaces
are those of dyadic non-commutative martingales. Among the consequences
given in this chapter, we mention the equivalence of $L^p(L^\infty (\Bbb{R}%
)\otimes \mathcal{M})$ and $\mathcal{H}_{cr}^p(\Bbb{R},\mathcal{M})$ for all 
$1<p<\infty .$

Chapter 6 deals with the interpolation for our Hardy spaces. As expected,
these spaces behave very well with respect to the complex and real
interpolations. This chapter also contains a result on Fourier multipliers.

We close this introduction by mentioning that throughout the paper the
letter $c$ will denote an absolute positive constant, which may vary from
lines to lines, and $c_p$ a positive constant depending only on $p.$

\chapter{Preliminaries}

\section{The non-commutative spaces $L^p(\mathcal{M},L_c^2(\Omega ))$}

\setcounter{theorem}{0} \setcounter{equation}{0} Let $\mathcal{M}$ be a von
Neumann algebra equipped with a normal semifinite faithful trace $\tau .$
Let $S_{\mathcal{M}}^{+}$ be the set of all positive $x\,$in $\mathcal{M}$
such that $\tau ($\textrm{supp} $x)<\infty ,$ where {$\mathrm{supp}$ $x$}
denotes the support of $x,$ that is, the least projection $e\in \mathcal{M}$
such that $ex=x $ (or $xe=x).$ Let $S_{\mathcal{M}}$ be the linear span of $%
S_{\mathcal{M}}^{+}. $ We define 
\[
\left\| x\right\| _p=(\tau |x|^p)^{\frac 1p},\quad \forall x\in S_{\mathcal{M%
}} 
\]
where $|x|=(x^{*}x)^{\frac 12}.$\ One can check that $\left\| \cdot \right\|
_p$ is well-defined and is a norm on $S_{\mathcal{M}}$ if $1\leq p<\infty .$
The completion of $(S_{\mathcal{M}},\left\| \cdot \right\| _p)$ is denoted
by $L^p(\mathcal{M})$ which is the usual non-commutative $L^p$ space
associated with $(\mathcal{M},\tau ).$ For convenience, we usually set $%
L^\infty (\mathcal{M})=\mathcal{M\,}$equipped with the operator norm $%
\left\| \cdot \right\| _{\mathcal{M}}.$ The elements in $L^p(\mathcal{M}%
,\tau )$ can also be viewed as closed densely defined operators on $H$ ($H$
being the Hilbert space on which $\mathcal{M}$ acts). We refer to \cite{F}
for more information on non-commutative $L^p$ spaces.

Let $(\Omega ,\mu )$ be a measurable space. We say $h\,$is a 
$S_{\mathcal{M}}$-valued simple function on $(\Omega ,\mu )$ if it can be
written as 
\begin{equation}
h=\sum_{i=1}^nm_i\cdot \chi _{A_i}  \label{hmn}
\end{equation}
where $m_i\in S_{\mathcal{M}}$ and $A_i$'s are measurable disjoint subsets
of $\Omega \,$with $\mu (A_i)<\infty .$ For such a function $h$ we define 
\[
\left\| h\right\| _{L^p(\mathcal{M},L_c^2(\Omega ))}=\left\| \left(
\sum_{i=1}^nm_i^{*}m_i\cdot \mu (A_i)\right) ^{\frac 12}\right\| _{L^p(%
\mathcal{M})} 
\]
and 
\[
\left\| h\right\| _{L^p(\mathcal{M},L_r^2(\Omega ))}=\left\| \left(
\sum_{i=1}^nm_im_i^{*}\cdot \mu (A_i)\right) ^{\frac 12}\right\| _{L^p(%
\mathcal{M})} 
\]
This gives two norms on the family of all such $h^{\prime }$s. To see that,
denoting by $B(L^2(\Omega ))$ the space of all bounded operators on $%
L^2(\Omega )$ with its usual trace $tr,\,$we consider the von Neumann
algebra tensor product $\mathcal{M}\otimes B(L^2(\Omega ))$ \thinspace with
the product trace $\tau \otimes tr.$ Given a set $A_0\subset \Omega $ with $%
\mu (A_0)=1,$ any element of the family of $h$'s above can be regarded as an
element in $L^p\left( \mathcal{M}\otimes B(L^2(\Omega ))\right) $ via the
following map: 
\begin{equation}
h\mapsto T(h)=\sum_{i=1}^nm_i\otimes (\chi _{A_i}\otimes \chi _{A_0})
\label{th}
\end{equation}
and 
\[
\left\| h\right\| _{L^p(\mathcal{M};L_c^2(\Omega ))}=\left\| T(h)\right\|
_{L^p(\mathcal{M}\otimes B(L^2(\Omega )))} 
\]
Therefore, $\left\| \cdot \right\| _{L^p(\mathcal{M};L_c^2(\Omega ))}\,$%
defines a norm on the family of the $h$'s. The corresponding completion (for 
$1\leq p<\infty $) is a Banach space, denoted by $L^p(\mathcal{M}%
;L_c^2(\Omega )).$ Then $L^p(\mathcal{M};L_c^2(\Omega ))$ is isometric to
the column subspace of $L^p(\mathcal{M}\otimes B(L^2(\Omega ))).$ For $%
p=\infty $ we let $L^\infty (\mathcal{M};L_c^2(\Omega ))$ be the Banach
space isometric by the above map $T$ to the column subspace of $L^\infty (%
\mathcal{M}\otimes B(L^2(\Omega ))).$

Similarly to $\left\| \cdot \right\| _{L^p(\mathcal{M};L_c^2(\Omega ))}$, $%
\left\| \cdot \right\| _{L^p(\mathcal{M};L_r^2(\Omega ))}\,$is also a norm
on the family of $S_{\mathcal{M}}$-valued simple functions and it defines
the Banach space $L^p(\mathcal{M};L_r^2(\Omega ))$ which is isometric to the
row subspace of $L^p(\mathcal{M}\otimes B(L^2(\Omega ))).$

Alternatively, we can fix an orthonormal basis of $L^2(\Omega )$. Then any
element of $L^p(\mathcal{M}\otimes B(L^2(\Omega )))$ can be identified with
an infinite matrix with entries in $L^p(\mathcal{M})$. Accordingly, $L^p(%
\mathcal{M};L_c^2(\Omega ))$ (resp. $L^p(\mathcal{M};L_r^2(\Omega ))$) can
be identified with the subspace of $L^p(\mathcal{M}\otimes B(L^2(\Omega )))$
consisting of matrices whose entries are all zero except those in the first
column (resp. row).

\begin{proposition}
Let $f\in L^p(\mathcal{M};L_c^2(\Omega )),g\in L^q(\mathcal{M};L_c^2(\Omega
))(1\leq p,q\leq \infty ),\frac 1r=\frac 1p+\frac 1q.$ Then $\langle
g,f\rangle $ exists as an element in $L^r(\mathcal{M})$ and 
\[
\left\| \langle g,f\rangle \right\| _{L^r(\mathcal{M})}\leq \left\|
g\right\| _{L^q(\mathcal{M};L_c^2(\Omega ))}\left\| f\right\| _{L^p(\mathcal{%
M};L_c^2(\Omega ))}, 
\]
where $\langle \,,\,\rangle $ denotes the scalar product in $L_c^2(\Omega )$%
. A similar statement also holds for row spaces.
\end{proposition}

\noindent\textbf{Proof.} This is clear from the discussion above via the
matrix representation of $L^p(\mathcal{M};L_c^2(\Omega ))$ (in an
orthonormal basis of $L^2(\Omega )$).\qed\medskip

\noindent\textbf{Remark.} Note that if $f$ and $g$ are $S_{\mathcal{M}}$%
-valued simple functions, then 
\[
\langle g ,f\rangle=\int_\Omega g^{*}fd\mu . 
\]
For general $f$ and $g$ as in Proposition 1.1, if one of $p$ and $q$ is
finite, one can easily prove that $\langle g ,f\rangle$ is the limit in $L^r(%
\mathcal{M})$ of a sequence $(\langle g_n ,f_n\rangle)_n$ with $S_{\mathcal{M%
}}$-valued simple functions $f_n,\ g_n$. Consequently, we can define $%
\int_\Omega g^{*}fd\mu$ as the limit of $\int_\Omega g_n^{*}f_nd\mu$. If
both $p$ and $q$ are infinite, this limit procedure is still valid but only
in the w*-sense.

\noindent\textbf{Convention.} Throughout this paper whenever we are in the
situation of Proposition 1.1, we will write $\langle g,f\rangle $ as the
integral $\int_\Omega g^{*}fd\mu $. Notationally, this is clearer. Moreover,
by the proceding remark this indeed makes sense in many cases.\medskip

Observe that the column and row subspaces of $L^p(\mathcal{M}\otimes
B(L^2(\Omega )))\,$ are 1-complemented subspaces. Therefore, from the
classical duality between $L^p(\mathcal{M}\otimes B(L^2(\Omega )))$ and $L^q(%
\mathcal{M}\otimes B(L^2(\Omega )))$ $(\frac 1p+\frac 1q=1,1\leq p<\infty )$
we deduce that 
\[
\left( L^p(\mathcal{M};L_c^2(\Omega ))\right) ^{*}=L^q(\mathcal{M}%
;L_c^2(\Omega )) 
\]
and 
\[
\left( L^p(\mathcal{M};L_r^2(\Omega ))\right) ^{*}=L^q(\mathcal{M}%
;L_r^2(\Omega )) 
\]
isometrically via the antiduality 
\[
(f,g)\mapsto \tau (\langle g,f\rangle )=\tau \int_\Omega g^{*}fd\mu . 
\]
Moreover, it is well known that (by the same reason), for $0<\theta <1$ and $%
1\leq p_0,p_1,p_{\theta} \leq \infty $ with $\frac 1{p_{\theta} }=\frac{%
1-\theta }{p_0}+\frac \theta {p_1},$ we have isometrically 
\begin{equation}
\left( L^{p_0}(\mathcal{M};L_c^2(\Omega )),L^{p_1}(\mathcal{M};L_c^2(\Omega
))\right) _{\theta} =L^{p_{\theta} }(\mathcal{M};L_c^2(\Omega )).
\label{interp}
\end{equation}

In the following, we are mainly interested in the spaces $L^p(\mathcal{M}%
;L_c^2(\Omega ))$ (resp. $L^p(\mathcal{M};L_r^2(\Omega )))$ with $(\Omega
,\mu )=\widetilde{\Gamma }=(\Gamma ,dxdy)\times (\{1,2\},\sigma ),\,$where $%
\Gamma =\{(x,y)\in \Bbb{R}_{+}^2,|x|<y\},$ $\sigma \{1\}=\sigma \{2\}=1.($
This cone $\Gamma $ is a fundamental subject used in the classical harmonic
analysis, see \cite{G}, \cite{fs}, \cite{Ko}, \cite{[St]} or any book on
Hardy spaces$).$ The presence of $\{1,2\}$ corresponds to our two variables $%
x,y,$ see below. We then denote them by $L^p(\mathcal{M},L_c^2(\widetilde{%
\Gamma }))$ $($resp. $L^p(\mathcal{M},L_r^2(\widetilde{\Gamma }))).$ For
simplicity, we will abbreviate them as $L^p(\mathcal{M},L_c^2)$ (resp. $L^p(%
\mathcal{M},L_r^2)$) if no confusion can arise.$\,$

\section{Operator valued Hardy spaces}

Let $1\leq p<\infty .$ For any $S_{\mathcal{M}}$-valued simple function $f$
on $\Bbb{R},$ we also use $f$ to denote its Poisson integral on the upper
half plane $\Bbb{R}_{+}^2=\{(x,y)|y>0\},$ 
\[
f(x,y)=\int_{\Bbb{R}}P_y(x-s)f(s)ds,\ \ (x,y)\in \Bbb{R}_{+}^2\ , 
\]
where $P_y(x)$ is the Poisson kernel (i.e. $P_y(x)=\frac 1\pi \frac
y{x^2+y^2}$). Note that $f(x,y)$ is a harmonic function still with values in 
$S_{\mathcal{M}},$ and so in $\mathcal{M}.$ Define the$\,\mathcal{H}_c^p(%
\Bbb{R},\mathcal{M})\,$ norm of $f$ by 
\[
\left\| f\right\| _{\mathcal{H}_c^p}=\left\| \nabla f(x+t,y)\chi _\Gamma
(x,y)\right\| _{L^p(L^\infty (\Bbb{R},dt\Bbb{)}\otimes \mathcal{M},L_c^2(%
\widetilde{\Gamma }))}, 
\]
where $\nabla f$ is the gradient of the Poisson integral $f(x,y)$ and $%
\widetilde{\Gamma }$ is defined as in the end of Section 1.1. In this paper,
we will always regard $\nabla f(x+t,y)\chi _\Gamma (x,y)$ as functions
defined on $\Bbb{R}\times \widetilde{\Gamma }$ with $t\in \Bbb{R},(x,y)\in
\Gamma $ and 
\[
\nabla f(x+t,y)(1)=\frac{\partial f}{\partial x}(x+t,y),\ \ \nabla
f(x+t,y)(2)=\frac{\partial f}{\partial y}(x+t,y). 
\]
And set 
\[
|\nabla f(x+t,y)|^2=|\frac{\partial f}{\partial x}(x+t,y)|^2+|\frac{\partial
f}{\partial y}(x+t,y)|^2. 
\]
Define the $\mathcal{H}_r^p(\Bbb{R},\mathcal{M})$ norm of $f$ by 
\[
\left\| f\right\| _{\mathcal{H}_r^p}=\left\| \nabla f(x+t,y)\chi _\Gamma
\right\| _{L^p(L^\infty (\Bbb{R)}\otimes \mathcal{M},L_r^2)}. 
\]
Set $\mathcal{H}_c^p(\Bbb{R},\mathcal{M})\,$(resp. $\mathcal{H}_r^p(\Bbb{R},%
\mathcal{M})$) to be the completion of the space of all $S_{\mathcal{M}}$%
-valued simple function $f$'s with finite $\mathcal{H}_c^p(\Bbb{R},
\mathcal{M})$(resp. $\mathcal{H}_r^p(\Bbb{R},\mathcal{M})$) norm. Equipped
respectively with the previous norms, $\mathcal{H}_c^p(\Bbb{R},\mathcal{M})$
and $\mathcal{H}_r^p(\Bbb{R},\mathcal{M})$ are Banach spaces. Define the
non-commutative analogues of the classical Lusin integral by 
\begin{eqnarray}
S_c(f)(t) &=&(\mathop{\displaystyle\int\!\!\!\int}_\Gamma |\nabla
f(x+t,y)|^2dxdy)^{\frac 12}  \label{Sc} \\
S_r(f)(t) &=&(\mathop{\displaystyle\int\!\!\!\int}_\Gamma |\nabla
f^{*}(x+t,y)|^2dxdy)^{\frac 12}.  \label{Sr}
\end{eqnarray}
Note that 
\[
|\nabla f(x,y)|^2=\int_{\{1,2\}}|\nabla f(x,y)(i)|^2d\sigma (i). 
\]
Then, for $f\in \mathcal{H}_c^p(\Bbb{R},\mathcal{M}),$%
\[
\left\| f\right\| _{\mathcal{H}_c^p}=\left\| S_c(f)\right\| _{L^p(L^\infty (%
\Bbb{R)}\otimes \mathcal{M})} 
\]
and the similar equality holds for $\mathcal{H}_r^p(\Bbb{R},\mathcal{M})$. $%
S_c(f)$ and $S_r(f)$ are the non-commutative analogues of the classical
Lusin square function. We will need the non-commutative analogues of the
classical Littlewood-Paley $g$-function, which are defined by 
\begin{eqnarray}
G_c(f)(t) &=&(\int_{\Bbb{R}_{+}}|\nabla f(t,y)|^2ydy)^{\frac 12}\ 
\label{Gc} \\
\ G_r(f)(t) &=&(\int_{\Bbb{R}_{+}}|\nabla f^{*}(t,y)|^2ydy)^{\frac 12}
\label{Gr}
\end{eqnarray}
We will see, in Chapters 2 and 4, that 
\begin{eqnarray*}
\left\| S_c(f)\right\| _{L^p(L^\infty (\Bbb{R)}\otimes \mathcal{M})}
&\backsimeq &\left\| G_c(f)\right\| _{L^p(L^\infty (\Bbb{R)}\otimes \mathcal{%
M})} \\
\left\| S_r(f)\right\| _{L^p(L^\infty (\Bbb{R)}\otimes \mathcal{M})}
&\backsimeq &\left\| G_r(f)\right\| _{L^p(L^\infty (\Bbb{R)}\otimes \mathcal{%
M})}
\end{eqnarray*}
for all $1\leq p<\infty .$

Define the Hardy spaces of non-commutative functions $f$ as follows: if $%
1\leq p<2,$%
\begin{equation}
\mathcal{H}_{cr}^p(\Bbb{R},\mathcal{M})=\mathcal{H}_c^p(\Bbb{R},\mathcal{M})+%
\mathcal{H}_r^p(\Bbb{R},\mathcal{M})  \label{Hcr}
\end{equation}
equipped with the norm 
\[
\left\| f\right\| _{\mathcal{H}_{cr}^p}=\inf \{\left\| g\right\| _{\mathcal{H%
}_c^p}+\left\| h\right\| _{\mathcal{H}_r^p}:f=g+h,g\in \mathcal{H}_c^p(\Bbb{R%
},\mathcal{M}),h\in \mathcal{H}_r^p(\Bbb{R},\mathcal{M})\} 
\]
and if $2\leq p<\infty ,$%
\begin{equation}
\mathcal{H}_{cr}^p(\Bbb{R},\mathcal{M})=\mathcal{H}_c^p(\Bbb{R},\mathcal{M}%
)\cap \mathcal{H}_r^p(\Bbb{R},\mathcal{M})  \label{Hcr2}
\end{equation}
equipped with the norm 
\[
\left\| f\right\| _{\mathcal{H}_{cr}^p}=\max \{\left\| f\right\| _{\mathcal{H%
}_c^p},\left\| f\right\| _{\mathcal{H}_r^p}\}. 
\]
\medskip

\noindent\textbf{Remark. }We have 
\[
\mathcal{H}_c^2(\Bbb{R},\mathcal{M})=\mathcal{H}_r^2(\Bbb{R},\mathcal{M})=%
\mathcal{H}_{cr}^2(\Bbb{R},\mathcal{M})=L^2(L^\infty (\Bbb{R)}\otimes 
\mathcal{M}). 
\]
In fact, notice that $\bigtriangleup |f|^2=2|\nabla f|^2$ and $%
f(x,y)(|x|+y)\rightarrow 0,\nabla f(x,y)(|x|+y)^2\rightarrow 0$ as $%
|x|+y\rightarrow 0,$ for $S_{\mathcal{M}}$-valued simple function $f$'s. By
the Green theorem $\,$%
\begin{eqnarray}
&&||\nabla f(t+x,y)\chi _\Gamma ||_{L^2(L^\infty (\Bbb{R)}\otimes \mathcal{M}%
,L_c^2)}^2  \nonumber \\
&=&2\tau \mathop{\displaystyle\int\!\!\!\int}_{\Bbb{R}_{+}^2}|\nabla
f|^2ydxdy  \nonumber \\
&=&\tau \mathop{\displaystyle\int\!\!\!\int}_{\Bbb{R}_{+}^2}\bigtriangleup
|f|^2ydxdy  \nonumber \\
&=&\tau \int_{\Bbb{R}}|f|^2ds=\left\| f\right\| _{L^2(L^\infty (\Bbb{R)}%
\otimes \mathcal{M}).}^2  \label{h2l2}
\end{eqnarray}
Similarly, $||f||_{\mathcal{H}_r^2}=\left\| f^{*}\right\| _{L^2(L^\infty (%
\Bbb{R)}\otimes \mathcal{M})}=\left\| f\right\| _{L^2(L^\infty (\Bbb{R)}%
\otimes \mathcal{M}).}$

Note we have also the following polarized version of (\ref{h2l2}), 
\begin{equation}
2\mathop{\displaystyle\int\!\!\!\int}_{\Bbb{R}_{+}^2}\nabla f(x,y)\nabla
g(x,y)ydxdy=\int_{\Bbb{R}}f(s)g(s)ds  \label{greenf}
\end{equation}
for $S_{\mathcal{M}}$-valued simple function $f,g$'s.

We will repeatedly use the following consequence of the convexity of the
operator valued function: $x\mapsto |x|^2$ (This convexity follows from the
convexity of $x\mapsto \langle x^*xh,h\rangle=\|xh\|^2$ for any $h$).
Let $f:(\Omega ,\mu )\rightarrow \mathcal{M}$ be a weak-* integrable
function, we have 
\begin{equation}
|\int_Af(t)d\mu (t)|^2\leq \mu (A)\int_A|f(t)|^2d\mu (t),\quad \forall
A\subset \Omega  \label{fu}
\end{equation}
In particular, set $d\mu (t)=g^2(t)dt,$ 
\begin{equation}
|\int_Af(t)g(t)dt|^2\leq \int_A|f(t)|^2dt\int_Ag^2(t)dt,\quad \forall
A\subset \Bbb{R}  \label{fg1}
\end{equation}
for every measurable function $g$ on $\Bbb{R},$ and 
\begin{equation}
|\int_Af(t)dt|^2\leq \int_A|f(t)|^2g^{-1}(t)dt\int_Ag(t)dt,\quad \forall
A\subset \Bbb{R}  \label{fg}
\end{equation}
for every positive measurable function $g$ on $\Bbb{R}.$

Let $H^p(\Bbb{R})$ ($1\le p<\infty $) denote the classical Hardy space on $%
\Bbb{R}.$ It is well known that 
\[
H^p(\Bbb{R})=\{f\in L^p(\Bbb{R}):S(f)\in L^p(\Bbb{R})\}, 
\]
where $S(f)$ is the classical Lusin integral function ($S(f)$ is equal to $%
S_c(f)$ above by taking $\mathcal{M}=\Bbb{C}$). In the following, $H^p(\Bbb{R%
})$ is always equipped with the norm $\left\| S(f)\right\| _{L^p(\Bbb{R})}.$

\begin{proposition}
Let $1\le p<\infty $, $f\in \mathcal{H}_c^p(\Bbb{R},\mathcal{M})$ and $m\in
L^q(\mathcal{M})$ (with $q$ the index conjugate to $p$). Then $\tau (mf)\in
H^p(\Bbb{R})$ and 
\[
\left\| \tau (mf)\right\| _{H^p}\leq \left\| m\right\| _{L^q(\mathcal{M}%
)}\left\| f\right\| _{\mathcal{H}_c^p}. 
\]
\end{proposition}

\noindent\textbf{Proof.} Note that 
\[
\nabla (\tau (mf)*P)=\tau (m(f*\nabla P))=\tau (m\nabla f), 
\]
here $P$ is the Poisson kernel (i.e. $P_y(x)=\frac 1\pi \frac y{x^2+y^2}$).
By (\ref{fg1}), we have 
\begin{eqnarray*}
&&\left\| \tau (mf)\right\| _{H^p}^p \\
&=&\int_{\Bbb{R}}(\mathop{\displaystyle\int\!\!\!\int}_\Gamma |\tau (m\nabla
f(x+t,y))|^2dxdy)^{\frac p2}dt \\
&\leq &\int_{\Bbb{R}}\sup_{\| g\| _{L^2(\widetilde{\Gamma })}\leq
1}\left| \mathop{\displaystyle\int\!\!\!\int}_\Gamma g\tau (m\nabla
f(x+t,y))dxdy\right| ^pdt \\
&=&\int_{\Bbb{R}}\sup_{\| g\| _{L^2(\widetilde{\Gamma })}\leq
1}\left| \tau \left[ m\mathop{\displaystyle\int\!\!\!\int}_\Gamma g_1%
\frac{\partial f}{\partial x}(x+t,y)+g_2\frac{\partial f}{\partial y}%
(x+t,y)dxdy\right] \right| ^p\!\!\!dt \\
&\leq &\int_{\Bbb{R}}\sup_{\| g\| _{L^2(\widetilde{\Gamma })}\leq
1}\Vert m\Vert _{L^q(\mathcal{M})}^p\left\| 
\mathop{\displaystyle\int\!\!\!\!
\int}_\Gamma g_1\frac{\partial f}{\partial x}(x+t,y)+g_2\frac{\partial f}{%
\partial y}(x+t,y)dxdy\right\| _{L^p(\mathcal{M})}^p\!\!\!dt \\
&\leq &\!\Vert m\Vert _{L^q(\mathcal{M})}^p\int_{\Bbb{R}}\sup_{\left\|
g\right\| _{L^2(\widetilde{\Gamma })}\leq 1}\left\| (\mathop{\displaystyle%
\int\!\!\!\!\int}_\Gamma |g|^2dxdy)^{\frac 12}(\mathop{\displaystyle\int \!\!\!\!
\int}_\Gamma |\nabla f(x+t,y)|^2dxdy)^{\frac 12}\right\| _{L^p(\mathcal{M}%
)}^p\! dt \\
&\leq&\!\Vert m\Vert _{L^q(\mathcal{M})}^p\tau \int_{\Bbb{R}}(%
\mathop{\displaystyle\int\!\!\!\!\int}_\Gamma |\nabla f(x+t,y)|^2dxdy)^{\frac
p2}\!dt \\
&=&\!\Vert m\Vert _{L^q(\mathcal{M})}^p\ \left\| f\right\| _{\mathcal{H}%
_c^p}^p.\qed
\end{eqnarray*}

\medskip

\noindent\textbf{Remark. }We should emphasize that for two functions $g,f$
defined on $\widetilde{\Gamma },$ we always set 
\[
gf(z)=g(z)(1)f(z)(1)+g(z)(2)f(z)(2). 
\]
Then in the above formula $|\tau (m\nabla f(x+t,y))|^2$ and $g\tau (m\nabla
f(x+t,y))$ etc. are functions defined on $\Gamma .$ We will use very often
such a product for ($\mathcal{M}$-valued) functions defined on $\widetilde{%
\Gamma }$.\medskip

\noindent\textbf{Remark.} (i) $\int fdt=0,\forall f\in \mathcal{H}_c^1(\Bbb{R%
},\mathcal{M}).$ In fact, if $f\in \mathcal{H}_c^1(\Bbb{R},\mathcal{M}),$ by
Proposition 1.2 and the classical property of $H^1($see \cite{[St]}, p.128$%
), $ we have $\tau (m\int fdt)=0,\forall m\in \mathcal{M}.$ Thus $\int
fdt=0. $

(ii) The collection of all $S_{\mathcal{M}}$-valued simple functions $f$
such that $\int fdt=0$ is a dense subset of $\mathcal{H}_c^p(\Bbb{R},%
\mathcal{M})(1<p<\infty ).$ Note that 
\[
\lim_{N\rightarrow \infty }\left\| \frac mN\chi _{[-N,N]}(t)\right\| _{%
\mathcal{H}_c^p(\Bbb{R},\mathcal{M})}=0,\quad \forall m\in S_{\mathcal{M}}. 
\]
For a simple function $f,$ let $f_N=$ $f-\frac{\int fdt}N\chi _{[-N,N]}.$
Then $\int f_N=0$ and $f_N\rightarrow f$ in $\mathcal{H}_c^p(\Bbb{R},%
\mathcal{M}).$ \medskip

\noindent\textbf{Remark.} See \cite{fs} and \cite{[St]} for the discussions
on the classical Lusin integral and the Littlewood-Paley $g$-function and
the fact that {a scalar valued function is in }$H^1$ {if and only if its
Lusin integral is in }$L^1.$ We define the non-commutative Hardy spaces $%
\mathcal{H}_{cr}^p(\Bbb{R},\mathcal{M})$ differently for the case $1\leq p<2$
and $p\geq 2$ (respectively by (\ref{Hcr}) and (\ref{Hcr2})) as Pisier and
Xu did for non-commutative martingales in \cite{JX}. This is to get the
expected equivalence between $\mathcal{H}_{cr}^p(\Bbb{R},\mathcal{M})$ and $%
L^p(\Bbb{R},\mathcal{M})$ for $1<p<\infty$ (see Chapter 5). And $\mathcal{H%
}_c^p(\Bbb{R},\mathcal{M})$ or $\mathcal{H}_r^p(\Bbb{R},\mathcal{M})$ alone
could be very far away from $L^p(\Bbb{R},\mathcal{M})$ for $p\neq 2$.

\section{Operator valued \textrm{BMO} spaces}

Now, we introduce the non-commutative analogue of \textrm{BMO }spaces. For
any interval $I$ on $\Bbb{R},$ we will denote its center by $C_I$ and its
Lebesgue measure by $|I|.$ Let ${\varphi }\in L^\infty (\mathcal{M},L_c^2(%
\Bbb{R},\frac{dt}{1+t^2})).$ By Proposition 1.1 (and our convention), for
every $g\in L^2(\Bbb{R},\frac{dt}{1+t^2}),\int_{\Bbb{R}}g{\varphi }\frac{dt}{%
1+t^2}\in \mathcal{M}.$ Then the mean value of ${\varphi }$ over $I$ ${%
\varphi }_I:=\frac 1{|I|}\int_I{\varphi }(s)ds$ exists as an element in $%
\mathcal{M}.$ And the Poisson integral of ${\varphi }$ 
\[
{\varphi }(x,y)=\int_{\Bbb{R}}P_y(x-s){\varphi }(s)ds 
\]
also exists as an element in $\mathcal{M}.$ $\,$Set 
\begin{equation}
\left\| {\varphi }\right\| _{\mathrm{BMO}_c}=\sup_{I\subset \Bbb{R}}\left\{
\left\| \left( \frac 1{|I|}\int_I|{\varphi }-{\varphi }_I|^2d\mu \right)
^{\frac 12}\right\| _{\mathcal{M}}\right\}  \label{dbmo}
\end{equation}
where again $|{\varphi }-{\varphi }_I|^2=({\varphi }-{\varphi }_I)^{*}({%
\varphi }-{\varphi }_I)$ and the supremum runs over all intervals $I\subset 
\Bbb{R}.($see $\,$Let $H$ be the Hilbert space on which $\mathcal{M\,}$acts.
Obviously, we have 
\begin{equation}
\left\| {\varphi }\right\| _{\mathrm{BMO}_c}=\sup_{e\in H,\left\| e\right\|
=1}\left\| {\varphi }e\right\| _{\mathrm{BMO}_2(\Bbb{R},H)}  \label{bmoh}
\end{equation}
where $\mathrm{BMO}_2(\Bbb{R},H)$ is the usual $H$-valued \textrm{BMO} space
on $\Bbb{R}.$ Thus $\left\| \cdot \right\| _{\mathrm{BMO}_c}\,$is a norm
modulo constant functions. Set $\mathrm{BMO}_c(\Bbb{R},\mathcal{M})\,$to be
the space of all ${\varphi }\in L^\infty (\mathcal{M},L_c^2(\Bbb{R},\frac{dt%
}{1+t^2}))$ such that $\left\| {\varphi }\right\| _{\mathrm{BMO}_c}<\infty .$
$\mathrm{BMO}_r(\Bbb{R},\mathcal{M})$ is defined as the space of all ${%
\varphi }$'s such that ${\varphi }^{*}\in \mathrm{BMO}_c(\Bbb{R},\mathcal{M}%
) $ with the norm $\left\| {\varphi }\right\| _{\mathrm{BMO}_r}=\left\| {%
\varphi }^{*}\right\| _{\mathrm{BMO}_c}.$ We define $\mathrm{BMO}_{cr}(\Bbb{R%
},\mathcal{M})$ as the intersection of these two spaces 
\[
\mathrm{BMO}_{cr}(\Bbb{R},\mathcal{M})=\mathrm{BMO}_c(\Bbb{R},\mathcal{M}%
)\,\cap \mathrm{BMO}_r(\Bbb{R},\mathcal{M}) 
\]
with the norm 
\[
\left\| {\varphi }\right\| _{\mathrm{BMO}_{cr}}=\max \{\left\| {\varphi }%
\right\| _{\mathrm{BMO}_c},\left\| {\varphi }\right\| _{\mathrm{BMO}_r}\}. 
\]
As usual, the constant functions are considered as zero in these BMO spaces,
and then these spaces are normed spaces (modulo constants). \medskip

Given an interval $I,$ we denote by $2^kI$ the interval $%
\{t:|t-C_I|<2^{k-1}|I|\}.$ The technique used in the proof of the following
Proposition is classical(see \cite{[St]}).

\begin{proposition}
Let ${\varphi }\in \mathrm{BMO}_c(\Bbb{R},\mathcal{M}).\,$Then 
\[
\left\| {\varphi }\right\| _{L^\infty (\mathcal{M},L_c^2(\Bbb{R},\frac{dt}{%
1+t^2}))}\leq c(\left\| {\varphi }\right\| _{\mathrm{BMO}c}+\left\| {\varphi 
}_{I_1}\right\| _{\mathcal{M}}) 
\]
where $I_1=(-1,1].$ Moreover, $\mathrm{BMO}_c(\Bbb{R},\mathcal{M}),\mathrm{%
BMO}_r(\Bbb{R},\mathcal{M}),\mathrm{BMO}_{cr}(\Bbb{R},\mathcal{M})$ are
Banach spaces.
\end{proposition}

\noindent\textbf{Proof.} Let $\,{\varphi }\in \mathrm{BMO}_c(\Bbb{R},%
\mathcal{M})$ and $I$ be an interval. Using (\ref{fu}), (\ref{fg}) we have 
\begin{eqnarray}
|{\varphi }_{2^nI}-{\varphi }_I|^2 &\leq &n\sum_{k=0}^{n-1}|{\varphi }%
_{2^kI}-{\varphi }_{2^{k+1}I}|^2  \nonumber \\
&=&n\sum_{k=0}^{n-1}\left| \frac 1{|2^kI|}\int_{2^kI}({\varphi }(s)-{\varphi 
}_{_{2^{k+1}I}})ds\right| ^2  \nonumber \\
&\leq &n\sum_{k=0}^{n-1}\frac 2{|2^{k+1}I|}\int_{2^{k+1}I}|{\varphi }(s)-{%
\varphi }_{_{2^{k+1}I}}|^2ds  \nonumber \\
&\leq &2n\left\| {\varphi }\right\| _{\mathrm{BMO}_c}^2.  \label{fai}
\end{eqnarray}
$\,$By (\ref{fg}), (\ref{fai}), 
\begin{eqnarray}
&&\left\| \int_{\Bbb{R}}\frac{|{\varphi }(t)|^2}{1+t^2}dt\right\| _{\mathcal{%
M}}  \nonumber \\
&=&\left\| \int_{I_1}\frac{|{\varphi }(t)|^2}{1+t^2}dt+\sum_{k=0}^\infty
\int_{2^{k+1}I_1/2^kI_1}\frac{|{\varphi }(t)|^2}{1+t^2}dt\right\| _{\mathcal{%
M}}  \nonumber \\
&\leq &2\left\| \int_{I_1}(|{\varphi }(t)-{\varphi }_{_{I_1}}|^2+|{\varphi }%
_{I_1}|^2)dt\right\| _{\mathcal{M}}  \nonumber \\
&&+4\left\| \sum_{k=0}^\infty \int_{2^{k+1}I_1/2^kI_1}\frac{|{\varphi }(t)-{%
\varphi }_{_{2^{k+1}I_1}}|^2+|{\varphi }_{2^{k+1}I_1}-{\varphi }_{I_1}|^2+|{%
\varphi }_{I_1}|^2}{2^{2k}}dt\right\| _{\mathcal{M}}  \nonumber \\
&\leq &c(\left\| |{\varphi }_{I_1}|^2\right\| _{\mathcal{M}}+\left\| {%
\varphi }\right\| _{\mathrm{BMO}_c}^2)  \label{1+t}
\end{eqnarray}
Thus 
\[
\left\| {\varphi }\right\| _{_{L^\infty (\mathcal{M},L_c^2(\Bbb{R},\frac{dt}{%
1+t^2}))}}=\left\| (\int_{\Bbb{R}}\frac{|{\varphi }(t)|^2}{1+t^2}dt)^{\frac
12}\right\| _{\mathcal{M}}\leq c(\left\| {\varphi }_{I_1}\right\| _{\mathcal{%
M}}+\left\| {\varphi }\right\| _{\mathrm{BMO}_c}) 
\]
And then $\mathrm{BMO}_c(\Bbb{R},\mathcal{M})$ is complete. Consequently, $%
\mathrm{BMO}_c(\Bbb{R},\mathcal{M}),$ $\mathrm{BMO}_r(\Bbb{R},\mathcal{M}),$ 
$\mathrm{BMO}_{cr}(\Bbb{R},\mathcal{M})$ are Banach spaces. \qed\medskip

It is classical that $\mathrm{BMO}$ functions are related with Carleson
measures(See \cite{G}, \cite{Ko}). The same relation still holds in the
present non-commutative setting. We say that an $\mathcal{M}$-valued measure 
$d{\lambda }\,$on $\Bbb{R}_{+}^2$ is a Carleson measure if

\[
N({\lambda })=\mathrm{sup}_I\left\{ \frac 1{|I|}\left\| \mathop{\displaystyle%
}\mathop{\displaystyle\int\!\!\!\int}_{T(I)}d{\lambda }\right\| _{\mathcal{M}%
}:I\in \Bbb{R}\mbox{ interval}\right\} <\infty , 
\]
where, as usual, $T(I)=I\times (0,|I|].$

\begin{lemma}
Let\textbf{\ }${\varphi }\in \mathrm{BMO}_c(\Bbb{R},\mathcal{M}).$ Then $d{%
\lambda }_\varphi =|\nabla {\varphi }|^2ydxdy\,$is an $\mathcal{M}$-valued
Carleson measure on $\Bbb{R}_{+}^2$ and $N({\lambda }_\varphi )\leq c\left\| 
{\varphi }\right\| _{\mathrm{BMO}_c}^2.$
\end{lemma}

\noindent\textbf{Proof. }The proof is very similar to the scalar situation
(see \cite{[St]}, p.160). For any interval $I$ on $\Bbb{R},$ write ${\varphi}
= $ ${\varphi} _1+$ ${\varphi} _2+{\varphi} _3,$ where ${\varphi} _1=({%
\varphi} -{\varphi} _{2I})\chi _{2I},{\varphi} _2=({\varphi} -{\varphi}
_{2I})\chi _{_{(2I)^c}}$ and ${\varphi} _3={\varphi} _{2I}.$ Set 
\[
d{\lambda} _{{\varphi} _1}=|\nabla {\varphi} _1|^2ydxdy,d{\lambda} _{{\varphi%
} _{_2}}=|\nabla {\varphi} _2|^2ydxdy. 
\]
Thus 
\[
N({\lambda} _{\varphi} )\leq 2(N({\lambda} _{{\varphi} _1})+N({\lambda} _{{%
\varphi} _2})). 
\]
We treat $N({\lambda} _{{\varphi} _1})$ first. Notice that $\bigtriangleup |{%
\varphi} _1|^2=2|\nabla {\varphi} _1|^2$ and ${\varphi} _1(x,y)(|x|+y)%
\rightarrow 0,\nabla {\varphi} _1(x,y)(|x|+y)^2\rightarrow 0$ as $%
|x|+y\rightarrow 0.$ By the Green theorem $\,$%
\begin{eqnarray}
\frac 1{|I|}\left\| \mathop{\displaystyle\int\!\!\!\int}_{T(I)}|\nabla {%
\varphi} _1|^2ydxdy\right\| _{\mathcal{M}} &\leq &\frac 1{|I|}\left\| %
\mathop{\displaystyle\int\!\!\!\int}_{\Bbb{R}_2^{+}}|\nabla {\varphi}
_1|^2ydxdy\right\| _{\mathcal{M}}  \label{gree} \\
&=&\frac 1{2|I|}\left\| \int_{\Bbb{R}}|{\varphi} _1|^2ds\right\| _{\mathcal{M%
}}  \nonumber \\
&=&\frac 1{2|I|}\left\| \int_{2I}|{\varphi} -{\varphi} _{2I}|^2ds\right\| _{%
\mathcal{M}}\leq \left\| {\varphi} \right\| _{\mathrm{BMO}_c}^2  \nonumber
\end{eqnarray}
To estimate $N({\lambda} _{{\varphi} _1})$, we note 
\[
|\nabla P_y(x-s)|^2\leq \frac 1{4(x-s)^4}\leq \frac 1{4|I|^42^{4k}},\quad
\forall s\in 2^{k+1}I/2^kI,\quad (x,y)\in T(I), 
\]
by (\ref{fg}) and (\ref{fai}) 
\begin{eqnarray*}
&&\frac 1{|I|}\left\| \mathop{\displaystyle\int\!\!\!\int}_{T(I)}|\nabla {%
\varphi} _2|^2ydxdy\right\| _{\mathcal{M}} \\
&=&\frac 1{|I|}\left\| \mathop{\displaystyle\int\!\!\!\int}_{T(I)}|\nabla
\int_{-\infty }^{+\infty }P_y(x-s){\varphi} _2(s)ds|^2ydxdy\right\| _{%
\mathcal{M}} \\
&\leq &\frac 1{|I|}\mathop{\displaystyle\int\!\!\!\int}_{T(I)}\sum_{k=1}^%
\infty \int\limits_{2^{k+1}I/2^kI}|\nabla
P_y(x-s)|^22^{2k}ds\sum_{k=1}^\infty \frac 1{2^{2k}}\left\|
\int\limits_{2^{k+1}I}|{\varphi} _2|^2ds\right\| _{\mathcal{M}}ydxdy \\
&\leq &\frac c{|I|}\mathop{\displaystyle\int\!\!\!\int}_{T(I)}\frac
1{|I|^2}\left\| {\varphi} \right\| _{\mathrm{BMO}_c}^2ydxdy\leq c\left\| {%
\varphi} \right\| _{\mathrm{BMO}_c}^2
\end{eqnarray*}
Therefore $N({\lambda} _{{\varphi} _i})\leq c\left\| {\varphi} \right\| _{%
\mathrm{BMO}_c}^2,i=1,2,$ and then $N({\lambda} _{\varphi} )\leq c\left\| {%
\varphi} \right\| _{\mathrm{BMO}_c}^2$.\qed \medskip

\noindent\textbf{Remark.} We will see later (Corollary 2.6) that the
converse to lemma 1.4 is also true. \thinspace

We will need the following elementary fact to make our later applications of
Green's theorem rigorous in Chapters 2 and 4.

\begin{lemma}
Suppose ${\varphi }\in \mathrm{BMO}_c(\Bbb{R},\mathcal{M})$ and suppose $I$
is an interval such that ${\varphi }_I=0.$ Let $3I$ be the interval
concentric with $I$ having length $3|I|.$ Then there is $\psi \in \mathrm{BMO%
}_c(\Bbb{R},\mathcal{M})$ such that $\psi ={\varphi }$ on $I,\psi =0$ on $%
\Bbb{R}\backslash 3I$ and 
\[
\left\| \psi \right\| _{\mathrm{BMO}_c}\leq c\left\| {\varphi }\right\| _{%
\mathrm{BMO}_c}. 
\]
\end{lemma}

\noindent\textbf{Proof.} This is well known for the classical $\mathrm{BMO}$
and a proof is outlined in \cite{G}, p. 269. One can check that the method
to construct $\psi $ mentioned there works as well for $\mathrm{BMO}_c(\Bbb{R%
},\mathcal{M}).$\qed
\medskip

\noindent\textbf{Remark. }We have seen that the non-commutative $\mathrm{BMO}_c(%
\Bbb{R},\mathcal{M})$ are well adapted to many generalizations of
classical results, such as Proposition 1.3 and Lemma 1.4, 1.5. We will also
prove an analogue of the classical Fefferman duality theorem between $%
\mathcal{H}^1$ and $\mathrm{BMO}$ in the next chapter. However, unlike the
classical case, we could not replace the power 2 by $p$ in the definition of
the non-commutative BMO norm ((\ref{dbmo})). In fact, $\sup_{I\subset \Bbb{R}%
}\left\| \left( \frac 1{|I|}\int_I|{\varphi }-{\varphi }_I|^pd\mu \right)
^{\frac 1p}\right\| _{\mathcal{M}}$ may not be a norm for $p\neq 2$ in the
non-commutative case (Note we do not have $|x_1+x_2|\leq |x_1|+|x_2|$ in
general for $x_1,x_2\in \mathcal{M}$). See the remark at the end of Chapter
6 for more information.

\chapter{The Duality between $\mathcal{H}^1$ and \textrm{\textbf{BMO}}}

\setcounter{theorem}{0}\setcounter{equation}{0}
The main result (Theorem 2.4) of this chapter is the analogue in our setting
of the famous Fefferman duality theorem between $H^1$ and \textrm{BMO. }

\section{The bounded map from $L^\infty (L^\infty (\Bbb{R)}\otimes \mathcal{M%
},L_c^2)$ to $\mathrm{BMO}_c(\Bbb{R},\mathcal{M})$}

As in the classical case, we will embeds $\mathcal{H}_c^1(\Bbb{R},\mathcal{M})$
into a larger space $L^1(L^\infty (\Bbb{R)}\otimes \mathcal{M},L_c^2)$,
which requires the following maps $\Phi ,\Psi .$

\begin{definition}
We define a map $\Phi $ from $\mathcal{H}_c^p(\Bbb{R},\mathcal{M})\ (1\leq
p<\infty )$ to $L^p(L^\infty (\Bbb{R)}\otimes \mathcal{M},L_c^2(\widetilde{%
\Gamma }))\,$by 
\[
\Phi (f)(x,y,t)=\nabla f(x+t,y)\chi _\Gamma (x,y) 
\]
and a map $\Psi $ for a sufficiently nice $h\in L^p(L^\infty (\Bbb{R)}%
\otimes \mathcal{M},L_c^2(\widetilde{\Gamma }))$ ($1\leq p\leq \infty $) by 
\begin{equation}
\Psi (h)(s)=\int_{\Bbb{R}}\mathop{\displaystyle\int\!\!\!\int}_\Gamma
h(x,y,t)Q_y(x+t-s)dydxdt;\quad \forall s\in \Bbb{R}  \label{hq}
\end{equation}
where, $Q_y(x)$ is defined as a function on $\Bbb{R\times }\widetilde{\Gamma 
}$ by 
\begin{equation}
Q_y(x)(1)=\frac{\partial P_y(x)}{\partial x},\quad Q_y(x)(2)=\frac{\partial
P_y(x)}{\partial y};\forall (x,y)\in \Gamma .  \label{qy}
\end{equation}
\end{definition}

Note that $\Phi $ is simply the natural embedding of $\mathcal{H}_c^p(\Bbb{R}%
,\mathcal{M})$ into $L^p(L^\infty (\Bbb{R)}\otimes \mathcal{M},L_c^2(%
\widetilde{\Gamma }))$. On the other hand, $\Psi $ is well defined for
sufficiently nice $h$, more precisely ''nice'' will mean that $%
h(x,y,t)=\sum_{i=1}^nm_if_i(t)\chi _{A_i}$ with $m_i\in S_{\mathcal{M}%
},A_i\in \widetilde{\Gamma },|A_i|<\infty $ and with scalar valued simple
functions $f_i.$ In this case, it is easy to check that $\Psi (h)\in L^p(%
\mathcal{M},L_c^2(\Bbb{R},\frac{dt}{1+t^2})).$

We will prove that $\Psi $ extends to a bounded map from $L^\infty (L^\infty
(\Bbb{R)}\otimes \mathcal{M},L_c^2(\widetilde{\Gamma }))$ to $\mathrm{BMO}_c(%
\Bbb{R},\mathcal{M})$ (see Lemma 2.2) and also from $L^p(L^\infty (\Bbb{R)}%
\otimes \mathcal{M},L_c^2(\widetilde{\Gamma }))$ to $\mathcal{H}_c^p(\Bbb{R},%
\mathcal{M})$ for all $1<p<\infty $ (see Theorem 4.8). The following
proposition, combined with Theorem 4.8 in Chapter 4, implies that $\Psi $ is
a projection of $L^p(L^\infty (\Bbb{R)}\otimes \mathcal{M},L_c^2(\widetilde{%
\Gamma }))$ onto $\mathcal{H}_c^p(\Bbb{R},\mathcal{M})$ if we identify $%
\mathcal{H}_c^p(\Bbb{R},\mathcal{M})$ with a subspace of $L^p(L^\infty (\Bbb{R)%
}\otimes \mathcal{M},L_c^2(\widetilde{\Gamma }))$ via $\Phi .$

\begin{proposition}
For any $f\in \mathcal{H}_c^p(\Bbb{R},\mathcal{M})$ $(1\leq p<\infty ),$ 
\[
\Psi \Phi (f)=f 
\]
\end{proposition}

\noindent\textbf{Proof.} We have 
\begin{eqnarray*}
&&\int_{-\infty }^{+\infty }\mathop{\displaystyle\int\!\!\!\int}_\Gamma \Phi
(f)\nabla g(t+x,y)dydxdt \\
&=&\int_{-\infty }^{+\infty }\int_{-\infty }^{+\infty }\mathop{\displaystyle%
\int\!\!\!\int}_\Gamma \Phi (f)Q_y(x+t-s)dydxdtg(s)ds.
\end{eqnarray*}
On the other hand, by (\ref{greenf}) we have 
\[
\int_{-\infty }^{+\infty }\mathop{\displaystyle\int\!\!\!\int}_\Gamma \Phi
(f)\nabla g(t+x,y)dydxdt=\int_{-\infty }^{+\infty }f(s)g(s)ds 
\]
for every $g\,$good enough. Therefore 
\[
\int_{-\infty }^{+\infty }\mathop{\displaystyle\int\!\!\!\int}_\Gamma \Phi
(f)Q_y(x+t-s)dydxdt=f(s) 
\]
almost everywhere. This is $\Psi \Phi (f)=f.\qed$

\medskip

We can also prove $\Psi \Phi ({\varphi} )={\varphi} $ by showing directly
the Poisson integral of $\Psi \Phi ({\varphi} )$ coincides with that of ${%
\varphi} ,$ namely 
\begin{equation}
\int_{\Bbb{R}}\Psi \Phi ({\varphi} )(w)P_v(u-w)dw=\int_{\Bbb{R}}{\varphi}
(w)P_v(u-w)dw,\quad \forall (u,v)\in \Bbb{R}_{+}^2.
\end{equation}
$\,$Indeed, using elementary properties of the Poisson kernel, we have 
\begin{eqnarray*}
&&\int_{\Bbb{R}}\Psi \Phi ({\varphi} )(h)P_v(u-h)dh \\
&=&\int_{\Bbb{R}}\int_{\Bbb{R}}\mathop{\displaystyle\int\!\!\!\int}_\Gamma
\int_{\Bbb{R}}{\varphi} (s)\nabla P_y(x+t-s)ds\nabla
P_y(x+t-h)dydxdtP_v(u-h)dh \\
&=&\int_{\Bbb{R}}{\varphi} (s)\mathop{\displaystyle\int\!\!\!\int}_\Gamma
\int_{\Bbb{R}}\int_{\Bbb{R}}\frac \partial {\partial y}P_y(x+t-s)\frac
\partial {\partial y}P_y(x+t-h)P_v(u-h)dtdhdxdyds \\
&&+\int_{\Bbb{R}}{\varphi} (s)\mathop{\displaystyle\int\!\!\!\int}_\Gamma
\int_{\Bbb{R}}\int_{\Bbb{R}}\frac \partial {\partial x}P_y(x+t-s)\frac
\partial {\partial x}P_y(x+t-h)P_v(u-h)dtdhdxdyds \\
&=&\int_{\Bbb{R}}{\varphi} (s)\int_{\Bbb{R}}\mathop{\displaystyle\int\!\!\!%
\int}_{\Bbb{R}_{+}^2}\frac \partial {\partial y}P_y(x-s)\frac \partial
{\partial y}P_y(x-h)2ydydxP_v(u-h)dhds \\
&&+\int_{\Bbb{R}}{\varphi} (s)\int_{\Bbb{R}}\int_{\Bbb{R}}\frac \partial
{\partial s}P_y(x-s)\frac \partial {\partial u}P_{y+v}(x-u)2ydxdyds \\
&=&\int_{\Bbb{R}}{\varphi} (s)\int_0^\infty 2y\frac{\partial ^2}{\partial v^2%
}P_{v+2y}(u-s)dyds-\int_{\Bbb{R}}{\varphi} (s)\int_0^\infty 2y\frac{\partial
^2}{\partial u^2}P_{v+2y}(u-s)dyds \\
&=&\int_{\Bbb{R}}{\varphi} (s)\int_0^\infty y\frac{\partial ^2}{\partial y^2}%
P_{v+2y}(u-s)dyds \\
&=&\int_{\Bbb{R}}{\varphi} (s)(0-\int_0^\infty \frac \partial {\partial
y}P_{v+2y}(u-s)dy)ds \\
&=&\int_{\Bbb{R}}{\varphi} (s)P_v(u-s)ds.\qed
\end{eqnarray*}

\begin{lemma}
$\Psi $ extends to a bounded map from $L^\infty (L^\infty (\Bbb{R)}\otimes 
\mathcal{M},L_c^2(\widetilde{\Gamma }))$ to $\mathrm{BMO}_c(\Bbb{R},\mathcal{%
M})$ of norm controlled by a universal constant.
\end{lemma}

\noindent\textbf{Proof.} Let $\mathcal{S}$ be the family of all $L^\infty (%
\Bbb{R)}\otimes \mathcal{M}$-valued simple functions $h$ which can written
as $h(x,y,t)=\sum_{i=1}^nm_if_i(t)\chi _{A_i}(x,y)$ with $m_i\in S_{\mathcal{%
M}}$, $f_i\in L^\infty (\Bbb{R)}\cap L^1(\Bbb{R)}$ and compact $A_i\subset 
\widetilde{\Gamma }$. (By compact $A_i$ we mean that the two components of $%
A_i$ are compact subsets in $\Gamma $.) Note that $\mathcal{S}$ is w*-dense
in $L^\infty (L^\infty (\Bbb{R)}\otimes \mathcal{M},L_c^2(\widetilde{\Gamma }%
))$ (in fact, the unit ball of $\mathcal{S}$ is w*-dense in the unit ball of 
$L^\infty (L^\infty (\Bbb{R)}\otimes \mathcal{M},L_c^2(\widetilde{\Gamma }))$%
). We will first show that 
\begin{eqnarray}
||\Psi (h)||_{\mathrm{BMO}_c}\leq c\left\| h\right\| _{L^\infty (L^\infty (%
\Bbb{R)}\otimes \mathcal{M},L_c^2)}\ ,\ \forall \;h\in \mathcal{S}. \label{Sh}
\end{eqnarray}
Fix $h\in \mathcal{S}$ and let ${\varphi }=\Psi (h)$. Then ${\varphi }\in
L^\infty (\mathcal{M},L_c^2(\Bbb{R},\frac{dt}{1+t^2}))$ by Proposition 1.3.
To estimate the $\mathrm{BMO}_c$-norm of ${\varphi }$, we fix an interval $I$
and set $h=h_1+h_2$ with 
\begin{eqnarray*}
h_1(x,y,t) &=&h(x,y,t)\chi _{_{2I}}(t) \\
h_2(x,y,t) &=&h(x,y,t)\chi _{_{(2I)^c}}(t).
\end{eqnarray*}
Let 
\[
B_I=\displaystyle\int_{-\infty }^{+\infty }\mathop{\displaystyle\iint }%
_\Gamma Q_Ih_2dydxdt 
\]
with the notation $Q_I(x,t)=\frac 1{|I|}\int_IQ_y(x+t-s)ds.$ Now 
\begin{eqnarray*}
&&\frac 1{|I|}\int_I|{\varphi }(s)-B_I|^2ds \\
&\leq &\frac 2{|I|}\int_I|\int_{(2I)^c}\mathop{\displaystyle\int\!\!\!\int}%
_\Gamma (Q_y(x+t-s)-Q_I)hdxdydt|^2ds \\
&&+\frac 2{|I|}\int_I|\int_{-\infty }^{+\infty }\mathop{\displaystyle\int\!%
\!\!\int}_\Gamma Q_y(x+t-s)h_1dxdydt|^2ds \\
&=&A+B
\end{eqnarray*}
Notice that 
\begin{eqnarray}
\mathop{\displaystyle\int\!\!\!\int}_\Gamma |Q_y(x+t-s)-Q_I|^2dxdy &\leq &c%
\mathop{\displaystyle\int\!\!\!\int}_\Gamma (\frac{|I|}{(|x+t-s|+y)^3})^2dxdy
\nonumber \\
&\leq &c|I|^2(t-C_I)^{-4}  \label{qt}
\end{eqnarray}
for every $t\in (2I)^c$ and $s\in I.$ By (\ref{fg}) 
\[
\left| \mathop{\displaystyle\int\!\!\!\int}_\Gamma
(Q_y(x+t-s)-Q_I)hdxdy\right| ^2\leq c|I|^2(t-C_I)^{-4}\mathop{\displaystyle%
\int \!\!\!\int }_\Gamma h^{*}hdxdy 
\]
and by (\ref{fg}) again, 
\begin{eqnarray*}
&&\left\| A\right\| _{\mathcal{M}} \\
&\leq &c||\int_{(2I)^c}(t-C_I)^{-2}dt\int_{(2I)^c}(t-C_I)^2%
\mathop{\displaystyle\int\!\!\!\int}_\Gamma h^{*}hdxdy|I|^2(t-C_I)^{-4}dt||_{%
\mathcal{M}} \\
&\leq &||\frac c{|I|}\int_{(2I)^c}|I|^2(t-C_I)^{-2}\mathop{\displaystyle\int%
\!\!\!\int}_\Gamma h^{*}hdxdydt||_{\mathcal{M}} \\
&\leq &c\left\| h\right\| _{L^\infty (L^\infty (\Bbb{R)}\otimes \mathcal{M}%
,L_c^2)}^2
\end{eqnarray*}
For the second term $B$, we have 
\begin{eqnarray*}
&&\left\| B\right\| _{\mathcal{M}} \\
&\leq &\frac 2{|I|}||\int_{\Bbb{R}}|\int_{\Bbb{R}}\mathop{\displaystyle\int%
\!\!\!\int}_\Gamma Q_y(x+t-s)h_1dxdydt|^2ds||_{\mathcal{M}} \\
&=&\frac 2{|I|}\sup_{\tau |a|=1}\tau (|a|\int_{\Bbb{R}}|\int_{\Bbb{R}}%
\mathop{\displaystyle\int\!\!\!\int}_\Gamma Q_y(x+t-s)h_1dxdydt|^2ds) \\
&=&\frac 2{|I|}\sup_{\tau |a|=1}\tau \int_{\Bbb{R}}|\int_{\Bbb{R}}%
\mathop{\displaystyle\int\!\!\!\int}_\Gamma Q_y(x+t-s)h_1|a|^{\frac
12}dxdydt|^2ds \\
&=&\frac 2{|I|}\sup_{\tau |a|=1}\sup_{||f||_{_{_{_{L^2(L^\infty (\Bbb{R)}%
\otimes \mathcal{M})}}}}=1}(\tau \int_{\Bbb{R}}f(s)\int_{\Bbb{R}}%
\mathop{\displaystyle\int\!\!\!\int}_\Gamma Q_y(x+t-s)h_1|a|^{\frac
12}dxdydtds)^2 \\
&=&\frac 2{|I|}\sup_{\tau |a|=1}\sup_{||f||_{L^2(L^\infty (\Bbb{R)}\otimes 
\mathcal{M})}=1}(\tau \int_{\Bbb{R}}\mathop{\displaystyle\int\!\!\!\int}%
_\Gamma \nabla f(t+x,y)h_1|a|^{\frac 12}dxdydt)^2
\end{eqnarray*}
Hence by Cauchy-Schwartz inequality and (\ref{h2l2}) 
\begin{eqnarray*}
\left\| B\right\| _{\mathcal{M}} &\leq &\frac 2{|I|}\sup_{\tau |a|=1}\tau
\int_{\Bbb{R}}\mathop{\displaystyle\int\!\!\!\int}_\Gamma h_1^{*}h_1|a|dxdydt
\\
&\leq &\frac 2{|I|}||\int_{\Bbb{R}}\mathop{\displaystyle\int\!\!\!\int}%
_\Gamma h_1^{*}h_1dxdydt||_{\mathcal{M}} \\
&=&\frac 2{|I|}||\int_{2I}\mathop{\displaystyle\int\!\!\!\int}_\Gamma
h^{*}hdxdydt||_{\mathcal{M}} \\
&\leq &4\left\| h\right\| _{L^\infty (L^\infty (\Bbb{R)}\otimes \mathcal{M}%
,L_c^2)}^2
\end{eqnarray*}
Thus 
\[
||{\varphi }||_{\mathrm{BMO}_c}\leq c\left\| h\right\| _{L^\infty (L^\infty (%
\Bbb{R)}\otimes \mathcal{M},L_c^2)}. 
\]
In particular, by Proposition 1.3, 
\[
\Vert {\varphi }\Vert _{L^\infty (\mathcal{M},L_c^2(\Bbb{R},\frac{dt}{1+t^2}%
))}\le c\Vert h\Vert _{L^\infty (L^\infty (\Bbb{R)}\otimes \mathcal{M}%
,L_c^2)}. 
\]

Thus we have proved the boundedness of $\Psi $ from the w*-dense vector
subspace $\mathcal{S}$ of $L^\infty (L^\infty (\Bbb{R)}\otimes \mathcal{M}%
,L_c^2(\widetilde{\Gamma }))$ to $\mathrm{BMO}_c(\Bbb{R},\Bbb{\mathcal{M})}$%
. Now we extend $\Psi $ to the whole $L^\infty (L^\infty (\Bbb{R)}\otimes 
\mathcal{M},L_c^2(\widetilde{\Gamma }))$. To this end we first extend $\Psi $
to a bounded map from $L^\infty (L^\infty (\Bbb{R)}\otimes \mathcal{M},L_c^2(%
\widetilde{\Gamma }))$ into $L^\infty (\mathcal{M},L_c^2(\Bbb{R},\frac{dt}{%
1+t^2}))$. By the discussion above, $\Psi $ is also bounded from $\mathcal{S}
$ to $L^\infty (\mathcal{M},L_c^2(\Bbb{R},\frac{dt}{1+t^2}))$. Let $H_0^1$
be the subspace of all $f\in H^1(\Bbb{R)}$ such that $(1+t^2)f(t)\in L^2(%
\Bbb{R)}$. Let $L^1(\mathcal{M})\otimes H_0^1$ denote the algebraic tensor
product of $L^1(\mathcal{M})$ and $H_0^1$. Note that 
\[
L^1(\mathcal{M})\otimes H_0^1\subset \mathcal{H}_c^1(\Bbb{R},\mathcal{\ M}%
),\ \ L^1(\mathcal{M})\otimes H_0^1\subset L^1(\mathcal{M},L_c^2(\Bbb{R},%
\frac{dt}{1+t^2})) 
\]
and $L^1(\mathcal{M})\otimes H_0^1$ is dense in both of the latter spaces.
Moreover, it is easy to see that for any $h\in \mathcal{S}$ and $f\in L^1(%
\mathcal{M})\otimes H_0^1$ 
\[
\tau \int_{-\infty }^{+\infty }\mathop{\displaystyle\int\!\!\!\int}_\Gamma
h^{*}(x,y,t)\nabla f(t+x,y)dydxdt=\tau \int_{-\infty }^{+\infty }\Psi
(h)^{*}(s)f(s)ds. 
\]
Then it follows that $\Psi $ is continuous from $\big(\mathcal{S},\ \sigma (%
\mathcal{S},L^1(L^\infty (\Bbb{R)}\otimes \mathcal{M},L_c^2(\widetilde{%
\Gamma })))\big)
$ to $\big(L^\infty (\mathcal{M},L_c^2(\Bbb{R},\frac{dt}{1+t^2})),\ \sigma
(L^\infty (\mathcal{M},L_c^2(\Bbb{R},\frac{dt}{1+t^2})),L^1(\mathcal{M}%
)\otimes H_0^1)\big)$.

Now given $f\in L^1(\mathcal{M})\otimes H_0^1$ we define $\Psi _{*}(f):%
\mathcal{S}\to \Bbb{C}$ by 
\[
\Psi _{*}(f)(h)=\tau \int_{-\infty }^{+\infty }\Psi (h)^{*}(s)f(s)ds. 
\]
Then $\Psi _{*}(f)$ is an anti-linear functional on $\mathcal{S}$ continuous
with respect to the w*-topology; hence $\Psi _{*}(f)$ extends to a
w*-continuous anti-linear functional on $L^\infty (L^\infty (\Bbb{R)}\otimes 
\mathcal{M},L_c^2(\widetilde{\Gamma })))$, i.e. an element in $L^1(L^\infty (%
\Bbb{R)}\otimes \mathcal{M},L_c^2(\widetilde{\Gamma })))$, still denoted by $%
\Psi _{*}(f)$. By the w*-density of $\mathcal{S}$ in $L^\infty (L^\infty (%
\Bbb{R)}\otimes \mathcal{M},L_c^2(\widetilde{\Gamma })))$, this extension is
unique. Therefore, we have defined a map 
\[
\Psi _{*}:L^1(\mathcal{M})\otimes H_0^1\to L^1(L^\infty (\Bbb{R)}\otimes 
\mathcal{M},L_c^2(\widetilde{\Gamma })). 
\]
The above uniqueness of the extension $\Psi _{*}(f)$ for any given $f$
implies that $\Psi _{*}$ is linear. On the other hand, by what we already
proved in the previous part, we have 
\begin{eqnarray*}
|\Psi _{*}(f)(h)| &\le &\Vert f\Vert _{L^1(\mathcal{M},L_c^2(\Bbb{R},\frac{dt%
}{1+t^2}))}\Vert \Psi (h)\Vert _{L^\infty (\mathcal{M},L_c^2(\Bbb{R},\frac{dt%
}{1+t^2}))} \\
&\le &c\,\Vert f\Vert _{L^1(\mathcal{M},L_c^2(\Bbb{R},\frac{dt}{1+t^2}%
))}\Vert h\Vert _{L^\infty (L^\infty (\Bbb{R)}\otimes \mathcal{M},L_c^2)}\ .
\end{eqnarray*}
Since the unit ball of $\mathcal{S}$ is w*-dense in the unit ball of $%
L^\infty (L^\infty (\Bbb{R)}\otimes \mathcal{M},L_c^2(\widetilde{\Gamma })))$%
, it follows that 
\[
\Psi _{*}:(L^1(\mathcal{M})\otimes H_0^1,\left\| \cdot \right\| _{L^1(%
\mathcal{M},L_c^2(\Bbb{R},\frac{dt}{1+t^2}))})\to L^1(L^\infty (\Bbb{R)}%
\otimes \mathcal{M},L_c^2(\widetilde{\Gamma })) 
\]
is bounded and its norm is majorized by $c$. This, together with the density
of $L^1(\mathcal{M})\otimes H_0^1$ in $L^1(\mathcal{M},L_c^2(\Bbb{R},\frac{dt%
}{1+t^2}))$ implies that $\Psi _{*}$ extends to a unique bounded map from $%
L^1(\mathcal{M},L_c^2(\Bbb{R},\frac{dt}{1+t^2}))$ into $L^1(L^\infty (\Bbb{R)%
}\otimes \mathcal{M},L_c^2(\widetilde{\Gamma })))$, still denoted by $\Psi
_{*}$. Consequently, the adjoint $(\Psi _{*})^{*}$ of $\Psi _{*}$ is bounded
from $L^\infty (L^\infty (\Bbb{R)}\otimes \mathcal{M},L_c^2(\widetilde{%
\Gamma })))$ to $L^\infty (\mathcal{M},L_c^2(\Bbb{R},\frac{dt}{1+t^2}))$
(noting that this adjoint is taken with respect to the anti-dualities). By
the very definition of $\Psi _{*}$, we have 
\[
(\Psi _{*})^{*}\big|_{\mathcal{S}}=\Psi . 
\]
This shows that $(\Psi _{*})^{*}$ is an extension of $\Psi $ from $L^\infty
(L^\infty (\Bbb{R)}\otimes \mathcal{M},L_c^2(\widetilde{\Gamma }))$ to $%
L^\infty (\mathcal{M},L_c^2(\Bbb{R},\frac{dt}{1+t^2})),$ which we denote by $%
\Psi $ again. Being an adjoint, $\Psi $ is w*-continuous.

It remains to show that the so extended map $\Psi $ really takes values in $%
\mathrm{BMO}_c(\Bbb{R},\Bbb{\mathcal{M})}$. Given a bounded interval $I\subset 
\Bbb{R},$ the
w*-topology of $L^\infty (\mathcal{M},L_c^2(\Bbb{R},\frac{dt}{1+t^2}))$
induces a topology in $L^\infty (\mathcal{M},L_c^2(I))$ equivalent to the
w*-topology in $L^\infty (\mathcal{M},L_c^2(I))$. Then by the w*-continuity
of $\Psi $, we deduce that, for every ${\varepsilon }>0,I\subset \Bbb{R}%
,f\in L^1(\mathcal{M},L_c^2(I)),$ there exists a $h\in \mathcal{S}$ such that
\begin{eqnarray}
&&\tau \int_If^{*}(\Psi (g)(t)-\Psi (g)_I)dt \nonumber\\
&\leq &\tau \int_If^{*}(\Psi (h)(t)-\Psi (h)_I)dt+{\varepsilon } \nonumber\\
&\leq &\left\| \Psi (h)(t)-\Psi (h)_I\right\| _{L^\infty (\mathcal{M}%
,L_c^2(I))}\left\| f\right\| _{L^1(\mathcal{M},L_c^2(I))}+{\varepsilon } 
\label{2g}
\end{eqnarray}
and
\begin{eqnarray}
\Vert h\Vert _{L^\infty (L^\infty (\Bbb{R)}\otimes \mathcal{M},L_c^2(%
\widetilde{\Gamma }))}&\le &\Vert g\Vert _{L^\infty (L^\infty (\Bbb{R)}\otimes 
\mathcal{M},L_c^2(\widetilde{\Gamma }))}+\varepsilon \label{2h}
\end{eqnarray}
Combining (\ref{2g}), (\ref{2h}) and (\ref{Sh}) we get
\begin{eqnarray*}
&&\int_If^{*}(\Psi (g)(t)-\Psi (g)_I)dt\\
&\le &c|I|\left\| h\right\| _{L^\infty (L^\infty (\Bbb{R)}\otimes \mathcal{M}
,L_c^2(\widetilde{\Gamma }))}\left\| f\right\| _{L^1(\mathcal{M},L_c^2(I))}+{%
\varepsilon } \\
&\le &c|I|(\left\| g\right\| _{L^\infty (L^\infty (\Bbb{R)}\otimes \mathcal{M}
,L_c^2(\widetilde{\Gamma }))}+\varepsilon) \left\| f\right\| _{L^1(\mathcal{M},L_c^2(I))}+{%
\varepsilon } 
\end{eqnarray*}
By letting ${\varepsilon }\rightarrow 0$ and taking supremum over all $\left\|
f\right\| _{L^1(L^\infty (\Bbb{R)}\otimes \mathcal{M},L_c^2(\widetilde{%
\Gamma }))}\leq 1$ and $I\subset \Bbb{R},$ we get $\Psi (g)\in \mathrm{BMO}%
_c(\Bbb{R},\Bbb{\mathcal{M})}$ and 
\[
||\Psi (g)||_{\mathrm{BMO}_c}\leq c\left\| g\right\| _{L^\infty (L^\infty (%
\Bbb{R)}\otimes \mathcal{M},L_c^2)}. 
\]
Therefore, we have extended $\Psi $ to a bounded map from $L^\infty
(L^\infty (\Bbb{R)}\otimes \mathcal{M},L_c^2(\widetilde{\Gamma }))$ to $%
\mathrm{BMO}_c(\Bbb{R},\mathcal{M})$, thus completing the proof of the lemma.%
\qed

\medskip \noindent\textbf{Remark.} We sketch an alternate proof of the fact
that ${\varphi} =\Psi (h)$ is in $\mathrm{BMO}_c(\Bbb{R},\Bbb{\mathcal{M})}$
for $h\in \mathcal{S}$. Let $H$ be the Hilbert space on which $\mathcal{M\,}$%
acts. Recall that $\mathcal{M}_{*}$ is a quotient space of $B(H)_{*}$ by the
preannihilator of $\mathcal{M}.$ Denote the quotient map by $q.$ For every $%
a,b\in H,$ denote $[a\otimes b]=q(a\otimes b).$ Note that $\tau
(m^{*}[a\otimes b])=\tau ([m^{*}(a\otimes b)])=\left\langle m(b),\overline{a}%
\right\rangle ,\forall m\in \mathcal{M}.$ From (\ref{bmoh}) and the
classical duality between $\mathrm{BMO}(\Bbb{R},H)$ and $H^1(\Bbb{R},H),$%
\begin{eqnarray}
||{\varphi} ||_{\mathrm{BMO}_c(\Bbb{R},\mathcal{M})} &=&\sup_{e\in H,\left\|
e\right\| _H=1}||{\varphi} e||_{\mathrm{BMO}(\Bbb{R},H).}  \nonumber \\
&\leq &c\sup_{e\in H,\left\| e\right\| _H=1}\sup_{\left\| g\right\| _{H^1(%
\Bbb{R},H)}=1}\left| \int_{-\infty }^{+\infty }\left\langle {\varphi} (e),%
\overline{g}\right\rangle dt\right|  \nonumber \\
&=&c\sup_{e\in H,\left\| e\right\| _H=1}\sup_{\left\| g\right\| _{H^1(\Bbb{R}%
,H)}=1}\left| \tau \int_{-\infty }^{+\infty }{\varphi} ^{*}[g\otimes
e]dt\right|  \label{ge}
\end{eqnarray}
Let $f=[g\otimes e].$ Noting that 
\[
|\nabla f|^2=\left\langle \nabla g,\nabla g\right\rangle [e\otimes
e]=|\nabla g|^2[e\otimes e], 
\]
we get 
\begin{equation}
\tau \left( S_c(f)(t)\right) =(\mathop{\displaystyle\int\!\!\!\int}_\Gamma
\left| \nabla g(t+x,y)\right| ^2dxdy)^{\frac 12}.  \label{gf}
\end{equation}
Thus $\left\| f\right\| _{\mathcal{H}_c^1(\Bbb{R},\mathcal{M})}=1\,$if $%
\left\| g\right\| _{H^1(\Bbb{R},H)}=1\,$and $\left\| e\right\| _H=1.\,$
Therefore 
\begin{eqnarray*}
||{\varphi} ||_{\mathrm{BMO}_c(\Bbb{R},\mathcal{M})} &\leq &c\;\sup_{\left\|
f\right\| _{\mathcal{H}_c^1(\Bbb{R},\mathcal{M})}=1}\left| \tau
\int_{-\infty }^{+\infty }{\varphi} ^{*}fdt\right| \\
&=&c\;\sup_{\left\| f\right\| _{\mathcal{H}_c^1(\Bbb{R},\mathcal{M}%
)}=1}\left| \tau \int_{-\infty }^{+\infty }\mathop{\displaystyle\int\!\!\!%
\int}_\Gamma h^{*}(x,y,t)\nabla f(t+x,y)dydxdt\right| \\
&\leq &c\left\| h\right\| _{L^\infty (L^\infty (\Bbb{R)}\otimes \mathcal{M}%
,L_c^2)}.
\end{eqnarray*}

\begin{corollary}
Let $f\in L^1(\mathcal{M},L_c^2(\Bbb{R},(1+s^2)ds))$ with $\int fds=0.$ Then 
$f\in \mathcal{H}_c^1(\Bbb{R},\mathcal{M})$ and 
\[
\left\| f\right\| _{\mathcal{H}_c^1}\leq c\left\| f\right\| _{L^1(\mathcal{M}%
,L_c^2(\Bbb{R},(1+s^2)ds))} 
\]
\end{corollary}

\noindent\textbf{Proof.} By Lemma 2.2, the assumption that $\int fds=0$ and
Proposition 1.3, we have 
\begin{eqnarray*}
\left\| f\right\| _{\mathcal{H}_c^1} &=&||\nabla f(t+x,y)\chi _\Gamma
||_{L^1(L^\infty (\Bbb{R)}\otimes \mathcal{M},L_c^2)} \\
&=&\sup_{\left\| h\right\| _{L^\infty (L^\infty (\Bbb{R)}\otimes \mathcal{M}%
,L_c^2)}\leq 1}\left| \tau \int \mathop{\displaystyle\int\!\!\!\int}_\Gamma
h^{*}\nabla f(t+x,y)dxdydt\right| \\
&=&\sup_{\left\| h\right\| _{L^\infty (L^\infty (\Bbb{R)}\otimes \mathcal{M}%
,L_c^2)}\leq 1}\left| \tau \int_{\Bbb{R}}(\Psi (h))^{*}(s)f(s)ds\right| \\
&\leq &c\sup_{\left\| {\varphi} \right\| _{\mathrm{BMO}_c(\Bbb{R},\mathcal{M}%
)}\leq 1}\left| \tau \int_{\Bbb{R}}{\varphi} ^{*}(s)f(s)ds\right| \\
&\leq &c\sup_{\left\| {\varphi} \right\| _{L^\infty (\mathcal{M},L_c^2(\Bbb{R%
},\frac{ds}{1+s^2}))}\leq 1}\left| \tau \int_{\Bbb{R}}{\varphi}
^{*}(s)(1+s^2)f(s)\frac{ds}{1+s^2}\right| \\
&\leq &c\left\| (1+s^2)f(s)\right\| _{L^1(\mathcal{M},L_c^2(\Bbb{R},\frac{ds%
}{1+s^2}))} \\
&=&c\left\| f\right\| _{L^1(\mathcal{M},L_c^2(\Bbb{R},(1+s^2)ds))}.\qed
\end{eqnarray*}

\medskip \noindent
\textbf{Remark.} In particular, every $S_{\mathcal{M}}$-valued simple
function $f$ with $\int fds=0$ is in $\mathcal{H}_c^1(\Bbb{R},\mathcal{M}).$
Consequently, by the remark before Proposition 1.3, $\mathcal{H}_c^1(\Bbb{R},%
\mathcal{M})\cap \mathcal{H}_c^p(\Bbb{R},\mathcal{M})$ is dense in $\mathcal{%
H}_c^p(\Bbb{R},\mathcal{M})\ (p>1).$

\section{The duality theorem of operator valued $\mathcal{H}^1$ and \textrm{BMO}}

Denote by $\mathcal{H}_{c0}^1(\Bbb{R},\mathcal{M})$ (resp. $\mathcal{H}%
_{r0}^1(\Bbb{R},\mathcal{M})$) the family of functions $f$ in $\mathcal{H}%
_c^1(\Bbb{R},\mathcal{M})$ (resp. $\mathcal{H}_r^1(\Bbb{R},\mathcal{M}),%
\mathcal{H}_{cr}^1(\Bbb{R},\mathcal{M}) $) such that $f\in L^1(\mathcal{M}%
,L_c^2(\Bbb{R},(1+t^2)dt))\,($resp. $L^1(\mathcal{M},$ $L_r^2(\Bbb{R}%
,(1+t^2)dt))\,.$ It is easy to see that $\mathcal{H}_{c0}^1(\Bbb{R},\mathcal{%
M})$ (resp. $\mathcal{H}_{r0}^1(\Bbb{R},\mathcal{M})$) is a dense subspace
of $\mathcal{H}_c^1(\Bbb{R},\mathcal{M})$ (resp. $\mathcal{H}_r^1(\Bbb{R},%
\mathcal{M}))$). Let 
\[
\mathcal{H}_{cr0}^1(\Bbb{R},\mathcal{M})=\mathcal{H}_{c0}^1(\Bbb{R},\mathcal{%
M})+\mathcal{H}_{r0}^1(\Bbb{R},\mathcal{\ M}). 
\]
Then $\mathcal{H}_{cr0}^1(\Bbb{R},\mathcal{M})$ is a dense subspace of $%
\mathcal{H}_{cr}^1(\Bbb{R},\mathcal{M})$. Recall that we have proved in
Chapter 1 that $\mathrm{BMO}_c(\Bbb{R},\mathcal{M})\subseteq L^\infty (%
\mathcal{M},L_c^2(\Bbb{R},\frac{dt}{1+t^2})).$ Thus by Proposition 1.1 $%
\langle {\varphi} ,f\rangle =\int_{-\infty }^{+\infty }{\varphi} ^{*}fdt$
exists in $L^1(\mathcal{M})$ for all ${\varphi} \in \mathrm{BMO}_c(\Bbb{R},%
\mathcal{M})$ and $f\in \mathcal{H}_{c0}^1(\Bbb{R},\mathcal{M})$ (see our
convention after Proposition 1.1).

\begin{theorem}
(a) We have $(\mathcal{H}_c^1(\Bbb{R},\mathcal{M}))^{*}=\mathrm{BMO}_c(\Bbb{R%
},\mathcal{M})$ with equivalent norms. More precisely, every ${\varphi }\in 
\mathrm{BMO}_c(\mathcal{M})\,$defines a continuous linear functional on $%
\mathcal{H}_c^1(\Bbb{R},\mathcal{M})\,$by 
\begin{equation}
l_\varphi (f)=\tau \int_{-\infty }^{+\infty }{\varphi }^{*}fdt;\qquad
\forall f\in \mathcal{H}_{c0}^1(\Bbb{R},\mathcal{M}).  \label{lfai}
\end{equation}
Conversely, every $l\in (\mathcal{H}_c^1(\Bbb{R},\mathcal{M}))^{*}$ can be
given as above by some ${\varphi }\in \mathrm{BMO}_c(\Bbb{R},\mathcal{M})$.
Moreover, there exists a universal constant $c>0$ such that 
\[
c^{-1}\left\| {\varphi }\right\| _{\mathrm{BMO}_c}\leq \left\| l_\varphi
\right\| _{(\mathcal{H}_c^1)^{*}}\leq c\left\| {\varphi }\right\| _{\mathrm{%
BMO}_c}\ . 
\]
Thus $(\mathcal{H}_c^1(\Bbb{R},\mathcal{M}))^{*}=\mathrm{BMO}_c(\Bbb{R},%
\mathcal{M})$ with equivalent norms.

(b) Similarly, $(\mathcal{H}_r^1(\Bbb{R},\mathcal{M}))^{*}=\mathrm{BMO}_r(%
\Bbb{R},\mathcal{M})$ with equivalent norms.

(c) $(\mathcal{H}_{cr}^1(\Bbb{R},\mathcal{M}))^{*}=\mathrm{BMO}_{cr}(\Bbb{R},%
\mathcal{M})$ with equivalent norms.
\end{theorem}

Our proof of Theorem 2.4 requires two technical variants of the square
functions $G_c(f)$ and $S_c(f).$ These are operator valued
functions defined as follows: 
\begin{eqnarray}
G_c(f)(x,y)=(\int_y^\infty |\nabla f(x,s)|^2sds)^{\frac 12}, \label{Gc2}
\\
S_c(f)(x,y)=(\mathop{\displaystyle\int\!\!\!\int}_{\Gamma (0,y)}|\nabla
f(t+x,s)|^2dtds)^{\frac 12} \label{Sc2}
\end{eqnarray}
where $y\geq 0, \Gamma (0,y)=\{(t,s):|t|<s-y,s\geq y\}$ and $f$ is $S_{\mathcal{M}}$-valued simple function.
Note that $%
G_c(f)(x,0)$ and $S_c(f)(x,0)$ are just $G_c(f)$ and $S_c(f)\ $defined in
Chapter 1.

\begin{lemma}
\[
G_c(f)(x,y)\leq 2\sqrt{2}S_c(f)(x,\frac y2)\ . 
\]
\end{lemma}

\noindent\textbf{Proof.} It suffices to prove this inequality for $x=0.\,$%
Let us denote by $B_s$ the ball centered at $(0,s)$ and tangent to the
boundary of $\Gamma (0,\frac y2),\forall s>y.$ By the harmonicity of $\nabla
f,$ we get 
\[
\nabla f(0,s)=\frac 2{\pi (s-\frac y2)^2}\int_{B_s}\nabla f(x,u)dxdu 
\]
By (\ref{fu}), 
\[
|\nabla f(0,s)|^2\leq \frac 8{\pi s^2}\int_{B_s}|\nabla f(x,u)|^2dxdu 
\]
Integrating this inequality, we obtain 
\begin{equation}
\int_y^\infty s|\nabla f(0,s)|^2ds\leq \int_y^\infty \frac 8{\pi
s}\int_{B_s}|\nabla f(x,u)|^2dxduds  \label{lem}
\end{equation}
However $(x,u)\in B_s$ clearly implies that $\frac u2\leq s\leq 4u.$ Thus,
the right hand side of (\ref{lem}) is majorized by 
\[
\int_{\Gamma (0,\frac y2)}|\nabla f(x,u)|^2\int_{\frac u2}^{4u}\frac 8{\pi
s}dsdxdu\leq 8S_c^2(f)(0,\frac y2) 
\]
Therefore $G_c(f)(0,y)\leq 2\sqrt{2}S_c(f)(0,\frac y2).\bigskip $\qed

\medskip

\noindent\textbf{Proof of Theorem 2.4.} (i) We will first prove 
\begin{equation}
|l_{\varphi} (f)|\leq c\left\| {\varphi} \right\| _{\mathrm{BMO}_c}\left\|
f\right\| _{\mathcal{H}_c^1}  \label{lf}
\end{equation}
when both $f$ and ${\varphi} $ have compact support. Once this is done, by
Lemma 1.5, we can see (\ref{lf}) holds for any ${\varphi} \in \mathrm{BMO}_c(%
\Bbb{R},\mathcal{M})$ and any compactly supported $f\in \mathcal{H}_{c0}^1(%
\Bbb{R},\mathcal{M})$. Then recall that by Proposition 1.3 
\[
\mathrm{BMO}_c(\Bbb{R},\mathcal{M})\subset L^\infty (\mathcal{M},L_c^2(\Bbb{R%
},\frac{dt}{1+t^2})) 
\]
and by Corollary 2.3 $\ $%
\[
\left\| f\right\| _{\mathcal{H}_c^1}\leq c\left\| f\right\| _{L^1(\mathcal{M}%
,L_c^2(\Bbb{R},(1+t^2)dt))},~\forall f\in \mathcal{H}_{c0}^1(\Bbb{R},%
\mathcal{M}), 
\]
we deduce (\ref{lf}) for all ${\varphi} \in \mathrm{BMO}_c(\Bbb{R},\mathcal{M%
}),f\in \mathcal{H}_{c0}^1(\Bbb{R},\mathcal{M})$ by choosing compactly
supported $f_n\in \mathcal{H}_{c0}^1(\Bbb{R},\mathcal{M})\rightarrow f$ in $%
L^1(\mathcal{M},L_c^2(\Bbb{R},(1+t^2)dt)).$ Finally, from the density of $%
\mathcal{H}_{c0}^1(\Bbb{R},\mathcal{M})$ in $\mathcal{H}_c^1(\Bbb{R},%
\mathcal{M}),$ $l_{\varphi} $ defined in (\ref{lfai}) extends to a
continuous functional on $\mathcal{H}_c^1(\Bbb{R},\mathcal{M}).$

\medskip

Let us now prove (\ref{lf}) for compactly supported $f\in \mathcal{H}_{c0}^1(%
\Bbb{R},\mathcal{M})$ and compactly supported ${\varphi }\in \mathrm{BMO}_c(%
\Bbb{R},\mathcal{M})$. By approximation we may assume that $\tau \,$is
finite and $G_c(f)(x,y)$ is invertible in $\mathcal{M}\,$for every $(x,y)\in 
\Bbb{R}_{+}^2.$ Recall that $\triangle ({\varphi }^{*}f)=2\nabla {\varphi }%
^{*}\nabla f.$ By the Green theorem and the Cauchy-Schwartz inequality 
\begin{eqnarray*}
&&|l_\varphi (f)| \\
&=&2|\tau \mathop{\displaystyle\int\!\!\!\int}_{\Bbb{R}_{+}^2}\nabla {%
\varphi }^{*}\nabla fydydx| \\
\ &\leq &2(\tau \mathop{\displaystyle\int\!\!\!\int}_{\Bbb{R}%
_{+}^2}G_c^{-\frac 12}(f)|\nabla f|^2G_c^{-\frac 12}(f)ydydx)^{\frac
12}(\tau \mathop{\displaystyle\int\!\!\!\int}_{\Bbb{R}_{+}^2}G_c^{\frac
12}(f)|\nabla {\varphi }|^2G_c^{\frac 12}(f)ydydx)^{\frac 12} \\
&=&2(\tau \mathop{\displaystyle\int\!\!\!\int}_{\Bbb{R}_{+}^2}G_c^{-1}(f)|%
\nabla f|^2ydydx)^{\frac 12}(\tau \mathop{\displaystyle\int\!\!\!\int}_{\Bbb{%
R}_{+}^2}G_c(f)|\nabla {\varphi }|^2ydydx)^{\frac 12} \\
&=&2I\bullet II,
\end{eqnarray*}
Note here $G_c(f)$ is the function of two variables defined by (\ref{Gc2}),
which is differentiable in the weak-* sense. For I we have 
\begin{eqnarray*}
I^2 &=&\tau \int_{-\infty }^{+\infty }\int_0^\infty -G_c^{-1}(f)\frac{%
\partial G_c^2(f)}{\partial y}dydx \\
\ &=&\tau \int_{-\infty }^{+\infty }\int_0^\infty (-G_c^{-1}(f)\frac{%
\partial G_c(f)}{\partial y}G_c(f)-\frac{\partial G_c(f)}{\partial y})dydx \\
\ &=&2\tau \int_{-\infty }^{+\infty }\int_0^\infty -\frac{\partial G_c(f)}{%
\partial y}dydx \\
&=&2\tau \int_{-\infty }^{+\infty }G_c(f)(x,0)dx \\
&\leq &4\sqrt{2}\tau \int_{-\infty }^{+\infty }S_c(f)(x,0)dx \\
&=&4\sqrt{2}\left\| f\right\| _{\mathcal{H}_c^1}.
\end{eqnarray*}
To estimate II, we create a square net partition in $\Bbb{R}_{+}^2$ as
follows:

\[
\sigma (i,j)=\{(x,y):(i-1)2^j<x\leq i2^j,2^j\leq y<2^{j+1}\},\quad \forall
i,j\in \Bbb{Z}. 
\]
Let $C_{i,j}$ denote the center of $\sigma (i,j).$ Define 
\begin{eqnarray*}
\widetilde{S}_c(f)(x,y) &=&S_c(f)(C_{i,j}),\quad \forall (x,y)\in \sigma
(i,j), \\
d_k(x) &=&\widetilde{S}_c(f)(x,2^k)-\widetilde{S}_c(f)(x,2^{k+1}),\quad
\forall x\in \Bbb{R}.
\end{eqnarray*}
It is easy to check that 
\begin{eqnarray}
S_c(f)(x,2y) &\leq &\widetilde{S}_c(f)(x,y)\leq S_c(f)(x,\frac y2), 
\nonumber \\
d_k(x) &\geq &0,\quad \forall x\in \Bbb{R},  \nonumber \\
\widetilde{S}_c(f)(x,y) &=&\sum_{k=j}^\infty d_k(x),\quad \forall 2^j\leq
y<2^{j+1},  \nonumber \\
S_c(f)(x,0) &=&\sum_{k=-\infty }^\infty d_k(x).  \label{sc}
\end{eqnarray}
Now by Lemma 2.5 and (\ref{sc}) 
\begin{eqnarray*}
II^2 &=&\tau \int_{-\infty }^{+\infty }\int_0^\infty G_c(f)(x,y)|\nabla {%
\varphi }|^2ydydx \\
&\leq &2\sqrt{2}\tau \int_{-\infty }^{+\infty }\int_0^\infty \widetilde{S}%
_c(f)(x,\frac y4)|\nabla {\varphi }|^2ydydx \\
&=&2\sqrt{2}\tau \int_{-\infty }^{+\infty }\sum_{k=-\infty }^\infty 
\widetilde{S}_c(f)(x,2^k)\int_{2^{k+2}}^{2^{k+3}}|\nabla {\varphi }|^2ydydx
\\
&=&2\sqrt{2}\tau \int_{-\infty }^{+\infty }\sum_{k=-\infty }^\infty
(\sum_{j=k}^\infty d_j(x))\int_{2^{k+2}}^{2^{k+3}}|\nabla {\varphi }|^2ydydx
\\
&=&2\sqrt{2}\tau \int_{-\infty }^{+\infty }\sum_{j=-\infty }^\infty
d_j(x)\int_0^{2^{j+3}}|\nabla {\varphi }|^2ydydx \\
&=&2\sqrt{2}\tau \sum_{i=-\infty }^\infty \sum_{j=-\infty }^\infty
d_j(i2^j)\int_{(i-1)2^j}^{i2^j}\int_0^{2^{j+3}}|\nabla {\varphi }|^2ydydx
\end{eqnarray*}
Hence by Lemma 1.4 
\begin{eqnarray*}
II^2 &\leq &c\tau \sum_{i=-\infty }^\infty \sum_{j=-\infty }^\infty
d_j(i2^j)2^j\left\| {\varphi }\right\| _{\mathrm{BMO}_c}^2 \\
&=&c\left\| {\varphi }\right\| _{\mathrm{BMO}_c}^2\tau \sum_{j=-\infty
}^\infty \int_{-\infty }^{+\infty }d_j(x)dx \\
&=&c\left\| {\varphi }\right\| _{\mathrm{BMO}_c}^2\tau \int_{-\infty
}^{+\infty }S_c(f)(x,0)dx \\
&=&c\left\| {\varphi }\right\| _{\mathrm{BMO}_c}^2\left\| f\right\| _{%
\mathcal{H}_c^1}.
\end{eqnarray*}
Combining the preceding estimates on I and II, we get $\,$ 
\[
|l_\varphi (f)|\leq c\left\| {\varphi }\right\| _{\mathrm{BMO}_c}\left\|
f\right\| _{\mathcal{H}_c^1}. 
\]
Therefore, $l_\varphi $ defines a continuous functional on $\mathcal{H}_c^1$
of norm smaller than $c\left\| {\varphi }\right\| _{\mathrm{BMO}_c}.\medskip 
$

(ii) Now suppose $l\in (\mathcal{H}_c^1(\Bbb{R},\mathcal{M}))^{*}.$ Then by
the Hahn-Banach theorem $l$ extends to a continuous functional on $%
L^1(L^\infty (\Bbb{R)}\otimes \mathcal{M},L_c^2(\widetilde{\Gamma }))$ of
the same norm. Thus by 
\[
(L^1(L^\infty (\Bbb{R)}\otimes \mathcal{M},L_c^2(\widetilde{\Gamma }%
)))^{*}=L^\infty (L^\infty (\Bbb{R)}\otimes \mathcal{M},L_c^2(\widetilde{%
\Gamma })) 
\]
there exists $g\in L^\infty (L^\infty (\Bbb{R)}\otimes \mathcal{M},L_c^2(%
\widetilde{\Gamma }))$ such that 
\[
||g||_{L^\infty (L^\infty (\Bbb{R)}\otimes \mathcal{M},L_c^2(\widetilde{%
\Gamma }))}^2=\sup_{t\in \Bbb{R}}||\mathop{\displaystyle\int\!\!\!\int}%
_\Gamma g^{*}(x,y,t)g(x,y,t)dydx||_{L^\infty (\Bbb{R)}\otimes \mathcal{M}%
}=||l||^2 
\]
and 
\[
l(f)=\tau \int_{-\infty }^{+\infty }\mathop{\displaystyle\int\!\!\!\int}%
_\Gamma g^{*}(x,y,t)\nabla f(t+x,y)dydxdt,\ \forall \;f\in \mathcal{H}%
_{c0}^1(\Bbb{R},\mathcal{M}). 
\]
Let ${\varphi }=\Psi (g)$, where $\Psi $ is the extension given by Lemma
2.2. By that lemma ${\varphi }\in \mathrm{BMO}_c(\Bbb{R},\mathcal{M})$ and 
\[
||{\varphi }||_{\mathrm{BMO}_c}\leq c||g||_{L^\infty (L^\infty (\Bbb{R)}%
\otimes \mathcal{M},L_c^2(\widetilde{\Gamma }))}=c\Vert l\Vert . 
\]
Then we must show that 
\[
l(f)=\tau \int_{-\infty }^{+\infty }{\varphi }^{*}(s)f(s)ds,\ \forall \;f\in 
\mathcal{H}_{c0}^1(\Bbb{R},\mathcal{M}). 
\]
But this follows from the second part of the proof of Lemma 2.2 in virtue of
the w*-continuity of $\Psi $. Therefore, we have accomplished the proof of
the theorem concerning $\mathcal{H}_c^1(\Bbb{R},\mathcal{M})$ and $\mathrm{%
BMO}_c(\Bbb{R},\mathcal{M}).$ Passing to adjoints yields the part on $%
\mathcal{H}_r^1(\Bbb{R},\mathcal{M})$ and $\mathrm{BMO}_r.\,$Finally, the
duality between $\mathcal{H}_{cr}^1(\Bbb{R},\mathcal{M})$ and $\mathrm{BMO}%
_{cr}(\Bbb{R},\mathcal{M})\mathrm{\,}$is obtained by the classical fact that
the dual of a sum is the intersection of the duals.\qed

\begin{corollary}
${\varphi }\in \mathrm{BMO}_c(\Bbb{R},\mathcal{M})$ if and only if $d{%
\lambda }_\varphi =|\nabla {\varphi }|^2ydxdy\,$is an $\mathcal{M}$-valued
Carleson measure on $\Bbb{R}_{+}^2,$ and $c^{-1}N({\lambda }_\varphi )\leq
\left\| {\varphi }\right\| _{\mathrm{BMO}_c}^2\leq cN({\lambda }_\varphi ).$
\end{corollary}

\noindent\textbf{Proof.} From the first part of the proof of Theorem 2.4, if 
${\varphi} $ is such that $d{\lambda} _{\varphi} =|\nabla {\varphi} |^2ydxdy$
is an $\mathcal{M}$-valued Carleson measure, then ${\varphi} $ defines a
continuous linear functional $l_{\varphi} =\tau \int_{-\infty }^{+\infty }{%
\varphi} ^{*}fdt$ on $\mathcal{H}_{c0}^1(\Bbb{R},\mathcal{M})$ and 
\[
\left\| l_{\varphi} \right\| _{(\mathcal{H}_c^1)^{*}}\leq cN^{\frac 12}({%
\lambda} _{\varphi} ) 
\]
Therefore by Theorem 2.4 again there exists a function ${\varphi} ^{\prime
}\in $\textrm{BMO}$_c(\Bbb{R},\mathcal{M})$ with $\left\| {\varphi} ^{\prime
}\right\| _{\mathrm{BMO}_c}^2\leq \left\| l_{\varphi} \right\| _{(\mathcal{H}%
_c^1)^{*}}^2\leq cN({\lambda} _{\varphi} )$ such that 
\[
\tau \int_{-\infty }^{+\infty }{\varphi} ^{*}fdt=\tau \int_{-\infty
}^{+\infty }{\varphi} ^{\prime *}fdt. 
\]
Thus ${\varphi} ={\varphi} ^{\prime }$ and ${\varphi} \in $\textrm{BMO}$_c(%
\Bbb{R},\mathcal{M})\,$with $\left\| {\varphi} \right\| _{\mathrm{BMO}%
_c}^2\leq cN({\lambda} _{\varphi} ).$ The converse had been already proved
in Lemma 1.4.\qed

\begin{corollary}
For $f\in \mathcal{H}_c^1(\Bbb{R},\mathcal{M}),$ we have 
\[
c^{-1}\left\| G_c(f)\right\| _1\leq \left\| S_c(f)\right\| _1\leq c\left\|
G_c(f)\right\| _1 
\]
\end{corollary}

\noindent\textbf{Proof. }By Theorem 2.4 and the first part of its proof, we
have 
\[
\left\| S_c(f)\right\| _1=\left\| f\right\| _{\mathcal{H}_c^1}\leq
c\sup_{\left\| {\varphi }\right\| _{\mathrm{BMO}_c}=1}\left| \tau \int f{%
\varphi }^{*}dt\right| \leq c\left\| G_c(f)\right\| _1^{\frac 12}\left\|
S_c(f)\right\| _1^{\frac 12} 
\]
Therefore 
\[
\left\| S_c(f)\right\| _1\leq c\left\| G_c(f)\right\| _1 
\]
The converse is an immediate consequence of Lemma 2.5.\qed
\medskip

\noindent\textbf{Remark. }The technique used in the proof of Lemma 2.5 is
classical (see \cite{st2}). The method to prove Theorem 2.4 is inspired by
the analogous one for martingales (see \cite{Gar}, \cite{Herz}, \cite{[PX2]}).

\section{The atomic decomposition of operator valued $\mathcal{H}^1$}

As in the classical case, the duality between $\mathcal{H}_c^1(\Bbb{R},%
\mathcal{M})$ and $\mathrm{BMO}_c(\Bbb{R},\mathcal{M})$ implies an atomic
decomposition of $\mathcal{H}_c^1(\Bbb{R},\mathcal{M}).$ The rest of this
chapter is devoted to this atomic decomposition. We say that a function $%
a\in L^1(\mathcal{M},L_c^2(\Bbb{R}))\,$is an $\mathcal{M}_c$-atom if

(i) $a$ is supported in a bounded interval $I;$

(ii) $\int_Iadt=0;$

(iii) $\tau (\int_I|a|^2dt)^{\frac 12}\leq |I|^{-\frac 12}.\smallskip $

Let $\mathcal{H}_c^{1,at}(\Bbb{R},\mathcal{M})$ be the space of all $f$
which admit a representation of the form 
\[
f=\sum_{i\in \Bbb{N}}{\lambda} _ia_i, 
\]
where the $a_i$'s are $\mathcal{M}_c$-atoms and ${\lambda} _i\in \Bbb{C}$
are such that $\sum_{i\in \Bbb{N}}|{\lambda} _i|<\infty .$ We equip $%
\mathcal{H}_c^{1,at}(\Bbb{R},\mathcal{M})$ with the following norm 
\[
\left\| f\right\| _{\mathcal{H}_c^{1,at}}=\inf \{\sum_{i\in \Bbb{N}}|{\lambda%
} _i|;f=\sum_{i\in \Bbb{N}}{\lambda} _ia_i;a_i\mbox{ are }\mathcal{M}_c%
\mbox{-\mbox{atoms}},{\lambda} _i\in \Bbb{C}\} 
\]
Similarly, we define $\mathcal{H}_r^{1,at}(\Bbb{R},\mathcal{M}).$ Then we
set 
\[
\mathcal{H}_{cr}^{1,at}(\Bbb{R},\mathcal{M})=\mathcal{H}_c^{1,at}(\Bbb{R},%
\mathcal{M})+\mathcal{H}_r^{1,at}(\Bbb{R},\mathcal{M}). 
\]

\begin{theorem}
$\mathcal{H}_c^{1,at}(\Bbb{R},\mathcal{M})=\mathcal{H}_c^1(\Bbb{R},\mathcal{M%
})$ with equivalent norms.
\end{theorem}

\noindent\textbf{Proof. }It is enough to prove $(\mathcal{H}_c^{1,at}(\Bbb{R}%
,\mathcal{M}))^{*}=\mathrm{BMO}_c(\Bbb{R},\mathcal{M}).$ Now, for any ${%
\varphi }\in \mathrm{BMO}_c(\Bbb{R},\mathcal{M})$ and $f\in \mathcal{H}%
_c^{1,at}(\Bbb{R},\mathcal{M})$ with $f=\sum_{i\in \Bbb{N}}{\lambda }_ia_i$
as above, by the Cauchy-Schwartz inequality we have 
\begin{eqnarray*}
|\tau \int {\varphi }^{*}fdt| &\leq &\sum_{i\in \Bbb{N}}|{\lambda }_i\tau
\int_{I_i}({\varphi }-{\varphi }_{I_i})^{*}a_idt| \\
&\leq &\sum_{i\in \Bbb{N}}|{\lambda }_i|\tau (\int_{I_i}|a_i|^2dt)^{\frac
12}\left\| (\int_{I_i}|{\varphi }-{\varphi }_{I_i}|^2dt)^{\frac 12}\right\|
_{\mathcal{M}} \\
&\leq &\left\| {\varphi }\right\| _{\mathrm{BMO}_c}\sum_{i\in \Bbb{N}}|{%
\lambda }_i|.
\end{eqnarray*}
Thus $\mathrm{BMO}_c(\Bbb{R},\mathcal{M})\subset (\mathcal{H}_c^{1,at}(\Bbb{R%
},\mathcal{M}))^{*}$ (a contractive inclusion). To prove the converse
inclusion, we denote by $L_0^1(\mathcal{M},L_c^2(I))$ the space of functions 
$f\in L^1(\mathcal{M},L_c^2(I))$ with $\int fdt=0.$ Notice that $L_0^1(%
\mathcal{M},L_c^2(I))\in \mathcal{H}_c^{1,at}(\Bbb{R},\mathcal{M})$ for
every bounded $I.$ Thus, every continuous functional $l$ on $\mathcal{H}%
_c^{1,at}(\Bbb{R},\mathcal{M})$ induces a continuous functional on $L_0^1(%
\mathcal{M},L_c^2(I))$ with norm smaller than $|I|^{\frac 12}\left\|
l\right\| _{(\mathcal{H}_c^{1,at})^{*}}.$ Consequently, we can choose a
sequence $({\varphi }_n)_{n\geq 1}$ satisfying the following conditions: 
\[
\begin{array}{l}
l(a)=\tau \int {\varphi }_n^{*}adt,\quad \forall \mathcal{M}_c\mbox{-}
\mbox{ atom
}a\mbox{ }\mbox{with}\ \mathrm{supp}\,a\subset (-n,n], \\ 
\left\| {\varphi }_n\right\| _{L^{^\infty }(\mathcal{M},L_c^2((-n,n]))}\leq c%
\sqrt{n}\left\| l\right\| _{(\mathcal{H}_c^{1,at})^{*}}; \\ 
{\varphi }_n\big|_{(-m,m]}={\varphi }_m,\ \forall n>m.
\end{array}
\]
Let ${\varphi }(t)={\varphi }_n(t),\forall t\in (-n,-n+1]\cup (n-1,n],n>0.$
We then have ${\varphi }\in L^\infty (\mathcal{M},L_c^2(\Bbb{R},\frac{dt}{%
1+t^2}))$ and 
\[
l(a)=\tau \int {\varphi }^{*}adt,\quad \forall \mathcal{M}_c\mbox{-}%
\mbox{
atom }a. 
\]
Considering $[g\otimes e]$ as defined in the remark after Lemma 2.2, by (\ref
{ge}) and (\ref{gf}) we have 
\begin{eqnarray*}
\left\| {\varphi }\right\| _{\mathrm{BMO}_c} &\leq &c\sup_{e\in H,\left\|
e\right\| _H=1}\sup_{\left\| g\right\| _{H^1(\Bbb{R},H)}=1}\left| \tau
\int_{-\infty }^{+\infty }{\varphi }^{*}[g\otimes e]dt\right| \\
&\leq &\sup_{\left\| f\right\| _{\mathcal{H}_c^{1,at}}=1}\left| \tau
\int_{-\infty }^{+\infty }{\varphi }^{*}fdt\right| \\
&=&\left\| l\right\| _{(\mathcal{H}_c^{1,at})^{*}}\ .\qed
\end{eqnarray*}

\begin{corollary}
$\mathcal{H}_r^{1,at}(\Bbb{R},\mathcal{M})=\mathcal{H}_r^1(\Bbb{R},\mathcal{M%
})$ and $\mathcal{H}_{cr}^{1,at}(\Bbb{R},\mathcal{M})=\mathcal{H}_{cr}^1(%
\Bbb{R},\mathcal{M})$ with equivalent norms. \medskip
\end{corollary}

\noindent\textbf{Remark. }The $\mathcal{M}$-atom considered in this section
is a non-commutative analogue of the classical $2$-atom for $H^1$ space. It
seems difficult to consider the non-commutative analogues of the classical $%
p- $atom for $p\neq 2.$

\medskip \noindent
\textbf{Remark.} We only considered the functions defined on $\Bbb{R}$ in
this chapter. However, one can check that all the proofs work well for the
functions defined on $\Bbb{R}^n$. And the analogous results can be proved
similarly for the functions defined on $\Bbb{T}^n$, where $\Bbb{T}$ is the
unit circle. Moreover, the relevant constants are independent
of $n$.

\chapter{The Maximal Inequality}

\section{The non-commutative Hardy-Littlewood maximal inequality}

\setcounter{theorem}{0} \setcounter{equation}{0} We recall the definition of
the non-commutative maximal function introduced in \cite{J} with an
inspiration from Pisier's non-commutative vector-valued space $L_p(N,\tau
;l_\infty )$ (see \cite{p2}). Let $0<p\leq \infty ,$ and let $(a_n)_{n\in 
\Bbb{Z}}$ be a sequence of elements in $L^p(\mathcal{M}).$ Set 
\begin{equation}
\left\| \sup_{n\in \Bbb{Z}}|a_n|\right\| _p=\inf_{a_n=ay_nb}\left\|
a\right\| _{2p}\left\| b\right\| _{2p}\sup_n\left\| y_n\right\| _{\mathcal{M}%
}  \label{pisier1}
\end{equation}
where the infimum is taken over all $a,b\in L_{2p}(\mathcal{M})$ and all
bounded sequences $(y_n)_{n\in \Bbb{Z}}\in \mathcal{M}$ such that $%
a_n=ay_nb. $ By convention, if $(a_n)_{n\in \Bbb{Z}}$ does not have such a
representation , we define $\left\| \sup_{n\in \Bbb{Z}}|a_n|\right\| _p$ as $%
+\infty .$ If $p>1$ and $(a_n)_{n\in \Bbb{Z}}$ is a sequence of positive
elements, it is proved by Junge (see \cite{J}, Remark 3.7) that(with $q$ the
index conjugate to $p$) 
\begin{equation}
\left\| \sup_{n\in \Bbb{Z}}|a_n|\right\| _p=\sup \left\{ \sum_{n\in \Bbb{Z}%
}\tau (a_nb_n):b_n\in L^q(\mathcal{M}),b_n\geq 0,\left\| \sum_{n\in \Bbb{Z}%
}b_n\right\| _q\leq 1\right\} .  \label{junge}
\end{equation}
Moreover, in this case, there exists a positive element $a\in L_{2p}(%
\mathcal{M})\,$and a sequence of positive elements $y_n$ such that $%
a_n=ay_na $ and $\left\| \sup_{n\in \Bbb{Z}}|a_n|\right\| _p=\left\|
a\right\| _{2p}^2\sup_n\left\| y_n\right\| _{\mathcal{M}}.$

We define similarly $\left\| \sup_{{\lambda }\in {\Lambda }}|a({\lambda }%
)|\right\| _p$ if ${\Lambda }$ is a countable set. If ${\Lambda }$ is
uncountable we set 
\begin{equation}
\left\| \sup_{{\lambda }\in {\Lambda }}|a({\lambda })|\right\| _p=\sup_{({%
\lambda }_n)_{n\in \Bbb{Z}}\in {\Lambda }}\left\| \sup_{n\in \Bbb{Z}}|a({%
\lambda }_n)|\right\| _p.  \label{uncount}
\end{equation}
Note that $\sup_{_\lambda }a_{_\lambda }$ does not make any sense in the
non-commutative setting and $\left\| \sup_{{\lambda }\in {\Lambda }}|a({%
\lambda })|\right\| _p$ is just a notation. Also note that 
\begin{equation}
\left\| \sup_{{\lambda }\in {\Lambda }}|a({\lambda })|\right\| _\infty
=\sup_{{\lambda }\in {\Lambda }}\left\| a({\lambda })\right\| _\infty .
\label{sup}
\end{equation}
To put the proceding definitions in proper perspective, 
we recall the following
identities satisfied by the norm of an $l_\infty(\Lambda)$-valued function
$a:\Bbb{R}\rightarrow l_\infty(\Lambda)$ in the classical space 
$L^p(\Bbb{R}, l_\infty(\Lambda))$ for an arbitrary index set $\Lambda$.

(a) 
\[
\left\| \sup_{{\lambda }\in {\Lambda }}|a({\lambda })|\right\| _p=\sup_{
J\subset {\Lambda } \mbox {finite}}\left\| \sup_{n\in J}|a({
\lambda }_n)|\right\| _p 
\]

(b) If $\left\| \sup_{{\lambda }\in {\Lambda }}|a({\lambda })|\right\|
_p<\infty ,$ then there exists $a\in L^p(\Bbb{R})$ such that 
$|a(\lambda )|\leq a, \forall \lambda\in \Lambda$
and 
\[
\left\| a\right\| _p=\left\| \sup_{{\lambda }\in {\Lambda }}|a({\lambda }%
)|\right\| _p. 
\]

The main result of this chapter is the non-commutative Hardy-Littlewood
maximal inequality. We will reduce it to the non-commutative Doob maximal
inequality for martingales already established by M. Junge \cite{H}. To this
end, we need to introduce two increasing filtration of dyadic $\sigma -$%
algebras on $\Bbb{R}.$ The key property of these $\sigma -$algebras is that
any interval of $\Bbb{R}$ is contained in an atom belonging to one of these $%
\sigma -$algebras with a comparable size (see Proposition 3.1 below). This
approach is also useful for other problems. We will see, for instance, that $%
\mathrm{BMO}_c(\Bbb{R},\mathcal{M})$ can be written as the intersection of
two dyadic \textrm{BMO }spaces. This approach is extremely simple and seems
new even in the scalar case (i.e. the classical case), see \cite{m}.

The two increasing filtrations of dyadic $\sigma -$algebras $\mathcal{D=}\{%
\mathcal{D}_n\}_{n\in \Bbb{Z}},\mathcal{D}^{\prime }=\{\mathcal{D}_n^{\prime
}\}_{n\in \Bbb{Z}}$ that we will need are defined as follows: The first one, 
$\mathcal{D=}\{\mathcal{D}_n\}_{n\in \Bbb{Z}}$, is simply the usual dyadic
filtration, that is, $\mathcal{\ D}_n$ is the $\sigma -$algebra generated by
the atoms 
\[
D_n^k=(k2^{-n},(k+1)2^{-n}];\quad k\in \Bbb{Z}. 
\]
The definition of $\mathcal{D}^{\prime }=\{\mathcal{D}_n^{\prime }\}_{n\in 
\Bbb{Z}}$ is a little more complicated. For an even integer $n,$ the atoms
of $\mathcal{D}_n^{\prime }$ are given by 
\[
D_n^{{\prime }\,k}=((k+\frac 13)2^{-n},(k+\frac 43)2^{-n}],\quad k\in \Bbb{Z}%
; 
\]
while for an odd integer $n,$ $\mathcal{D}_n^{\prime }$ is generated by the
atoms 
\[
D_n^{{\prime }\,k}=((k+\frac 23)2^{-n},(k+\frac 53)2^{-n}],\quad k\in \Bbb{Z}%
. 
\]
It is easy to see that $\mathcal{D}^{\prime }=\{\mathcal{D}_n^{\prime
}\}_{n\in \Bbb{Z}}$ is indeed an increasing filtration.\medskip

The following simple observation is the key of our approach.

\begin{proposition}
For any interval $I\subset \Bbb{R},$ there exist $k_I,N\in \Bbb{Z}$ such
that $I\subset $ $D_N^{k_I}$ and $|D_N^{k_I}|\leq 6|I|$ or $I\subset $ $%
D_N^{^{\prime }k_I}$ and $|D_N^{^{\prime }k_I}|\leq 6|I|,$ the constant $N$
only depends on the length of $I.$
\end{proposition}

\noindent\textbf{Proof.} $\,$To see this, choose $N\in \Bbb{Z}$ such that $%
\frac{2^{-N-1}}3\leq $ $|I|<\frac{2^{-N}}3.$ Denote 
\[
A_N=\{(k2^{-N});k\in \Bbb{Z}\},\quad A_N^{\prime }=\{((k+\frac
13)2^{-N},(k+\frac 23)2^{-N});k\in \Bbb{Z\}}. 
\]
Note that for any two points $a,b\in A_N\cup A_N^{^{\prime }},$ we have $%
|a-b|\geq \frac 132^{-N}>|I|.$ Thus there is no more than one element of $%
A_N\cup A_N^{^{\prime }}$ in $I.$ Then $I\cap A_N=\phi $ or $I\cap
A_N^{^{\prime }}=\phi .$ Therefore, $I$ must be contained in some $D_N^{k_I}$
or $D_N^{^{\prime }\ k_I}.$\qed \medskip \newline
\noindent\textbf{Remark.} See \cite{m} for a generalization of Proposition
3.1.\medskip \newline
\noindent\textbf{Remark.} If an $\mathcal{M}_c$-atom defined in Chapter 2
admits its supporting interval as $D_N^k$ (resp. $D_N^{\prime }{}^k$) for
some $k,N\in \Bbb{Z},$ we call it $\mathcal{M}_c$-$\mathcal{D}$-atom (resp. $%
\mathcal{M}_c$-$\mathcal{D}^{\prime }$-atom). Proposition 3.1 implies that
an $\mathcal{M}_c$-atom is either an $\mathcal{M}_c$-$\mathcal{D}$-atom or
an $\mathcal{M}_c$--$\mathcal{D}^{\prime }$-atom up to a fixed factor.
Therefore the atomic Hardy space $\mathcal{H}_c^{1,at}(\Bbb{R},\mathcal{M})$
defined in Chapter 2 can be characterized only by $\mathcal{M}_c$-$\mathcal{D%
}$-atoms and $\mathcal{M}_c$-$\mathcal{D}^{\prime }$-atoms. A similar remark
applies to the atomic row Hardy space $\mathcal{H}_r^{1,at}(\Bbb{R},\mathcal{%
M})$. See Chapter 5 for more results of this type.

\medskip

The proof of the following Proposition (as well as that of Theorem 3.3)
illustrates well our approach to reduce problems on functions to those on
martingales. Put 
\[
f_h(t)=\frac 1{h_1+h_2}\int_{t-h_1}^{t+h_2}f(x)dx,\ \forall h=(h_1,h_2)\in 
\Bbb{R}^{+}\Bbb{\times R}^{+}. 
\]

\begin{proposition}
Let $(a_n)_{n\in \Bbb{Z}}$ be a positive sequence in $L^p(L^\infty (\Bbb{R}%
)\otimes \mathcal{M})$ and $h_n=(h_{n,1},h_{n,2})\in \Bbb{R}^{+}\Bbb{\times R%
}^{+},n\in \Bbb{Z}.$

(i)\quad If $1\leq p<\infty $, 
\begin{equation}
\left\| \sum_{n\in \Bbb{Z}}(a_n)_{h_n}\right\| _{L^p(L^\infty (\Bbb{R}%
)\otimes \mathcal{M})}\leq c_p\left\| \sum_{n\in \Bbb{Z}}a_n\right\|
_{L^p(L^\infty (\Bbb{R})\otimes \mathcal{M})}.  \label{stein}
\end{equation}

(ii)\quad If $1<p\leq \infty $, 
\begin{equation}
\left\| \sup_{n\in \Bbb{Z}}|(a_n)_{h_n}|\right\| _{L^p(L^\infty (\Bbb{R}%
)\otimes \mathcal{M})}\leq c_p\left\| \sup_{n\in \Bbb{Z}}|a_n|\right\|
_{L^p(L^\infty (\Bbb{R})\otimes \mathcal{M})}.  \label{max2}
\end{equation}
\end{proposition}

\noindent\textbf{Proof.} From Proposition 3.1, $\forall n\in \Bbb{Z},$ for
every $t\in \Bbb{R},$ there exist some $k_t,N_n\in \Bbb{Z}$ such that $%
(t-h_{n,1},t+h_{n,2})$ is contained in $D_{N_n}^{k_t}$ or $D_{N_n}^{\prime
\,k_t}$ and 
\[
|D_{N_n}^{k_t}|=|D_{N_n}^{{\prime }\,k_t}|\leq 6(h_{n,1}+h_{n,2}). 
\]
Thus 
\begin{equation}
(a_n)_{h_n}\leq 6(E(a_n|\mathcal{D}_{N_{_n}})+E(a_n|\mathcal{D}%
_{N_n}^{\prime })),\quad \forall n\in \Bbb{Z},  \label{anhn}
\end{equation}
where $E(\cdot \;|\mathcal{D}_{N_{_n}})($resp. $E(\cdot \;|\mathcal{D}%
_{N_n}^{\prime }))\,$denotes the conditional expectation with respect to $%
\mathcal{D}_{N_{_n}}($resp. $\mathcal{D}_{N_n}^{\prime }).$ Then (\ref{stein}%
)$\,$follows from Theorem 0.1 of \cite{J}. By (\ref{junge}) and (\ref{stein}%
), 
\begin{eqnarray*}
&&\left\| \sup_{n\in \Bbb{Z}}|(a_n)_{h_n}|\right\| _{L^p(L^\infty (\Bbb{R}%
)\otimes \mathcal{M})} \\
&=&\sup \{\sum_{n\in \Bbb{Z}}\tau \int_{\Bbb{R}}\frac 1{h_{n,1}+h_{n,2}}{%
\int }_{t-h_{n,1}}^{t+h_{n,2}}a_n(x)dxb_n(t)dt:\left\| \sum_{n\in \Bbb{Z}%
}b_n\right\| _{L^q(L^\infty (\Bbb{R})\otimes \mathcal{M})}\leq 1\} \\
&=&\sup \{\sum_{n\in \Bbb{Z}}\tau \int_{\Bbb{R}}\frac 1{h_{n,1}+h_{n,2}}{%
\int }_{x-h_{n,2}}^{x+h_{n,1}}b_n(t)dta_n(x)dx:\left\| \sum_{n\in \Bbb{Z}%
}b_n\right\| _{L^q(L^\infty (\Bbb{R})\otimes \mathcal{M})}\leq 1\} \\
&\leq &\sup \{\sum_{n\in \Bbb{Z}}\tau \int_{\Bbb{R}}b_n(x)a_n(x)dx:\left\|
\sum_{n\in \Bbb{Z}}b_n\right\| _{L^q(L^\infty (\Bbb{R})\otimes \mathcal{M}%
)}\leq c_p\} \\
&\leq &c_p\left\| \sup_{n\in \Bbb{Z}}|a_n|\right\| _{L^p(L^\infty (\Bbb{R}%
)\otimes \mathcal{M})}
\end{eqnarray*}
This is (\ref{max2}).\qed

\medskip

The following is our non-commutative Hardy-Littlewood maximal inequality.
Denote by $\mathcal{P}(\mathcal{M})$ the family of all projections of a von
Neumann algebra $\mathcal{M}.$

\begin{theorem}
(i)\quad Let $f\in L^1(L^\infty (\Bbb{R})\otimes \mathcal{M})$ and ${\lambda 
}>0.$ Then there exists $e^{\lambda }\in \mathcal{P}(L^\infty (\Bbb{R})\otimes 
\mathcal{M})\,$such that 
\begin{equation}
\sup_{h\in \Bbb{R}^{+}\Bbb{\times R}^{+}}\left\| e^{\lambda }f_he^{\lambda}
\right\| _{L^{^\infty }(\Bbb{R})\otimes \mathcal{M}}\leq {\lambda },\quad
\left[ \tau \otimes \int \right] (1-e^{\lambda })<\frac{c_1\left\| f\right\| _1%
}\lambda .  \label{11}
\end{equation}

(ii)\quad Let $1<p\leq \infty $ and $f\in L^p(L^\infty (\Bbb{R})\otimes 
\mathcal{M}).$ Then 
\begin{equation}
\left\| \sup_{h\in \Bbb{R}^{+}\Bbb{\times R}^{+}}|f_h|\right\| _p\leq
c_p\left\| f\right\| _p.  \label{MAXC}
\end{equation}
Moreover, for every positive $f\in L^p(L^\infty (\Bbb{R})\otimes 
\mathcal{M}),$ there exists a positive $F\in L^p(L^\infty (\Bbb{R})\otimes 
\mathcal{M})$ such that $f_h\leq F$ for all $h$ and 
\begin{equation}
\left\| F\right\| _p\leq
c_p\left\| f\right\| _p.  \label{MAXCp}
\end{equation}
\end{theorem}

\noindent\textbf{Proof.} By decomposing $f=f_1-f_2+i(f_3-f_4)$ with positive 
$f_k$, we can assume $f$ positive. To prove (i), for given $f,{\lambda }%
,(h_n)_{n\in \Bbb{Z}}\in \Bbb{R}^{+}\Bbb{\times R}^{+},$ let $\mathcal{D}%
_{N_n},\mathcal{D}_{N_n}^{\prime }$be as in the proof of Proposition 3.2. By
the weak type (1,1) inequality of non-commutative martingales in \cite{C} we
have $\forall {\lambda }>0,$ $\exists e^{\lambda },e^{\prime {\lambda }}\in 
\mathcal{P}(L^\infty (\Bbb{R})\otimes \mathcal{M})$ $\,$such that 
\[
\sup_n\left\| e^{\lambda }E(f|\mathcal{D}_{N_{_n}})e^{\lambda }\right\|
_{L^{^\infty }(\Bbb{R})\otimes \mathcal{M}}\leq \frac \lambda {12},\quad
\tau \otimes \int (1-e^{\lambda })<\frac{c\left\| f\right\| _1}\lambda 
\]
and 
\[
\sup_n\left\| e^{\lambda }E(f|\mathcal{D}_{N_{_n}}^{\prime })e^{\lambda}
\right\| _{L^{^\infty }(\Bbb{R})\otimes \mathcal{M}}\leq \frac \lambda
{12},\quad \tau \otimes \int (1-e^{\prime \,{\lambda }})<\frac{c\left\|
f\right\| _1}\lambda 
\]
for every $f\in L^1(L^\infty (\Bbb{R})\otimes \mathcal{M})$ and $(h_n)_{n\in 
\Bbb{Z}}\in \Bbb{R}^{+}\Bbb{\times R}^{+}.$ Let $\widetilde{e^{\lambda }}%
=e^{\lambda }\wedge e^{\prime \,{\lambda }},$ then 
\[
\tau \otimes \int (1-\widetilde{e^{\lambda }})<\frac{2c\left\| f\right\| _1}%
\lambda . 
\]
By Proposition 3.1, we have 
\[
\widetilde{e^{\lambda }}f_{h_n}\widetilde{e^{\lambda }}\leq 6(e^{\lambda }E(f|%
\mathcal{D}_{N_n})e^{\lambda }+e^{{\prime }{\lambda }}E(f|\mathcal{D}%
_{N_{h_n}}^{\prime })e^{{\prime }{\lambda }}). 
\]
Therefore, 
\begin{eqnarray*}
&&\sup_{h\in \Bbb{R}^{+}\Bbb{\times R}^{+}}\left\| \widetilde{e^{\lambda }}f_h%
\widetilde{e^{\lambda }}\right\| _{L^\infty (\Bbb{R})\otimes \mathcal{M}} \\
&=&\sup_{(h_n)_{n\in \Bbb{Z}}}\sup_n\left\| \widetilde{e^{\lambda }}f_{h_n}%
\widetilde{e^{\lambda }}\right\| _{L^\infty (\Bbb{R})\otimes \mathcal{M}} \\
&\leq &6\sup_n\left\| e^{^{\prime }{\lambda }}E(f|\mathcal{D}_{N_n}^{\prime
})e^{{\prime }{\lambda }}\right\| _{L^\infty (\Bbb{R})\otimes \mathcal{M}%
}+6\sup_n\left\| e^{\lambda }E(f|\mathcal{D}_{N_n})e^{\lambda }\right\|
_{L^\infty (\Bbb{R})\otimes \mathcal{M}} \\
&\leq &{\lambda }.
\end{eqnarray*}
This is (\ref{11}).\ To prove (ii), by (\ref{junge}) and (\ref{stein}), 
\begin{eqnarray*}
&&\left\| \sup_{n\in \Bbb{Z}}|f_{h_n}|\right\| _{L^p(L^\infty (\Bbb{R}%
)\otimes \mathcal{M})} \\
&=&\sup \{\sum_{n\in \Bbb{Z}}\tau \int_{\Bbb{R}}\frac 1{h_{n,1}+h_{n,2}}{%
\int }_{t-h_{n,1}}^{t+h_{n,2}}f(x)dxb_n(t)dt:\left\| \sum_{n\in \Bbb{Z}%
}b_n\right\| _{L^q(L^\infty (\Bbb{R})\otimes \mathcal{M})}\leq 1\} \\
&=&\sup \{\sum_{n\in \Bbb{Z}}\tau \int_{\Bbb{R}}\frac 1{h_{n,1}+h_{n,2}}{%
\int }_{x-h_{n,2}}^{x+h_{n,1}}b_n(t)dtf(x)dx:\left\| \sum_{n\in \Bbb{Z}%
}b_n\right\| _{L^q(L^\infty (\Bbb{R})\otimes \mathcal{M})}\leq 1\} \\
&\leq &\sup \{\tau \int_{\Bbb{R}}f(x)\sum_{n\in \Bbb{Z}}a_n(x)dx:\left\|
a_n\right\| _{L^q(L^\infty (\Bbb{R})\otimes \mathcal{M})}\leq c_p\} \\
&\leq &c_p\left\| f\right\| _{L^p(L^\infty (\Bbb{R})\otimes \mathcal{M})}
\end{eqnarray*}
where, we set $a_n=\frac 1{h_{n,1}+h_{n,2}}{\int }%
_{x-h_{n,2}}^{x+h_{n,1}}b_n(t).$ This yields (\ref{MAXC}).
The inequality (\ref{MAXCp}) follows from Theorem 0.2 of \cite{J} because of (\ref{anhn}). 
\qed

Using standard arguments and Theorem 3.3 we can easily obtain the
non-commutative analogue of the classical non-tangential maximal inequality. Recall, 
as in Chapter 1, we also use $f$ to denote its Poisson integral on the upper half plane.

\begin{theorem}
(i)\quad Let $f\in L^1(L^\infty (\Bbb{R})\otimes \mathcal{M}).$ Then $%
\forall {\lambda }>0,\exists $ $e^{\lambda }\in \mathcal{P}(L^\infty (\Bbb{R}%
)\otimes \mathcal{M}),$ such that 
\begin{equation}
\sup_{(t,y)\in \Gamma }\left\| e^{\lambda }f(x+t,y)e^{\lambda }\right\|
_{L^\infty (\Bbb{R})\otimes \mathcal{M}}\leq {\lambda },\quad \tau \otimes
\int (1-e^{\lambda })<\frac{c_1\left\| f\right\| _1}\lambda ,\forall {\lambda }%
>0  \label{p11}
\end{equation}
(ii)\quad Let $f\in L^p(L^\infty (\Bbb{R})\otimes \mathcal{M}),1<p\leq
\infty .$ Then 
\begin{equation}
\left\| \sup_{(t,y)\in \Gamma }|f(x+t,y)|\right\| _p\leq c_p\left\|
f\right\| _p.\quad  \label{pmax}
\end{equation}
Moreover, for every positive $f\in L^p(L^\infty (\Bbb{R})\otimes 
\mathcal{M}),$ there exists a positive $F\in L^p(L^\infty (\Bbb{R})\otimes 
\mathcal{M})$ such that $f(\cdot+t,y)\leq F$ for all $(t,y)\in \Gamma$ and 
\begin{equation}
\left\| F\right\| _p\leq
c_p\left\| f\right\| _p.  \label{pmaxp}
\end{equation}
\end{theorem}

\noindent\textbf{Proof.} Notice that 
\[
P_y(x)=\frac 1\pi \frac y{x^2+y^2}\leq \frac 1\pi \frac 1{2^{2(k-1)}y+y},\
\forall 2^{k-1}y\leq |x|. 
\]
We have, for every positive $f$ and any $(t,y)\in \Gamma ,$ 
\begin{eqnarray}
&&f(x+t,y)  \nonumber \\
&=&\int_{\Bbb{R}}f(s)P_y(x+t-s)ds  \nonumber \\
&\leq &\frac 1\pi \int_{|x+t-s|\leq y}f(s)\frac 1yds+\frac 1\pi
\sum_{k=1}^\infty \int_{2^{k-1}y\leq |x+t-s|\leq 2^ky}f(s)\frac
1{2^{2(k-1)}y+y}ds  \nonumber \\
&\leq &\frac 1\pi \sum_{k=0}^\infty \frac 8{2^k}\frac
1{2^{k+1}y}\int_{|x+t-s|\leq 2^ky}f(s)ds.  \label{Ff}
\end{eqnarray}
Considering $h_{k,y}=(2^ky-t,2^ky+t)\in \Bbb{R}^{+}\times \Bbb{R}^{+},$
we get (\ref {pmaxp}) from (\ref{MAXCp}). And by (\ref{MAXC}), 
\begin{eqnarray*}
\left\| \sup_{(t,y)\in \Gamma }|f(x+t,y)|\right\| _p &\leq &\frac 1\pi
\sum_{k=0}^\infty \frac 8{2^k}\left\| \sup_{h_{k,y}}|f_{h_{k,y}}|\right\| _p
\\
&\leq &c_p\left\| f\right\| _p.
\end{eqnarray*}
Decomposing $f=f_1-f_2+i(f_3-f_4)$ with positive $f_k,$ we get (\ref{pmax})
for all $f\in L^p(L^\infty (\Bbb{R})\otimes \mathcal{M}).$ We can prove (\ref
{p11}) similarly.\qed \medskip

\section{The non-commutative Lebesgue differentiation theorem and
\\ non-tangential limit of Poisson integrals}

We end this chapter with the non-commutative Lebesgue differentiation
theorem and non-tangential limit of Poisson integrals. These are
consequences of Theorem 3.3 and Theorem 3.4. To this end, we first need to
recall the non-commutative version of the almost everywhere convergence. Let 
$(f_{\lambda} )_{{\lambda} \in {\lambda} }$ be a family of elements in $L^p(%
\mathcal{M},\tau ).\,$We say $(f_{\lambda} )_{{\lambda} \in {\lambda} }\,$%
converges to $f$ almost uniformly, abbreviated as $f_{\lambda} \stackrel{a.u%
}{\rightarrow }f,$ if for every ${\varepsilon} >0,$ there exists $%
e_{\varepsilon} \in \mathcal{P}(\mathcal{M})\,$such that $\tau
(1-e_{\varepsilon} )<{\varepsilon} \,$and 
\[
\lim_{{\lambda} \rightarrow {\lambda} _0}\left\| e_{\varepsilon}
(f_{\lambda} -f)\right\| _\infty =0. 
\]
Moreover, we say $(f_{\lambda} )_{{\lambda} \in {\lambda} }$ converges to $%
f\,$ bilaterally almost uniformly, abbreviated as $f_{\lambda} \stackrel{%
b.a.u}{\rightarrow }f,\,$if for every ${\varepsilon} >0,$ there exists $%
e_{\varepsilon} \in \mathcal{P}(\mathcal{M})\,$such that $\tau
(1-e_{\varepsilon} )<{\varepsilon} \,$and 
\[
\lim_{{\lambda} \rightarrow {\lambda} _0}\left\| e_{\varepsilon}
(f_{\lambda} -f)e_{\varepsilon} \right\| _\infty =0. 
\]
Obviously, $f_{\lambda} \stackrel{a.u}{\rightarrow }f$ implies $f_{\lambda} 
\stackrel{b.a.u}{\rightarrow }f.$

Recall that the map $x\mapsto x^p$ $(1\leq p\leq 2)$ is convex on the
positive cone $\mathcal{M}_{+}$ of $\mathcal{M}$ (see \cite{B}). Thus, for $%
f\in L^p(L^\infty (\Bbb{R})\otimes \mathcal{M})$ $(1\leq p\leq 2),$ we get 
\begin{equation}
\int_A|f|dt\leq (\int_A|f|^pdt)^{\frac 1p},\quad \forall A\subseteq \Bbb{R}
,~|A|=1. \label{p<2}
\end{equation}
Note that for any $x,y\in \mathcal{M}_+,$ $x\leq y$ implies $x^{q}\leq y^{q}, 
\forall 0<q\leq 1$.
Using (\ref{p<2}) successively, we get the following Lemma.

\begin{lemma}
For $f\in L^p(L^\infty (\Bbb{R})\otimes \mathcal{M}),\,1\leq p<\infty ,$%
\begin{equation}
\int_A|f|dt\leq (\int_A|f|^pdt)^{\frac 1p},\quad \forall A\subseteq \Bbb{R}%
,~|A|=1.  \label{fp}
\end{equation}
\end{lemma}

And recall that for any bounded linear operators $a,b$ on a Hilbert space $%
H, $ $a$ positive and $\left\| b\right\| \leq 1,$ if $T$ is an operator
monotone function defined for positive operators (for example, $%
T(a)=a^{\frac 1p},p\geq 1$) then 
\begin{equation}
b^{*}T(a)b\leq T(b^{*}ab).
\end{equation}
This is the so-called Hansen's inequality (see \cite{H}). In
particular, we have 
\begin{equation}
b^{*}ab\leq (b^{*}a^pb)^{\frac 1p}.  \label{hansen}
\end{equation}

\begin{theorem}
(i)\quad Let $1\leq p<2.$ We have $f_h\stackrel{b.a.u}{\rightarrow }f$ as $%
h\rightarrow 0$ for any $f\in L^p(L^\infty (\Bbb{R})\otimes \mathcal{M}).$

(ii)\quad Let $2\leq p<\infty .$ We have $f_h\stackrel{a.u}{\rightarrow }f$
as $h\rightarrow 0$ for any $f\in L^p(L^\infty (\Bbb{R})\otimes \mathcal{M}%
). $
\end{theorem}

\noindent\textbf{Proof. }(i) Without loss of generality, we can assume $f$
selfadjoint. For any given $f\in L^p(L^\infty (\Bbb{R})\otimes \mathcal{M})$
and ${\varepsilon }>0,$ choose $f^n=\sum_{k=1}^{N_n}{\varphi }_kx_k,\,$
where $x_k\in S_{\mathcal{M}}^{+}$ and where ${\varphi }_k:\Bbb{R}%
\rightarrow \Bbb{C}$ are continuous functions with compact support, such
that 
\begin{equation}
\left\| |f-f^n|^p\right\| _1=\left\| f-f^n\right\| _p^p<(\frac
1{2^n})^p\frac \varepsilon {2^n}.  \label{ffn}
\end{equation}
$\,$Choose $e_{1,n}^{\varepsilon} \in \mathcal{P}(L^\infty (\Bbb{R})\otimes 
\mathcal{M})$ such that 
\[
\tau \otimes \int (1-e_{1,n}^{\varepsilon} )<\frac \varepsilon {2^n}\,\quad %
\mbox{and\quad }\left\| e_{1,n}^{\varepsilon} |f^n-f|^pe_{1,n}^{\varepsilon}
\right\| _{L^\infty (\Bbb{R})\otimes \mathcal{M}}<(\frac 1{2^n})^p.\, 
\]
Set $e_1^{\varepsilon} =\wedge _ne_{1,n}^{\varepsilon} .$ We have $\tau \otimes
\int (1-e_1^{\varepsilon} )<{\varepsilon }$ and by (\ref{hansen}), 
\begin{eqnarray}
\left\| e_1^{\varepsilon} (f^n-f)e_1^{\varepsilon} \right\| _{L^\infty (\Bbb{R}%
)\otimes \mathcal{M}} &\leq &\left\| e_1^{\varepsilon} |f^n-f|e_1^{\varepsilon}
\right\| _{L^\infty (\Bbb{R})\otimes \mathcal{M}}  \nonumber \\
&\leq &\left\| e_1^{\varepsilon} |f^n-f|^pe_1^{\varepsilon} \right\| _{L^\infty (%
\Bbb{R})\otimes \mathcal{M}}^{\frac 1p}  \nonumber \\
&<&\frac 1{2^n},~~~\forall n\geq 1.  \label{e1}
\end{eqnarray}
On the other hand, by (\ref{11}) and (\ref{ffn}) we can find a sequence $%
(e_{2,n}^{\varepsilon} )_{n\geq 0}\subset \mathcal{P}(L^\infty (\Bbb{R}%
)\otimes \mathcal{M})\,$such that 
\begin{eqnarray}
\tau \otimes \int (1-e_{2,n}^{\varepsilon} ) &<&\frac \varepsilon {2^n}~ 
\nonumber \\
~\left\| e_{2,n}^{\varepsilon} (|f^n-f|^p)_he_{2,n}^{\varepsilon} \right\|
_{L^\infty (\Bbb{R})\otimes \mathcal{M}} &<&(\frac 1{2^n})^p,~\ \forall h\in 
\Bbb{R}^{+}\Bbb{\times R}^{+}.\,  \label{2n}
\end{eqnarray}
Set $e_2^{\varepsilon} =\wedge _ne_{2,n}^{\varepsilon} ,$ we have $\tau \otimes
\int (1-e_2^{\varepsilon} )<{\varepsilon }.$ By (\ref{fp}), (\ref{hansen}) and
(\ref{2n}) 
\begin{eqnarray}
\left\| e_2^{\varepsilon} (f_h^n-f_h)e_2^{\varepsilon} \right\| _{L^\infty (\Bbb{%
R})\otimes \mathcal{M}} &\leq &\left\| e_{2,n}^{\varepsilon}
(|f^n-f|)_he_{2,n}^{\varepsilon} \right\| _{L^\infty (\Bbb{R})\otimes \mathcal{%
M}}  \nonumber \\
&\leq &\left\| e_{2,n}^{\varepsilon} (|f^n-f|^p)_h^{\frac
1p}e_{2,n}^{\varepsilon} \right\| _{L^\infty (\Bbb{R})\otimes \mathcal{M}} 
\nonumber \\
&\leq &(\left\| e_{2,n}^{\varepsilon} (|f^n-f|^p)_he_{2,n}^{\varepsilon}
\right\| _{L^\infty (\Bbb{R})\otimes \mathcal{M}})^{\frac 1p}  \nonumber \\
&<&\frac 1{2^n},\qquad \forall n\geq 0,h\in \Bbb{R}^{+}\Bbb{\times R}^{+}.
\label{enk}
\end{eqnarray}
Recall that by the classical Lebesgue differentiation theorem, 
\[
\lim_{h\rightarrow 0}\left\| {\varphi }_h-{\varphi }\right\| _\infty =0 
\]
if ${\varphi }:\Bbb{R}\rightarrow \Bbb{C\,}$is continuous with compact
support. Then by the choice of $f_{n\mbox{ }}$we deduce that 
\[
\lim_{h\rightarrow 0}\left\| f_h^n-f^n\right\| _{L^\infty (\Bbb{R})\otimes 
\mathcal{M}}=0,\forall n\geq 1. 
\]
Let $e^{\varepsilon} =e_1^{\varepsilon} \wedge e_2^{\varepsilon} ,$ then $\tau
\otimes \int (1-e^{\varepsilon} )<2{\varepsilon }.$ For any $n>0,$ choose $%
S_n>0$ such that $\left\| f_h^n-f^n\right\| _\infty <\frac 1{2^n}$ for any $%
h\in \Bbb{R}^{+}\Bbb{\times R}^{+}\ $such that $h_1+h_2<S_n.$ Then, for any $%
h\in \Bbb{R}^{+}\Bbb{\times R}^{+}\ $such that $h_1+h_2<S_n,$ 
\begin{eqnarray*}
\left\| e^{\varepsilon} (f_h-f)e^{\varepsilon} \right\| _\infty &\leq &\left\|
e^{\varepsilon} (f^n-f)e^{\varepsilon} \right\| _\infty +\left\|
f_h^n-f^n\right\| _\infty +\left\| e^{\varepsilon} (f_h^n-f_h)e^{\varepsilon}
\right\| _\infty \\
&\leq &\left\| e_1^{\varepsilon} (f^n-f)e_1^{\varepsilon} \right\| _\infty
+\left\| f_h^n-f^n\right\| _\infty +\left\| e_2^{\varepsilon}
(f_h^n-f_h)e_2^{\varepsilon} \right\| _\infty \\
&\leq &\frac 3{2^n}.
\end{eqnarray*}
Thus $\lim_{h\rightarrow 0}\left\| e^{\varepsilon} (f_h-f)e^{\varepsilon}
\right\| _\infty \rightarrow 0.\,$This completes the proof of (i).

(ii) The proof of (i) works well for the part (ii) of the theorem with some
minor changes. Let $(f^n)_{n\in \Bbb{N}}$ and $e_1^{\varepsilon}
,e_2^{\varepsilon} ,e^{\varepsilon} $ be as above. Since $p\geq 2,$ instead
of (\ref{e1}), (\ref{enk}), $\,$by (\ref{fp}) and (\ref{hansen}) we have 
\begin{equation}
\left\| e_1^{\varepsilon} (f^n-f)\right\| _\infty =\left\| e_1^{\varepsilon}
|f^n-f|^2e_1^{\varepsilon} \right\| _\infty ^{\frac 12}\leq \left\|
e_1^{\varepsilon} |f^n-f|^pe_1^{\varepsilon} \right\| _\infty ^{\frac
1p}<\frac 1{2^n},~\forall n\geq 1;  \label{e1'}
\end{equation}
and also 
\begin{eqnarray}
\left\| e_2^{\varepsilon} (f_h^n-f_h)\right\| _\infty &=&\left\|
e_2^{\varepsilon} |f_h^n-f_h|^2e_2^{\varepsilon} \right\| _\infty ^{\frac 12}
\nonumber \\
&\leq &(\left\| e_2^{\varepsilon} (|f^n-f|^2)_he_2^{\varepsilon} \right\|
_\infty )^{\frac 12}  \nonumber \\
&\leq &(\left\| e_2^{\varepsilon} (|f^n-f|^p)_he_2^{\varepsilon} \right\|
_\infty )^{\frac 1p}<\frac 1{2^n},\quad \forall n\geq 1.  \label{e2'}
\end{eqnarray}
Then we can conclude as in the proof of (i). $\qed$

\begin{theorem}
$\,$(i)\quad Let $1\leq p<2,f\in L^p(L^\infty (\Bbb{R})\otimes \mathcal{M}).$
We have $f(\cdot +u,y)\stackrel{b.a.u}{\rightarrow }f$ as $\Gamma \ni
(u,y)\rightarrow 0.$

(ii)\quad Let $2\leq p<\infty ,f\in L^p(L^\infty (\Bbb{R})\otimes \mathcal{M}%
).$ We have $f(\cdot +u,y)\stackrel{a.u}{\rightarrow }f$ as $\Gamma \ni
(u,y)\rightarrow 0.$
\end{theorem}

\noindent\textbf{Proof.} We can assume $f\geq 0$ by decomposing $f$ into
four positive parts. Given ${\varepsilon }>0,$ let $f^n,e_{i,n}^{\varepsilon}
,e_i^{\varepsilon} $ $(i=1,2)$ be as in the proof of Theorem 3.6. We use
the same notation $f^n$ for the Poisson integral of $f^n.$ It is easy to see
that 
\[
\lim_{(u,y)\rightarrow 0.}\left\| f^n(\cdot +u,y)-f^n\right\| _\infty
\rightarrow 0,\qquad \forall n\geq 0,\quad (u,y)\in \Gamma 
\]
Let $e^{\varepsilon} =e_1^{\varepsilon} \wedge e_2^{\varepsilon} .$ For any $n>0,$
choose $Y_n>0$ such that 
\[
\left\| f^n(\cdot +u,y)-f^n\right\| _\infty <\frac 1{2^n} 
\]
for any $(u,y)\in \Gamma ,|u|+y\leq Y_n.$ To prove (i), from (\ref{e1}), (%
\ref{enk}) we have, for any $(u,y)\in \Gamma ,|u|+y\leq Y_n,$ 
\begin{eqnarray*}
&&\left\| e^{\varepsilon} (f(\cdot +u,y)-f(\cdot ))e^{\varepsilon} \right\|
_\infty \\
&\leq &\left\| e^{\varepsilon} (f^n-f)e^{\varepsilon} \right\| _\infty +\left\|
f^n(\cdot +u,y)-f^n\right\| _\infty \\
&&+\left\| e^{\varepsilon} (\int_{\Bbb{R}}(f-f^n)(s)P_y(x+u-s)ds)e^{\varepsilon}
\right\| _\infty \\
&\leq &\frac 1{2^n}+\frac 1{2^n}+\sum_{k=0}^\infty \left\| e^{\varepsilon}
(\int_{|x+u-s|\leq 2^ky}|f-f^n|\frac 2{2^{2(k-1)}y+y}ds)e^{\varepsilon}
\right\| _\infty \\
&\leq &\frac 2{2^n}+\sum_{k=0}^\infty \frac 8{2^k}\left\| e_2^{\varepsilon}
(\frac 1{2^ky}\int_{|x+u-s|\leq 2^ky}|f-f^n|ds)e_2^{\varepsilon} \right\|
_\infty \\
&\leq &\frac 2{2^n}+\sum_{k=0}^\infty \frac 8{2^k}\left\| e_2^{\varepsilon}
(|f-f^n|)_{h_{k,y}}e_2^{\varepsilon} \right\| _\infty \\
&\leq &\frac 2{2^n}+\frac 8{2^n},
\end{eqnarray*}
where $h_{k,y}=(2^ky-t,2^ky+t)\in \Bbb{R}^{+}\times \Bbb{R}^{+}.$ Thus 
\[
\lim_{(u,y)\rightarrow 0}\left\| e^{\varepsilon} (f(\cdot
+ty,y)-f)e^{\varepsilon} \right\| _\infty =0,\forall {\varepsilon }>0, 
\]
and then $f(\cdot +u,y)\stackrel{b.a.u}{\rightarrow }f$ when $\Gamma \ni
(u,y)\rightarrow 0.$ This is (i). Using (\ref{e1'}) and (\ref{e2'}) instead
of (\ref{e1}) and (\ref{enk}), we can prove (ii) similarly.\qed

\medskip \noindent
\textbf{Remark.} All the results carried out in this chapter
can be generalized to the case of functions defined on 
$\Bbb{R}^n$ or $\Bbb{T}^n$. Unfortunately, the relevant constants there will depend on 
$n$ because the constant in Proposition 2.5 of \cite{m} depends on $n$.
This could be corrected if we could find a direct proof for Theorem 3.3.

\medskip
\noindent\textbf{Remark. } When $p=\infty ,$ the corresponding convergence
problems discussed in this section are still open.

\chapter{The Duality between $\mathcal{H}^p $ and $\mathbf{\mathrm{BMO}}%
^q,1<p<2.$}

\setcounter{theorem}{0} \setcounter{equation}{0} In this chapter, we
describe the dual of $\mathcal{H}_c^p(\Bbb{R},\mathcal{M}),$ which is 
\textrm{BMO}$_c^q(\Bbb{R},\mathcal{M})\ (q$ being the conjugate index of $%
p), $ the latter is the $L^q$-space analogue of \textrm{BMO} space already
considered in Chapters 1 and 2. These $\mathrm{BMO}_c^q(\Bbb{R},\mathcal{M})$
spaces not only are used to describe the dual of $\mathcal{H}_c^p(\Bbb{R},%
\mathcal{M})$ but also play an important role for all results in the sequel.
In particular, we will use it to prove the map $\Psi $ introduced in Chapter
3 extends to a bounded map from $L^p(L^\infty (\Bbb{R})\otimes \mathcal{M}%
,L_c^2(\widetilde{\Gamma })) $ to $\mathcal{H}_c^p(\Bbb{R},\mathcal{M})$ for
all $1<p<\infty .$ Consequently, $\mathcal{H}_c^p(\Bbb{R},\mathcal{M})$ can
be considered as a complemented subspace of $L^p(L^\infty (\Bbb{R})\otimes 
\mathcal{M},L_c^2(\widetilde{\Gamma })).$ For the most part, our results in
Chapter 4 are extension to the function space setting of results proved for
non-commutative martingales in \cite{JX}.

\section{Operator valued $\mathrm{BMO}^q$ $(q>2)$}

We will now introduce a useful operator inequality. Let $H$ be a Hilbert
space with the inner product $\langle \cdot ,\cdot \rangle $, let $a,b\in
B(H),$ then 
\begin{equation}
|a+b|^2\leq (1+t)|a|^2+(1+\frac 1t)|b|^2,\forall t>0,t\in \Bbb{R}.
\label{pisier}
\end{equation}
In fact, by Cauchy-Schwartz inequality, we have, for every $h\in H,$%
\begin{eqnarray*}
\langle |a+b|^2h,h\rangle &=&\langle (a+b)h,(a+b)h\rangle \\
&\leq &\langle ah,ah\rangle +\langle bh,bh\rangle +2\langle ah,ah\rangle
^{\frac 12}\langle bh,bh\rangle ^{\frac 12} \\
&\leq &(1+t)\langle |a|^2h,h\rangle +(1+\frac 1t)\langle |b|^2h,h\rangle
;\quad \forall t>0,t\in \Bbb{R}.
\end{eqnarray*}
\medskip

$\,$Let ${\varphi }\in L^q(\mathcal{M},L_c^2(\Bbb{R},\frac{dt}{1+t^2})).$
For $h\in \Bbb{R}^{+}\Bbb{\times R}^{+},$ denote $I_{h,t}=(t-h_1,t+h_2].$
Let 
\[
{\varphi }_h^{\#}(t)=\frac 1{h_1+h_2}\int_{I_{h,t}}|{\varphi }(x)-{\varphi }%
_{I_{h,t}}|^2dx 
\]
$\,$Set, for $2<q\leq \infty ,$ 
\[
\left\| {\varphi }\right\| _{\mathrm{BMO}_c^q}=\left\| \sup_{h\in \Bbb{R}^{+}%
\Bbb{\times R}^{+}}|{\varphi }_h^{\#}|\right\| _{L^{\frac q2}(L^\infty (\Bbb{%
R})\otimes \mathcal{M})}^{\frac 12} 
\]
and 
\[
\left\| {\varphi }\right\| _{\mathrm{BMO}_r^q}=\left\| {\varphi }%
^{*}\right\| _{\mathrm{BMO}_c^q}. 
\]
$\,$It is easy to check by (\ref{pisier}) that $\left\| {\cdot }\right\| _{%
\mathrm{BMO}_r^q}$ and $\left\| {\cdot }\right\| _{\mathrm{BMO}_c^q}$ are
norms. Let$\,\mathrm{BMO}_c^q(\Bbb{R},\mathcal{M})$ (resp. $\mathrm{BMO}_r^q(%
\Bbb{R},\mathcal{M}))$ be the space of all ${\varphi }\in L^q(\mathcal{M}%
,L_c^2(\Bbb{R},\frac{dt}{1+t^2}))$ $($resp. $L^q(\mathcal{M},L_r^2(\Bbb{R},%
\frac{dt}{1+t^2})))$ such that $\left\| {\varphi }\right\| _{\mathrm{BMO}%
_c^q}<\infty ($resp. $\left\| {\varphi }\right\| _{\mathrm{BMO}_r^q}<\infty
).$ $\mathrm{BMO}_{cr}^q(\Bbb{R},\mathcal{M})$ is defined as the
intersection of these two spaces 
\[
\mathrm{BMO}_{cr}^q(\Bbb{R},\mathcal{M})=\mathrm{BMO}_c^q(\Bbb{R},\mathcal{M}%
)\,\cap \mathrm{BMO}_r^q(\Bbb{R},\mathcal{M}) 
\]
equipped with the norm 
\[
\left\| {\varphi }\right\| _{\mathrm{BMO}_{cr}^q}=\max \{\left\| {\varphi }%
\right\| _{\mathrm{BMO}_c^q},\left\| {\varphi }\right\| _{\mathrm{BMO}%
_r^q}\}. 
\]
If $q=\infty ,$ all these spaces coincide with those introduced in Chapter
2. And if $\mathcal{M}=\Bbb{C},$ all these spaces coincide with the
classical $\mathrm{BMO}^q$. As in the case of $\mathrm{BMO}(\Bbb{R},\mathcal{%
M})$, we regard $\mathrm{BMO}_c^q(\Bbb{R},\mathcal{M})$ $($resp. $\mathrm{BMO%
}_r^q(\Bbb{R},\mathcal{M}),$ $\mathrm{BMO}_r^q(\Bbb{R},\mathcal{M}))$ as
normed spaces modulo constants. The following is the analogue for $\mathrm{%
BMO}_c^q(\Bbb{R},\mathcal{M})$ of Proposition 1.3. Recall that $%
I_t^n=(t-2^{n-1},t+2^{n-1}]$ for $t\in \Bbb{R}$ and $n\in \Bbb{Z}.$ Note
that we have trivially 
\begin{equation}
\left\| \frac 1{2^k}\int_{I_t^k}|{\varphi }(s)-{\varphi }_{_{I_t^k}}|^2ds%
\right\| _{L^{\frac q2}(L^\infty (\Bbb{R})\otimes \mathcal{M})}^{\frac
12}\leq \left\| {\varphi }\right\| _{\mathrm{BMO}_c^q}  \label{tri}
\end{equation}

\begin{proposition}
Let $2<q\leq \infty .$ Let ${\varphi }\in \mathrm{BMO}_c^q(\Bbb{R},\mathcal{M%
}).\,$Then 
\[
\left\| {\varphi }\right\| _{L^q(\mathcal{M},L_c^2(\Bbb{R},\frac{dt}{1+t^2}%
))}\leq c\left( \left\| {\varphi }\right\| _{\mathrm{BMO}_c^q}+\left\| {%
\varphi }_{I_0^1}\right\| _{L^q(\mathcal{M})}\right) .
\]
Moreover, $\mathrm{BMO}_c^q(\Bbb{R},\mathcal{M}),\mathrm{BMO}_r^q(\Bbb{R},%
\mathcal{M}),\mathrm{BMO}_{cr}^q(\Bbb{R},\mathcal{M})$ are Banach spaces.
\end{proposition}

\noindent\textbf{Proof.} The proof is similar to that of Proposition 1.3. By
(\ref{fu}) we have 
\begin{eqnarray}
|{\varphi }_{I_t^n}-{\varphi }_{I_0^1}|^2 &\leq &n(\sum_{k=3}^n|{\varphi }%
_{I_t^k}-{\varphi }_{I_t^{k-1}}|^2+|{\varphi }_{I_t^2}-{\varphi }_{I_0^1}|^2)
\nonumber \\
&\leq &n(\sum_{k=3}^n\frac 1{2^{k-1}}\int_{I_t^{k-1}}|{\varphi }(s)-{\varphi 
}_{_{I_t^k}}|^2ds+\frac 12\int_{I_0^1}|{\varphi }(s)-{\varphi }%
_{_{I_t^2}}|^2ds)  \nonumber \\
&\leq &n(\sum_{k=3}^n\frac 2{2^k}\int_{I_t^k}|{\varphi }(s)-{\varphi }%
_{_{I_t^k}}|^2ds+\frac 24\int_{I_t^2}|{\varphi }(s)-{\varphi }%
_{_{I_t^2}}|^2ds)  \nonumber \\
&=&2n\sum_{k=2}^n\frac 1{2^k}\int_{I_t^k}|{\varphi }(s)-{\varphi }%
_{_{I_t^k}}|^2ds,\quad \forall n>1,t\in [-1,1].  \label{tfai}
\end{eqnarray}
Thus by (\ref{tri}) 
\begin{equation}
\left\| |{\varphi }_{I_t^n}-{\varphi }_{I_0^1}|^2\right\| _{L^{\frac
q2}(L^\infty (\Bbb{R})\otimes \mathcal{M})}\leq 2n^2\left\| {\varphi }%
\right\| _{\mathrm{BMO}_c^q}^2,\quad \forall n>1,t\in [-1,1].  \label{faiq}
\end{equation}
To control ${\varphi }$'s $L^q(\mathcal{M},L_c^2(\Bbb{R},\frac{dt}{1+t^2}))$
norm by its $\mathrm{BMO}_c^q$ norm, we write 
\begin{eqnarray}
&&\left\| {\varphi }\right\| _{L^q(\mathcal{M},L_c^2(\Bbb{R},\frac{dt}{1+t^2}%
))}^2  \nonumber \\
&=&\left\| \int_{\Bbb{R}}\frac{|{\varphi }(s)|^2}{1+s^2}ds\right\|
_{L^{\frac q2}(\mathcal{M})}  \nonumber \\
&=&\left\| \chi _{[-\frac 12,\frac 12]}(t)\int_{\Bbb{R}}\frac{|{\varphi }%
(s)|^2}{1+s^2}ds\right\| _{L^{\frac q2}(L^\infty (\Bbb{R})\otimes \mathcal{M}%
)}  \nonumber \\
&\leq &\left\| \chi _{[-\frac 12,\frac 12]}(t)(\sum_{n=0}^\infty
\int_{I_t^{n+1}/I_t^n}\frac{|{\varphi }(s)|^2}{1+s^2}ds+\int_{I_0^1}\frac{|{%
\varphi }(s)|^2}{1+s^2}ds)\right\| _{L^{\frac q2}(L^\infty (\Bbb{R})\otimes 
\mathcal{M})}  \nonumber \\
&\leq &c(\left\| \chi _{[-\frac 12,\frac 12]}(t)(\sum_{n=2}^\infty
\int_{I_t^n}\frac{|{\varphi }(s)|^2}{2^{2n}}ds+\int_{I_0^1}|{\varphi }%
(s)|^2ds)\right\| _{L^{\frac{^q}2}(L^\infty (\Bbb{R})\otimes \mathcal{M})} 
\nonumber
\end{eqnarray}
hence by (\ref{faiq}) 
\begin{eqnarray}
\left\| {\varphi }\right\| _{L^q(\mathcal{M},L_c^2(\Bbb{R},\frac{dt}{1+t^2}%
))}^2 &\leq &c(\left\| \sum_{n=2}^\infty \chi _{[-\frac 12,\frac
12]}(t)\int_{I_t^n}\frac{|{\varphi }(s)-{\varphi }_{_{I_t^n}}|^2}{2^{2n}}%
ds\right\| _{L^{\frac{^q}2}(L^\infty (\Bbb{R})\otimes \mathcal{M})} 
\nonumber \\
&&+\left\| \sum_{n=1}^\infty \frac{|{\varphi }_{I_0^1}|^2}{2^n}\right\|
_{L^{\frac q2}(\mathcal{M})}+\sum_{n=1}^\infty \frac{n^2\left\| {\varphi }%
\right\| _{\mathrm{BMO}_c^q}^2}{2^n}  \nonumber \\
&\leq &c\sum_{n=1}^\infty \frac{(n^2+1)\left\| {\varphi }\right\| _{\mathrm{%
BMO}_c^q}^2}{2^n}+c\left\| {\varphi }_{I_0^1}\right\| _{L^q(\mathcal{M})}^2 
\nonumber \\
&<&\infty .
\end{eqnarray}
Thus $\mathrm{BMO}_c^q(\Bbb{R},\mathcal{M})$ is a Banach space. Passing to
adjoints we get that $\mathrm{BMO}_r^q(\Bbb{R},\mathcal{M})$ is a Banach
spaces and then so is $\mathrm{BMO}_{cr}^q(\Bbb{R},\mathcal{M}).\qed $

Put 
\[
{\lambda} _{\varphi} ^{n,\#}(t)=\frac 1{2^n}\mathop{\displaystyle\int
\!\!\!\int }_{T(I_t^n)}|\nabla {\varphi} |^2ydxdy. 
\]

\begin{lemma}
Let\textbf{\ }${\varphi }\in \mathrm{BMO}_c^q(\Bbb{R},\mathcal{M})$ $%
(2<q<\infty ).$ Then$\,\exists c>0$ such that 
\[
\left\| \sup_{n\in \Bbb{Z}}|{\lambda }_\varphi ^{n,\#}|\right\| _{L^{\frac
q2}(L^\infty (\Bbb{R})\otimes \mathcal{M})}\leq c\left\| {\varphi }\right\|
_{\mathrm{BMO}_c^q}^2.
\]
\end{lemma}

\noindent\textbf{Proof.} The proof is similar to that of Lemma 1.4 but more
complicated. For any $n\in \Bbb{Z},t\in \Bbb{R},$ write ${\varphi }=$ ${%
\varphi }_1^{n,t}+$ ${\varphi }_2^{n,t}+{\varphi }_3^{n,t},$ where ${\varphi 
}_1^{n,t}=({\varphi }-{\varphi }_{I_t^{n+1}})\chi _{I_t^{n+1}},$ ${\varphi }%
_2^{n,t}=({\varphi }-{\varphi }_{I_t^{n+1}})\chi _{_{(I_t^{n+1})^c}},$ and ${%
\varphi }_3^{n,t}={\varphi }_{I_t^{n+1}}.$ Set 
\[
{\lambda }_i^{n,\#}(t)=\frac 1{2^n}\mathop{\displaystyle\int \!\!\!\int }%
_{T(I_t^n)}|\nabla {\varphi }_i^{n,t}|^2ydxdy,\quad i=1,2. 
\]
Thus 
\begin{eqnarray*}
&&\left\| \sup_{n\in \Bbb{Z}}|{\lambda }_\varphi ^{n,\#}|\right\| _{L^{\frac
q2}(L^\infty (\Bbb{R})\otimes \mathcal{M})}\\
&\leq & 2\left\| \sup_{n\in \Bbb{Z}}|%
{\lambda }_1^{n,\#}|\right\| _{L^{\frac q2}(L^\infty (\Bbb{R})\otimes 
\mathcal{M})}+2\left\| \sup_{n\in \Bbb{Z}}|{\lambda }_2^{n,\#}|\right\|
_{L^{\frac q2}(L^\infty (\Bbb{R})\otimes \mathcal{M})}.\, 
\end{eqnarray*}
We treat ${\lambda }_1^{n,\#}$ first. Arguing as earlier for (\ref{gree}),
by the Green theorem we have 
\[
\frac 1{2^n}\mathop{\displaystyle\int\!\!\!\int}_{T(I_t^n)}|\nabla {\varphi }%
_1^{n,t}|^2ydxdy\leq \frac 1{2^n}\int_{-\infty }^{+\infty }|{\varphi }%
_1^{n,t}|^2ds. 
\]
Therefore, 
\begin{eqnarray}
&&\left\| \sup_{n\in \Bbb{Z}}|\frac 1{2^n}\mathop{\displaystyle\int\!\!\!%
\int}_{T(I_t^n)}|\nabla {\varphi }_1^{n,t}|^2ydxdy|\right\| _{L^{\frac
q2}(L^\infty (\Bbb{R})\otimes \mathcal{M})}  \nonumber \\
&\leq &\left\| \sup_{n\in \Bbb{Z}}|\frac 1{2^n}\int_{-\infty }^{+\infty }|{%
\varphi }_1^{n,t}|^2ds|\right\| _{L^{\frac q2}(L^\infty (\Bbb{R})\otimes 
\mathcal{M})}  \nonumber \\
&=&\left\| \sup_{n\in \Bbb{Z}}|\frac 1{2^n}\int_{I_t^{n+1}}|{\varphi }-{%
\varphi }_{I_t^{n+1}}|^2ds|\right\| _{L^{\frac q2}(L^\infty (\Bbb{R})\otimes 
\mathcal{M})}  \nonumber \\
&\leq &2\left\| {\varphi }\right\| _{\mathrm{BMO}_c^q}^2  \label{lbd1}
\end{eqnarray}
To deal with ${\lambda }_2^{n,\#},$ we note that 
\[
|\nabla P_y(x-s)|^2\leq \frac 1{4(x-s)^4}\leq \frac c{2^{4(n+k)}},\quad
\forall s\in I_t^{n+k+1}/I_t^{n+k},\quad (x,y)\in T(I_t^n). 
\]
Let $A_k=I_t^{n+k+1}/I_t^{n+k}.$ Then by (\ref{fg}), (\ref{fai}) and (\ref
{tfai}) 
\begin{eqnarray*}
&&\frac 1{2^n}\mathop{\displaystyle\int\!\!\!\int}_{T(I_t^n)}|\nabla {%
\varphi }_2^{n,t}|^2ydxdy \\
&=&\frac 1{2^n}\mathop{\displaystyle\int\!\!\!\int}_{T(I_t^n)}|\nabla
\int_{-\infty }^{+\infty }P_y(x-s){\varphi }_2^{n,t}(s)ds|^2ydxdy \\
&\leq &\frac 1{2^n}\mathop{\displaystyle\int\!\!\!\int}_{T(I_t^n)}\left(
\sum_{k=1}^\infty \int\limits_{A_k}|\nabla
P_y(x-s)|^22^{2k}ds\sum_{k=1}^\infty \int\limits_{A_k}\frac 1{2^{2k}}|{%
\varphi }_2^{n,t}(s)|^2dsy\right) dxdy \\
&\leq &\frac c{2^n}\mathop{\displaystyle\int\!\!\!\int}_{T(I_t^n)}\frac
1{2^{3n}}\sum_{k=1}^\infty \int\limits_{A_k}\frac 1{2^{2k}}|{\varphi }-{%
\varphi }_{I_t^{n+1}}|^2dsydxdy \\
&\leq &\frac c{2^n}\sum_{k=1}^\infty \int\limits_{A_k}\frac 2{2^{2k}}(|{%
\varphi }-{\varphi }_{I_t^{n+k+1}}|^2+|{\varphi }_{I_t^{n+k+1}}-{\varphi }%
_{I_t^{n+1}}|^2)ds \\
&\leq &c\sum_{k=1}^\infty \frac 1{2^{2k+n}}\int\limits_{A_k}|{\varphi }-{%
\varphi }_{I_t^{n+k+1}}|^2ds+\sum_{k=1}^\infty \frac c{2^k}\sum_{i=1}^k\frac{%
2k}{2^{n+i}}\int_{I_t^{n+i}}|{\varphi }(u)-{\varphi }_{I_t^{n+i}}|^2du \\
&\leq &cX_n+cY_n
\end{eqnarray*}
where 
\begin{eqnarray*}
X_n &=&\sum_{k=1}^\infty \frac 1{2^{2k+n}}\int\limits_{A_k}|{\varphi }-{%
\varphi }_{I_t^{n+k+1}}|^2ds, \\
Y_n &=&\sum_{k=1}^\infty \frac k{2^k}\sum_{i=1}^k\frac
1{2^{n+i}}\int_{I_t^{n+i}}|{\varphi }(s)-{\varphi }_{I_t^{n+i}}|^2ds.
\end{eqnarray*}
$X_n,Y_n$ are estimated as follows. For $X_n$ we have 
\begin{eqnarray*}
&&\left\| \sup_{n\in \Bbb{Z}}|X_n|\right\| _{L^{\frac q2}(L^\infty (\Bbb{R}%
)\otimes \mathcal{M})} \\
&=&\left\| \sup_{n\in \Bbb{Z}}|\sum_{k=1}^\infty \frac 1{2^k}\frac
1{2^{n+k}}\int\limits_{A_k}|{\varphi }-{\varphi }_{I_t^{n+k+1}}|^2ds|\right%
\| _{L^{\frac q2}(L^\infty (\Bbb{R})\otimes \mathcal{M})} \\
&\leq &\sum_{k=1}^\infty \frac 1{2^k}\left\| \sup_{n\in \Bbb{Z}}|\frac
1{2^{n+k}}\int\limits_{I_t^{n+k+1}}|{\varphi }-{\varphi }%
_{I_t^{n+k+1}}|^2ds|\right\| _{L^{\frac q2}(L^\infty (\Bbb{R})\otimes 
\mathcal{M})} \\
&\leq &2\left\| {\varphi }\right\| _{\mathrm{BMO}_c^q}^2.
\end{eqnarray*}
On the other hand, 
\begin{eqnarray*}
&&\left\| \sup_{n\in \Bbb{Z}}|Y_n|\right\| _{L^{\frac q2}(L^\infty (\Bbb{R}%
)\otimes \mathcal{M})} \\
&\leq &\sum_{k=1}^\infty \frac k{2^k}\sum_{i=1}^k\left\| \sup_{n\in \Bbb{Z}%
}|\frac 1{2^{n+i}}\int_{I_t^{n+i}}|{\varphi }(s)-{\varphi }%
_{I_t^{n+i}}|^2ds|\right\| _{L^{\frac q2}(L^\infty (\Bbb{R})\otimes \mathcal{%
M})} \\
&\leq &\sum_{k=1}^\infty \frac{k^2}{2^k}\left\| {\varphi }\right\| _{\mathrm{%
BMO}_c^q}^2 \\
&=&6\left\| {\varphi }\right\| _{\mathrm{BMO}_c^q}^2.
\end{eqnarray*}
Combining the preceding inequalities we get 
\[
\left\| \sup_{n\in \Bbb{Z}}|{\lambda }_{{\varphi }_2}^{n,\#}|\right\|
_{L^{\frac q2}(L^\infty (\Bbb{R})\otimes \mathcal{M})}\leq c\left\| {\varphi 
}\right\| _{\mathrm{BMO}_c^q}^2, 
\]
which, together with (\ref{lbd1}), yields 
\[
\left\| \sup_{n\in \Bbb{Z}}|{\lambda }_\varphi ^{n,\#}|\right\| _{L^{\frac
q2}(L^\infty (\Bbb{R})\otimes \mathcal{M})}\leq c\left\| {\varphi }\right\|
_{\mathrm{BMO}_c^q}^2.\qed
\]
Set 
\[
{\varphi }_n^{\#}(t)=\frac 1{2^n}\int_{I_t^n}|{\varphi }(x)-{\varphi }%
_{I_t^n}|^2dx 
\]
$\,$Notice that for every $h\in \Bbb{R}^{+}\times \Bbb{R}^{+}$ there exists $%
n\in \Bbb{Z}$ such that $(t-h_1,t+h_2)\in I_t^n$ for every $t\in \Bbb{R}$
and $2^n\leq 4(h_1+h_2),$ we have 
\begin{equation}
\frac 14\left\| {\varphi }\right\| _{\mathrm{BMO}_c^q}\leq \left\| \sup_n{%
\varphi }_n^{\#}\right\| _{L^{\frac q2}(L^{^\infty }(\Bbb{R})\otimes 
\mathcal{M})}^{\frac 12}\leq \left\| {\varphi }\right\| _{\mathrm{BMO}_c^q}.
\label{n}
\end{equation}

\begin{lemma}
The operator $\Psi $ defined in Chapter 2 extends to a bounded map from $%
L^q(L^\infty (\Bbb{R)}\otimes \mathcal{M},L_c^2(\widetilde{\Gamma }))\
(2<q<\infty )$ into $\mathrm{BMO}_c^q(\Bbb{R},\mathcal{M})$ and there exists 
$c_q>0$ such that 
\begin{equation}
\left\| \Psi (h)\right\| _{\mathrm{BMO}_c^q}\leq c_q\left\| h\right\|
_{L^q(L^\infty (\Bbb{R)}\otimes \mathcal{M},L_c^2)}.  \label{psibmoq}
\end{equation}
\end{lemma}

\noindent\textbf{Proof.} The pattern of this proof is similar to that of
Lemma 2.2. One new thing we need is the non-commutative Hardy-Littlewood
maximal inequality proved in the previous chapter.

Let $\mathcal{S}$ be the family of functions introduced in the proof of
Lemma 2.2. Since $\mathcal{S}$ is dense in $L^q(L^\infty (\Bbb{R)}\otimes 
\mathcal{M},L_c^2(\widetilde{\Gamma }))$, we need only to prove (\ref
{psibmoq}) for all $h\in \mathcal{S}$. Fix $h\in \mathcal{S}$ and set ${%
\varphi }=\Psi (h).$ Then ${\varphi }\in L^q(\mathcal{M},L_c^2(\Bbb{R},\frac{%
ds}{1+s^2}))$. Let $u\in \Bbb{R}$ and $n\in \Bbb{Z}.$ Set 
\begin{eqnarray*}
h_1^u(x,y,t) &=&h(x,y,t)\chi _{_{I_u^{n+1}}}(t), \\
h_2^u(x,y,t) &=&h(x,y,t)\chi _{_{(I_u^{n+1})^c}}(t)
\end{eqnarray*}
and 
\[
B_{I_u^n}=\displaystyle\int_{-\infty }^{+\infty }\mathop{\displaystyle}%
\mathop{\displaystyle\int\!\!\!\int}_\Gamma Q_{I_u^n}h_2^udydxdt, 
\]
where 
\[
Q_{I_u^n}(x,y,t)=\frac 1{2^n}\int_{I_u^n}Q_y(x+t-s)ds 
\]
(recall that $Q_y(x)$ is defined by (\ref{qy}) as the gradient of the
Poisson kernel). Then 
\begin{eqnarray*}
{\varphi }_n^{\#}(u) &\leq &\frac 4{2^n}\int_{I_u^n}|{\varphi }%
(s)-B_{I_t^n}|^2ds \\
&\leq &\frac 8{2^n}\int_{I_u^n}|\int_{(I_u^{n+1})^c}\mathop{\displaystyle%
\int\!\!\!\int}_\Gamma (Q_y(x+t-s)-Q_{I_u^n})hdxdydt|^2ds \\
&&+\frac 8{2^n}\int_{I_u^n}|\int_{-\infty }^{+\infty }\mathop{\displaystyle%
\int\!\!\!\int}_\Gamma Q_y(x+t-s)h_1^udxdydt|^2ds \\
&=&8A_n+\frac 8{2^n}\int_{I_u^n}|\int_{I_u^{n+1}}\mathop{\displaystyle\int\!%
\!\!\int}_\Gamma Q_y(x+t-s)hdxdydt|^2ds
\end{eqnarray*}
Recall that, as noted earlier in (\ref{qt}), 
\[
\mathop{\displaystyle\int\!\!\!\int}_\Gamma |Q_y(x+t-s)-Q_{I_u^n}|^2dxdy\leq
c2^{2n}(t-u)^{-4} 
\]
for $t\in (I_u^{n+1})^c$ and $s\in I_u^n.$ By (\ref{fg}), we have 
\begin{eqnarray*}
A_n &=&\frac 1{2^n}\int_{I_u^n}|\int_{(I_u^{n+1})^c}\mathop{\displaystyle%
\int\!\!\!\int}_\Gamma (Q_y(x+t-s)-Q_{I_u^n})hdxdydt|^2ds \\
&\leq &\int_{(I_u^{n+1})^c}c2^{2n}(t-u)^{-2}dt\int_{(I_u^{n+1})^c}(t-u)^{-2}%
\mathop{\displaystyle\int\!\!\!\int}_\Gamma |h|^2dxdydt \\
&=&c2^n\int_{(I_u^{n+1})^c}(t-u)^{-2}\mathop{\displaystyle\int\!\!\!\int}%
_\Gamma |h|^2dxdydt
\end{eqnarray*}
Then, for any positive $(a_n)_{n\in \Bbb{Z}}$ such that $\left\| \sum_{k\in 
\Bbb{Z}}a_n\right\| _{L^{(\frac q2)^{\prime }}(L^\infty (\Bbb{R})\otimes 
\mathcal{M)}}\leq 1,$%
\begin{eqnarray*}
&&\tau \sum_{n\in \Bbb{Z}}\int_{-\infty }^{+\infty }{\varphi }%
_n^{\#}(u)a_n(u)du \\
&\leq &\sum_{n\in \Bbb{Z}}\tau \int_{-\infty }^{+\infty
}c2^n\int_{(I_u^{n+1})^c}(t-u)^{-2}\mathop{\displaystyle\int\!\!\!\int}%
_\Gamma |h|^2dxdydta_n(u)du \\
&&+\sum_{n\in \Bbb{Z}}\tau \int_{-\infty }^{+\infty }\frac
8{2^n}\int_{I_u^n}|\int_{I_u^{n+1}}\mathop{\displaystyle\int\!\!\!\int}%
_\Gamma Q_y(x+t-s)hdxdydt|^2dsa_n(u)du \\
&=&A+B
\end{eqnarray*}
By the non-commutative H\"{o}lder inequality, 
\begin{eqnarray*}
A &=&\sum_{n\in \Bbb{Z}}\tau \int_{-\infty }^{+\infty
}c2^n\int_{(I_t^{n+1})^c}(t-u)^{-2}a_n(u)du\mathop{\displaystyle\int\!\!\!%
\int}_\Gamma |h|^2dxdydt \\
&\leq &\left\| \mathop{\displaystyle\int\!\!\!\int}_\Gamma |h|^2dxdy\right\|
_{L^{\frac q2}(L^\infty (\Bbb{R})\otimes \mathcal{M)}}\left\| \sum_{n\in 
\Bbb{Z}}c2^n\int_{(I_t^n)^c}(t-u)^{-2}a_n(u)du\right\| _{L^{(\frac
q2)^{\prime }}(L^\infty (\Bbb{R})\otimes \mathcal{M)}} \\
&\leq &\left\| h\right\| _{L^q(L^\infty (\Bbb{R)}\otimes \mathcal{M},L_c^2(%
\widetilde{\Gamma }))}^2\left\| \sum_{n\in \Bbb{Z}}\sum_{k=n}^{+\infty
}2^n\int_{I_t^{k+1}}\frac 1{2^{2k}}a_n(u)du\right\| _{L^{(\frac q2)^{\prime
}}(L^\infty (\Bbb{R})\otimes \mathcal{M)}}.
\end{eqnarray*}
Let us estimate the second factor in the last term. By (\ref{stein}), 
\begin{eqnarray*}
&&\left\| \sum_{n\in \Bbb{Z}}\sum_{k=n+1}^{+\infty }2^n\int_{I_t^{k+1}}\frac
1{2^{2k}}a_n(u)du\right\| _{L^{(\frac q2)^{\prime }}(L^\infty (\Bbb{R}%
)\otimes \mathcal{M)}} \\
&=&\left\| \sum_{k\in \Bbb{Z}}\frac 1{2^k}\int_{I_t^{k+1}}\sum_{n=-\infty
}^{k-1}\frac{2^n}{2^k}a_n(u)du\right\| _{L^{(\frac q2)^{\prime }}(L^\infty (%
\Bbb{R})\otimes \mathcal{M)}} \\
&\leq &c_q\left\| \sum_{k\in \Bbb{Z}}\sum_{n=-\infty }^{k-1}\frac{2^n}{2^k}%
a_n\right\| _{L^{(\frac q2)^{\prime }}(L^\infty (\Bbb{R})\otimes \mathcal{M)}%
} \\
&\leq &c_q\left\| \sum_{n\in \Bbb{Z}}a_n\right\| _{L^{(\frac q2)^{\prime
}}(L^\infty (\Bbb{R})\otimes \mathcal{M)}}\leq c_q.
\end{eqnarray*}
Thus 
\[
A\leq c_q\left\| h\right\| _{L^q(L^\infty (\Bbb{R)}\otimes \mathcal{M}%
,L_c^2)}^2. 
\]
For the term $B,$ by (\ref{stein}), (\ref{h2l2}) and Cauchy-Schwartz
inequality, 
\begin{eqnarray*}
B &\leq &\sum_{n\in \Bbb{Z}}\tau \int_{\Bbb{R}}\frac 8{2^n}\int_{\Bbb{R}%
}|\int_{I_u^{n+1}}\mathop{\displaystyle\int\!\!\!\int}_\Gamma
Q_y(x+t-s)hdxdydt|^2dsa_n(u)du \\
&=&\sum_{n\in \Bbb{Z}}\int_{\Bbb{R}}\frac 8{2^n}\sup_{\|f\|_{L^2(L^\infty (
\Bbb{R})\otimes \mathcal{M})}=1}(\tau \!\int_{\Bbb{R}}\!\!\int_{I_u^{n+1}}\!
\!\mathop{\displaystyle\int\!\!\!\!\int}_\Gamma Q_y(x+t-s)ha_n^{\frac
12}(u)dxdydtf(s)ds)^2du \\
&=&\sum_{n\in \Bbb{Z}}\int_{\Bbb{R}}\frac 8{2^n}\sup_{||f||_{L^2(L^\infty (%
\Bbb{R)}\otimes \mathcal{M})}=1}(\tau \int_{I_u^{n+1}}\mathop{\displaystyle%
\int\!\!\!\int}_\Gamma ha_n^{\frac 12}(u)\nabla f(t+x,y)dxdydt)^2du \\
&\leq &\sum_{n\in \Bbb{Z}}\int_{\Bbb{R}}\frac 8{2^n}\tau \int_{I_u^{n+1}}%
\mathop{\displaystyle\int\!\!\!\int}_\Gamma |h|^2a_n(u)dxdydtdu \\
&=&\sum_{n\in \Bbb{Z}}\tau \int_{\Bbb{R}}\mathop{\displaystyle\int\!\!\!\int}%
_\Gamma |h|^2dxdy\frac 8{2^n}\int_{I_t^{n+1}}a_n(u)dudt \\
&\leq &||\mathop{\displaystyle\int\!\!\!\int}_\Gamma |h|^2dxdy||_{L^{\frac
q2}(L^\infty (\Bbb{R})\otimes \mathcal{M)}}\left\| \sum_{n\in \Bbb{Z}}\frac{%
16}{2^n}\int_{I_t^n}a_n(u)du\right\| _{L^{(\frac q2)^{\prime }}(L^\infty (%
\Bbb{R})\otimes \mathcal{M)}} \\
&\leq &c_q\left\| h\right\| _{L^q(L^\infty (\Bbb{R)}\otimes \mathcal{M}%
,L_c^2)}^2.
\end{eqnarray*}
Thus 
\[
\left\| \sup_n|{\varphi }_n^{\#}|\right\| _{L^{\frac q2}(L^\infty (\Bbb{R}%
)\otimes \mathcal{M})}\leq c_q\left\| h\right\| _{L^q(L^\infty (\Bbb{R)}%
\otimes \mathcal{M},L_c^2)}^2 
\]
and then 
\[
||\Psi (h)||_{\mathrm{BMO}_c^q}\leq c_q\left\| h\right\| _{L^q(L^\infty (%
\Bbb{R)}\otimes \mathcal{M},L_c^2)}.\qed
\]
\medskip

\medskip \noindent
\textbf{Remark.} \medskip It seems difficult to define non-commutative $%
\mathrm{BMO}^q$ for $q<2.$

\section{The duality theorem of $\mathcal{H}^p$ and $\mathrm{BMO}^q (1<p<2)$}

Denote by $\mathcal{H}_{c0}^p(\Bbb{R},\mathcal{M})$ (resp. $\mathcal{H}%
_{r0}^p(\Bbb{R},\mathcal{M})$) the functions $f$ in $\mathcal{H}_c^p(\Bbb{R},%
\mathcal{M})$ (resp. $\mathcal{H}_r^p(\Bbb{R},\mathcal{M})$) such that $f\in
L^p(\mathcal{M},L_c^2(\Bbb{R},(1+t^2)dt))$ (resp. $L^p(\mathcal{M},$ $L_r^2(%
\Bbb{R},(1+t^2)dt)$) and $\int fdt=0.$ Set 
\[
\mathcal{H}_{cr0}^p(\Bbb{R},\mathcal{M})= \mathcal{H}_{c0}^p(\Bbb{R},%
\mathcal{M}) +\mathcal{H}_{r0}^p(\Bbb{R},\mathcal{M}). 
\]
It is easy to see that $\mathcal{H}_{c0}^p(\Bbb{R},\mathcal{M})$ (resp. $%
\mathcal{H}_{r0}^p(\Bbb{R},\mathcal{M}),$ $\mathcal{H}_{cr0}^p(\Bbb{R},%
\mathcal{M})$) is a dense subspace of $\mathcal{H}_c^p(\Bbb{R},\mathcal{M})$
(resp. $\mathcal{H}_r^p(\Bbb{R}$, $\mathcal{H}_{cr0}^p(\Bbb{R},\mathcal{M})$%
). By Propositions 1.1 and 4.1, $\int_{-\infty }^{+\infty }{\varphi} ^{*}fdt$
exists as an element in $L^1(\mathcal{M})$ for any ${\varphi}\in \mathrm{BMO}%
_c^q(\Bbb{R},\mathcal{M})$ and $f\in \mathcal{H}_{c0}^p(\Bbb{R},\mathcal{M})$
.

\begin{theorem}
Let $1<p<2,$ $q=\frac p{p-1}.$ Then

(a) $(\mathcal{H}_c^p(\Bbb{R},\mathcal{M}))^{*}=\mathrm{BMO}_c^q(\Bbb{R},%
\mathcal{M})$ with equivalent norms. More precisely, every ${\varphi }\in 
\mathrm{BMO}_c^q(\mathcal{M})\,$defines a continuous linear functional on $%
\mathcal{H}_c^p(\Bbb{R},\mathcal{M})\,$ by 
\begin{equation}
l_\varphi (f)=\tau \int_{-\infty }^{+\infty }{\varphi }^{*}fdt;\qquad
\forall f\in \mathcal{H}_{c0}^p(\Bbb{R},\mathcal{M})  \label{theorem4.1}
\end{equation}
Conversely every $l\in (\mathcal{H}_c^p(\Bbb{R},\mathcal{M}))^{*}$ can be
given as above by some ${\varphi }\in \mathrm{BMO}_c^q(\Bbb{R},\mathcal{M})\,
$and there exist constants $c,c_q>0$ such that 
\[
c_q\left\| {\varphi }\right\| _{\mathrm{BMO}_c^q}\leq \left\| l_\varphi
\right\| _{(\mathcal{H}_c^p)^{*}}\leq c\left\| {\varphi }\right\| _{\mathrm{%
BMO}_c^q}
\]
Thus $(\mathcal{H}_c^p(\Bbb{R},\mathcal{M}))^{*}=\mathrm{BMO}_c^q(\Bbb{R},%
\mathcal{M})$ with equivalent norms.

(b) Similarly, $(\mathcal{H}_r^p(\Bbb{R},\mathcal{M}))^{*}=\mathrm{BMO}_r^q(%
\Bbb{R},\mathcal{M})$ with equivalent norms.

(c) $(\mathcal{H}_{cr}^p(\Bbb{R},\mathcal{M}))^{*}=\mathrm{BMO}_{cr}^q(\Bbb{R%
},\mathcal{M})$ with equivalent norms.
\end{theorem}

\noindent\textbf{Proof. }(i) Let ${\varphi }\in \mathrm{BMO}_c^q(\Bbb{R},%
\mathcal{M})$ and $f\in \mathcal{H}_{c0}^p(\Bbb{R},\mathcal{M})$. As in the
proof of Theorem 2.4, we assume ${\varphi }$ and $f$ compactly supported.
Let $G_c(f)\,$ and $\widetilde{S}_c(f)$ be as in the proof of Theorem 2.4.
Similar to what we have explained there, $G_c(f)(x,y)\,$ can be assumed to
be invertible in $\mathcal{M}\,$for every $(x,y)\in \Bbb{R}_{+}^2.$ By the
Green theorem and the Cauchy-Schwartz inequality (see the corresponding part
of the proof of Theorem 2.4 to see why the Green theorem works well), 
\begin{eqnarray*}
|l_\varphi (f)| &=&2|\tau \int_{-\infty }^{+\infty }\int_0^\infty \nabla {%
\varphi }^{*}\nabla fydydx| \\
\ &\leq &2(\tau \int_{-\infty }^{+\infty }\int_0^\infty
G_c^{p-2}(f)(x,y)|\nabla f|^2(x,y)ydydx)^{\frac 12} \\
&&\bullet (3\tau \int_{-\infty }^{+\infty }\int_0^\infty \widetilde{S}%
_c^{2-p}(f)(x,\frac y4)|\nabla {\varphi }|^2ydydx)^{\frac 12} \\
&=&2I\bullet II
\end{eqnarray*}
Noting that $G_c^{p-1}(f)(x,y)\leq G_c^{p-1}(f)(x,0)\,$, we have 
\begin{eqnarray*}
I^2 &=&\tau \int_{-\infty }^{+\infty }\int_0^\infty -G_c^{p-2}(f)(x,y)\frac{%
\partial G_c^2(f)}{\partial y}(x,y)dydx \\
&=&\tau \int_{-\infty }^{+\infty }\int_0^\infty (-G_c^{p-2}(f)(x,y)\frac{%
\partial G_c(f)}{\partial y}G_c(f)(x,y) \\
&&-G_c^{p-1}(f)\frac{\partial G_c(f)}{\partial y}(x,y))dydx \\
\ &=&2\tau \int_{-\infty }^{+\infty }\int_0^\infty -G_c^{p-1}(f)(x,y)\frac{%
\partial G_c(f)}{\partial y}dydx \\
&\leq &2\tau \int_{-\infty }^{+\infty }\int_0^\infty -G_c^{p-1}(f)(x,0)\frac{%
\partial G_c(f)}{\partial y}(x,y)dxdy \\
&\leq &2\tau \int_{-\infty }^{+\infty }G_c^p(f)(x,0)dx \\
&\leq &6\tau \int_{-\infty }^{+\infty }S_c^p(f)(x)dx \\
&=&6\left\| f\right\| _{\mathcal{H}_c^p}^p
\end{eqnarray*}
Define 
\[
\delta ^k(x)=\widetilde{S}_c^{2-p}(f)(x,2^k)-\widetilde{S}%
_c^{2-p}(f)(x,2^{k+1}),\quad \forall x\in \Bbb{R}. 
\]
Then $\delta ^k\in L^{\frac p{2-p}}(L^{^\infty }(\Bbb{R})\otimes \mathcal{M}%
) $ is positive. Note that $(\frac q2)^{\prime }=\frac p{p-2}.$ Moreover, 
\begin{eqnarray*}
\delta ^k(x) &=&\delta ^k(x^{\prime }),\forall (i-1)2^j<x,x^{\prime }\leq
i2^j \\
\sum_{k=-\infty }^\infty \delta ^k(x) &=&\widetilde{S}_c^{2-p}(f)(x,0)
\end{eqnarray*}
Arguing as earlier for Theorem 2.4, we have 
\begin{eqnarray*}
II^2 &=&3\tau \int_{-\infty }^{+\infty }\sum_{k=-\infty }^\infty \widetilde{S%
}_c^{2-p}(f)(x,2^k)\int_{2^{k+2}}^{2^{k+3}}|\nabla {\varphi }|^2ydydx \\
&=&3\tau \int_{-\infty }^{+\infty }\sum_{k=-\infty }^\infty
(\sum_{j=k}^\infty \delta ^j(x))\int_{2^{k+2}}^{2^{k+3}}|\nabla {\varphi }%
|^2ydydx \\
&=&3\tau \int_{-\infty }^{+\infty }\sum_{j=-\infty }^\infty 2^j\delta
^j(x)\frac 1{2^j}\int_0^{2^{j+3}}|\nabla {\varphi }|^2ydydx \\
&\leq &3\tau \int_{-\infty }^{+\infty }\sum_{j=-\infty }^\infty
\int_{x-2^j}^{x+2^j}\delta ^j(t)dt\frac 1{2^j}\int_0^{2^{j+3}}|\nabla {%
\varphi }|^2ydydx \\
&=&24\tau \sum_{j=-\infty }^\infty \int_{-\infty }^{+\infty }\delta
^j(t)\frac 1{2^{j+3}}\int_{t-2^j}^{t+2^j}\int_0^{2^{j+3}}|\nabla {\varphi }%
|^2ydydxdt
\end{eqnarray*}
hence by (\ref{junge}) and Lemma 4.2 
\begin{eqnarray*}
II^2 &\leq &24\left\| \sum_{j=-\infty }^\infty \delta ^j(t)\right\|
_{L^{(\frac q2)^{\prime }}}\left\| \sup_j|\frac
1{2^{j+3}}\int_{t-2^j}^{t+2^j}\int_0^{2^{j+3}}|\nabla {\varphi }%
|^2ydydx|\right\| _{L^{\frac q2}} \\
&\leq &c\left\| f\right\| _{\mathcal{H}_c^p}^{2-p}\left\| {\varphi }\right\|
_{\mathrm{BMO}_c^q}^2.
\end{eqnarray*}
Combining the preceding estimates on I and II, we get 
\[
|l_\varphi (f)|\leq c\left\| {\varphi }\right\| _{\mathrm{BMO}_c^q}\left\|
f\right\| _{\mathcal{H}_c^p}. 
\]
Therefore, $l_\varphi $ defines a continuous functional on $\mathcal{H}_c^p$
of norm smaller than $c\left\| {\varphi }\right\| _{\mathrm{BMO}_c^q}.$%
\medskip

(ii) Now suppose $l\in (\mathcal{H}_c^p)^{*}.$ Then by the Hahn-Banach
theorem $l$ extends to a continuous functional on $L^p(L^\infty (\Bbb{R)}%
\otimes \mathcal{M},L_c^2(\widetilde{\Gamma }))$ of the same norm. Thus by 
\[
(L^p(L^\infty (\Bbb{R)}\otimes \mathcal{M},L_c^2(\widetilde{\Gamma }%
)))^{*}=L^q(L^\infty (\Bbb{R)}\otimes \mathcal{M},L_c^2(\widetilde{\Gamma }%
)) 
\]
there exists $h\in L^q(L^\infty (\Bbb{R)}\otimes \mathcal{M},L_c^2(%
\widetilde{\Gamma }))$ such that 
\[
||h||_{L^q(L^\infty (\Bbb{R)}\otimes \mathcal{M},L_c^2(\widetilde{\Gamma }%
))}^2=||\mathop{\displaystyle\int\!\!\!\int}_\Gamma
h^{*}(x,y,t)h(x,y,t)dydx||_{L^{\frac q2}(L^\infty (\Bbb{R)}\otimes \mathcal{%
M)}}=||l||^2 
\]
and 
\begin{eqnarray*}
l(f) &=&\tau \int_{-\infty }^{+\infty }\mathop{\displaystyle\int\!\!\!\int}%
_\Gamma h^{*}(x,y,t)\nabla f(t+x,y)dydxdt \\
&=&\tau \int_{-\infty }^{+\infty }\Psi ^{*}(h)f(s)ds.
\end{eqnarray*}
Let 
\begin{equation}
{\varphi} =\Psi (h)  \label{faih}
\end{equation}
Then 
\[
l(f)=\tau \int_{-\infty }^{+\infty }{\varphi} ^{*}(s)f(s)ds 
\]
and by Lemma 4.3 $||{\varphi} ||_{\mathrm{BMO}_c^q}\leq c_q||l||.$ This
finishes the proof of the theorem concerning $\mathcal{H}_c^p$ and $\mathrm{%
BMO}_c^q.$ Passing to adjoints yields the part on $\mathcal{H}_r^p$ and $%
\mathrm{BMO}_r^q.\,$Finally, the duality between $\mathcal{H}_{cr}^p$ and $%
\mathrm{BMO}_{cr}^q\mathrm{\,}$is obtained by the classical fact that the
dual of a sum is the intersection of the duals.\qed

\begin{corollary}
${\varphi }\in \mathrm{BMO}_c^q(\Bbb{R},\mathcal{M})$ if and only if 
\[
\left\| \sup_{n\in \Bbb{Z}}|{\lambda }_\varphi ^{n,\#}|\right\| _{L^{\frac
q2}(L^{^\infty }(\Bbb{R})\otimes \mathcal{M})}<\infty \,
\]
and there exist $c,c_q>0$ such that 
\[
c_q\left\| {\varphi }\right\| _{\mathrm{BMO}_c^q}^2\leq \left\| \sup_{n\in 
\Bbb{Z}}|{\lambda }_\varphi ^{n,\#}|\right\| _{L^{\frac q2}(L^\infty (\Bbb{R}%
)\otimes \mathcal{M})}\leq c\left\| {\varphi }\right\| _{\mathrm{BMO}_c^q}^2.
\]
\end{corollary}

\noindent\textbf{Proof.} From the proof of Theorem 4.4, if ${\varphi} $ is
such that 
\[
\left\| \sup_n|{\lambda} _{\varphi} ^{n,\#}|\right\| _{L^{\frac q2}(L^\infty
(\Bbb{R})\otimes \mathcal{M})}<\infty , 
\]
then ${\varphi} $ defines a continuous linear functional on $\mathcal{H}%
_{c0}^p $ by $l_{\varphi} =\tau \int_{-\infty }^{+\infty }{\varphi} ^{*}fdt$
and 
\[
\left\| l_{\varphi} \right\| _{(\mathcal{H}_c^p)^{*}}\leq c\left\| \sup_n|{%
\lambda} _{\varphi} ^{n,\#}|\right\| _{L^{\frac q2}(L^\infty (\Bbb{R}%
)\otimes \mathcal{M})}^{\frac 12} 
\]
and then by Theorem 4.4 again, there exists a function ${\varphi} ^{{\prime }%
}\in $\textrm{BMO}$_c^q(\Bbb{R},\mathcal{M})$ with 
\[
\left\| {\varphi} ^{\prime }\right\| _{\mathrm{BMO}_c^q}^2\leq c_q\left\|
l_{\varphi} \right\| _{(\mathcal{H}_c^p)^{*}}^2\leq c_q\left\| \sup_n{\lambda%
} _{\varphi} ^{n,\#}\right\| _{L^{\frac q2}(L^\infty (\Bbb{R})\otimes 
\mathcal{M})} 
\]
such that 
\[
\tau \int_{-\infty }^{+\infty }{\varphi} ^{*}fdt=\tau \int_{-\infty
}^{+\infty }{\varphi} ^{\prime *}fdt. 
\]
Thus ${\varphi} \in $\textrm{BMO}$_c^q(\Bbb{R},\mathcal{M})\,$and $\left\| {%
\varphi} \right\| _{\mathrm{BMO}_c^q}^2\leq c_q\left\| \sup_n{\lambda}
_{\varphi} ^{n,\#}\right\| _{L^{\frac q2}(L^{^\infty }(\Bbb{R})\otimes 
\mathcal{M})}.$ Combining this with Lemma 4.2, we get the desired assertion.%
\qed

\medskip

Now we are in a position to show that as in the classical case, the Lusin
square function and the Littlewood-Paley $g$-function have equivalent $L^p$%
-norm in the non-commutative setting. The case $p=1$ was already obtained in
Chapter 2.

\begin{theorem}
For $f\in \mathcal{H}_c^p(\Bbb{R},\mathcal{M})($resp. $\mathcal{H}_r^p(\Bbb{R%
},\mathcal{M})),$ $1\leq p<\infty ,$ we have 
\begin{eqnarray}
c_p^{-1}\left\| G_c(f)\right\| _p\leq \left\| S_c(f)\right\| _p\leq
c_p\left\| G_c(f)\right\| _p;  \label{GSC} \\
c_p^{-1}\left\| G_r(f)\right\| _p\leq \left\| S_r(f)\right\| _p\leq
c_p\left\| G_r(f)\right\| _p.
\end{eqnarray}
\end{theorem}

\noindent\textbf{Proof.} We need only to prove the second inequality of (\ref
{GSC}). The case of $p=2$ is obvious. The case of $p=1$ is Corollary 2.7 and
the part of $1<p<2$ can be proved similarly by using the following
inequality already obtained during the proof of Theorem 4.4 
\[
|\tau \int {\varphi} ^{*}fdt|\leq c\left\| {\varphi} \right\| _{\mathrm{BMO}%
_c^q}\left\| G_c(f)\right\| _p^{\frac p2}\left\| S_c(f)\right\| _p^{1-\frac
p2}. 
\]
For $p>2,$ let $g$ be a positive element in $L^{_{^{(\frac p2)^{{\prime }%
}}}}(L^\infty (\Bbb{R})\otimes \mathcal{M})$ with $\left\| g\right\|
_{(\frac p2)^{{\prime }}}\leq 1.$ By (\ref{junge}) and (\ref{MAXC}) we have 
\begin{eqnarray*}
&&\left| \tau \int_{\Bbb{R}}\mathop{\displaystyle\int\!\!\!\int}_\Gamma
|\nabla f(x+t,y)|^2dxdyg(t)dt\right| \\
&=&\left| \tau \mathop{\displaystyle\int\!\!\!\int}_{\Bbb{R}_{+}^2}|\nabla
f(x,y)|^2y\frac 1y\int_{x-y}^{x+y}g(t)dtdxdy\right| \\
&\leq &4\left| \tau \int_{\Bbb{R}}\sum_{n=-\infty }^{+\infty
}\int_{2^{n-1}}^{2^{n}}|\nabla f(x,y)|^2ydy\frac
1{2^{n+1}}\int_{x-2^{n}}^{x+2^{n}}g(t)dtdx\right| \\
&\leq &4\left\| \int_{\Bbb{R}_{+}}\!|\nabla f(x,y)|^2ydy\right\| _{L^{\frac
p2}(L^\infty (\Bbb{R})\otimes \mathcal{M})}\left\| \sup_n|\frac
1{2^{n+1}}\int_{x-2^n}^{x+2^n}\!\!g(t)dt|\right\| _{L^{(\frac p2)^{{%
\prime }}}(L^\infty (\Bbb{R})\otimes \mathcal{M})} \\
&\leq &c_p\left\| G_c(f)\right\| _p^2
\end{eqnarray*}
Therefore, taking the supremum over all $g$ as above, we obtain 
\[
\left\| S_c(f)\right\| _p^2\leq c_p\left\| G_c(f)\right\| _p^2.\qed
\]

\section{The equivalence of $\mathcal{H}^q$ and $\mathrm{BMO}^q (q>2)$}

The following is the analogue for functions of a result for non-commutative
martingales proved in \cite{JX}.

\begin{theorem}
$\mathcal{H}_c^p(\Bbb{R},\mathcal{M})=\mathrm{BMO}_c^p(\Bbb{R},\mathcal{M})$
with equivalent norms for $2<p<\infty .$
\end{theorem}

\noindent\textbf{Proof.} Note that for every ${\varphi} \in \mathcal{H}_c^p(%
\Bbb{R},\mathcal{M})$ and every $g\in \mathcal{H}_c^{p^{\prime }}(\Bbb{R},%
\mathcal{M})$ ($p^{\prime }=\frac p{p-1}$) 
\begin{eqnarray*}
&&|\tau \int_{-\infty }^{+\infty }\mathop{\displaystyle\int\!\!\!\int}%
_\Gamma \nabla g(x+t,y)\nabla {\varphi} ^{*}(x+t,y)dxdydt| \\
&\leq &\left\| \nabla g(x+t,y)\right\| _{L^{p^{\prime }}(L^\infty (\Bbb{R)}%
\otimes \mathcal{M},L_c^2(\widetilde{\Gamma }))}\left\| \nabla {\varphi}
(x+t,y)\right\| _{L^p(L^\infty (\Bbb{R)}\otimes \mathcal{M},L_c^2(\widetilde{%
\Gamma }))} \\
&\leq &\left\| g\right\| _{\mathcal{H}_c^{p^{\prime }}}\left\| {\varphi}
\right\| _{\mathcal{H}_c^p}.
\end{eqnarray*}
Then by Theorem 4.4 
\begin{equation}
\left\| {\varphi} \right\| _{\mathrm{BMO}_c^p}\leq c_p\sup_{\left\|
g\right\| _{\mathcal{H}_c^{p^{\prime }}}\leq 1}|\tau \int g{\varphi}
^{*}dt|\leq c_p\left\| {\varphi} \right\| _{\mathcal{H}_c^p}.  \label{bxh}
\end{equation}
To prove the converse, we consider the following tent space $T_c^p.$ Denote $%
\widetilde{\Bbb{R}_{+}^2}=(\Bbb{R}_{+}^2,\frac{dxdy}{y^2})\times
(\{1,2\},\sigma )$ with $\sigma \{1\}=\sigma \{2\}=1.$ For $f\in L^p(%
\mathcal{M},L_c^2(\widetilde{\Bbb{R}_{+}^2})),$ set 
\[
A_c(f)(t)=(\mathop{\displaystyle\int\!\!\!\int}_{|x|<y}|f(x+t,y)|^2dx\frac{dy%
}{y^2})^{\frac 12}. 
\]
Define, for $1<p<\infty ,$ 
\begin{equation}
T_c^p=\{f\in L^p(\mathcal{M},L_c^2(\widetilde{\Bbb{R}_{+}^2})),\left\|
f\right\| _{T_c^p}=\left\| A_c(f)\right\| _{L^p(L^\infty (\Bbb{R})\otimes 
\mathcal{M})}<\infty \}.  \label{tcp}
\end{equation}
We will prove that, for $p>2$ and ${\varphi} $ $\in \mathrm{BMO}_c^p(\Bbb{R},%
\mathcal{M}),\,{\varphi} $ induces a linear functional on $T_c^{p^{\prime }}$
defined by 
\[
l_{\varphi} (f)=\tau {\int \int }_{\Bbb{R}_{+}^2}\nabla {\varphi}
^{*}(x,y)yf(x,y)dxdy/y 
\]
and 
\begin{equation}
\left\| {\varphi} \right\| _{\mathcal{H}_c^p}\leq c_p\left\| l_{\varphi}
\right\| \leq c_p\left\| {\varphi} \right\| _{\mathrm{BMO}_c^p}.
\label{lhbmo}
\end{equation}
We first prove the second inequality of (\ref{lhbmo}). Set 
\begin{eqnarray*}
A_c(f)(t,y) &=&(\mathop{\displaystyle\int\!\!\!\int}%
_{s>y,|x|<s-y}|f(x+t,s)|^2dx\frac{ds}{s^2})^{\frac 12} \\
\overline{A}_c(f)(t,y) &=&(\mathop{\displaystyle\int\!\!\!\int}%
_{s>y,|x|<\frac s4}|f(x+t,s)|^2dx\frac{ds}{s^2})^{\frac 12}.
\end{eqnarray*}
It is easy to see that 
\begin{eqnarray}
\overline{A}_c^2(f)(t,y) &\leq &\overline{A}_c^2(f)(t,0)\leq A_c^2(f)(t),
\label{Ac1} \\
\overline{A}_c^2(f)(t+x,y) &\leq &A_c^2(f)(t,\frac y2),\quad \forall
|x|<\frac y4,(t,y)\in \Bbb{R}_{+}^2.  \label{Ac}
\end{eqnarray}
For nice $f$ and by approximation, we can assume $A_c(f)(t,y)$ is invertible
for all $(t,y)\in \Bbb{R}_{+}^2.$ Thus by Cauchy-Schwartz inequality 
\begin{eqnarray*}
l_{\varphi} (f) &=&\tau {\int \int }_{\Bbb{R}_{+}^2}f(t,y)\nabla {\varphi}
^{*}(t,y)ydt\frac{dy}y \\
&\leq &(\tau \mathop{\displaystyle\int\!\!\!\int}_{\Bbb{R}%
_{+}^2}A_c^{p^{\prime }-2}(f)(t,\frac y2)|f|^2ydt\frac{dy}{y^2})^{\frac
12}(\tau \mathop{\displaystyle\int\!\!\!\int}_{\Bbb{R}_{+}^2}A_c^{2-p^{%
\prime }}(f)(t,\frac y2)|\nabla {\varphi} |^2ydtdy)^{\frac 12} \\
&=&I\cdot II
\end{eqnarray*}
Similarly to the proof of Theorem 4.4, we have 
\[
II^2\leq c\left\| {\varphi} \right\| _{\mathrm{BMO}_c^p}^2\left\| f\right\|
_{T_c^{p^{\prime }}}^{2-p^{\prime }} 
\]
Concerning the factor $I,$ by (\ref{Ac}) we have (recall $p^{\prime }-2<0$) 
\begin{eqnarray*}
I^2 &\leq &\tau \mathop{\displaystyle\int\!\!\!\int}_{\Bbb{R}%
_{+}^2}2\int_{t-\frac y4}^{t+\frac y4}\overline{A}_c^{p^{\prime
}-2}(f)(x,y)dx|f(t,y)|^2dt\frac{dy}{y^2} \\
&\leq &2\tau \mathop{\displaystyle\int\!\!\!\int}_{\Bbb{R}_{+}^2}\overline{A}%
_c^{p^{\prime }-2}(f)(x,y)\int_{x-\frac y4}^{x+\frac y4}|f(t,y)|^2dtdx\frac{%
dy}{y^2} \\
&\leq &-2\tau \mathop{\displaystyle\int\!\!\!\int}_{\Bbb{R}_{+}^2}\overline{A%
}_c^{p^{\prime }-2}(f)(x,y)\frac{\partial \overline{A}_c^2(f)}{\partial y}%
(x,y)dydx \\
&=&-4\tau \mathop{\displaystyle\int\!\!\!\int}_{\Bbb{R}_{+}^2}\overline{A}%
_c^{p^{\prime }-1}(f)(x,y)\frac{\partial \overline{A}_c(f)}{\partial y}%
(x,y)dydx \\
&\leq &-4\tau \int_{\Bbb{R}}\overline{A}_c^{p^{\prime }-1}(f)(x,0)\int_{\Bbb{%
R}^{+}}\frac{\partial \overline{A}_c(f)}{\partial y}(x,y)dydx \\
&\leq &4\left\| f\right\| _{T_c^{p^{\prime }}}^{p^{\prime }}
\end{eqnarray*}
Thus 
\begin{equation}
\left\| l_{\varphi} \right\| \leq c\left\| {\varphi} \right\| _{\mathrm{BMO}%
_c^p}.  \label{lbmo}
\end{equation}
Next we prove that $\left\| {\varphi} \right\| _{\mathcal{H}_c^p}\leq
c_p\left\| l_{\varphi} \right\| .$ Since we can regard $T_c^{p^{\prime }}$
as a closed subspace of $L^{p^{\prime }}(L^\infty (\Bbb{R)}\otimes \mathcal{M%
},L_c^2(\widetilde{\Bbb{R}_{+}^2}))$ via the map\thinspace $%
f(x,y)\rightarrow f(x,y)\chi _{\{|x-t|<y\}}.$ $\,l_{\varphi} $ extends to a
linear functional on $L^{p^{\prime }}(L^\infty (\Bbb{R)}\otimes \mathcal{M}%
,L_c^2(\widetilde{\Bbb{R}_{+}^2}))$ with the same norm. Then there exists $%
h\in L^p(L^\infty (\Bbb{R)}\otimes \mathcal{M},L_c^2(\widetilde{\Bbb{R}_{+}^2%
}))\,$such that $\left\| h\right\| _{L^p(L^\infty (\Bbb{R)}\otimes \mathcal{M%
},L_c^2(\widetilde{\Bbb{R}_{+}^2}))}$ $\leq \left\| l_{\varphi} \right\| \,$%
and 
\begin{eqnarray*}
l_{\varphi} (f) &=&\tau \int_{\Bbb{R}}\mathop{\displaystyle\int\!\!\!\int}%
_{|x-t|<y}f(x,y)h^{*}(x,y,t)dx\frac{dy}{y^2}dt \\
&=&\tau \mathop{\displaystyle\int\!\!\!\int}_{\Bbb{R}_{+}^2}f(x,y)%
\int_{x-y}^{x+y}h^{*}(x,y,t)dtdx\frac{dy}{y^2}.
\end{eqnarray*}
for every $f(x,y)\in T_c^{p^{\prime }}.$ Thus 
\begin{equation}
\nabla {\varphi} (x,y)y=\frac 1y\int_{x-y}^{x+y}h(x,y,t)dt.  \label{pfai}
\end{equation}
Then 
\begin{eqnarray*}
\left\| {\varphi} \right\| _{\mathcal{H}_c^p}^2 &=&(\tau \int_{\Bbb{R}}(%
\mathop{\displaystyle\int\!\!\!\int}_{_\Gamma }|\frac
1y\int_{x+s-y}^{x+s+y}h(x+s,y,t)dt|^2dx\frac{dy}{y^2})^{\frac p2}ds)^{\frac
2p} \\
&\leq &(\tau \int_{\Bbb{R}}(\mathop{\displaystyle\int\!\!\!\int}_{_{\Bbb{R}%
_{+}^2}}\frac 1y\int_{s-2y}^{s+2y}|h(x,y,t)|^2dtdx\frac{dy}{y^2})^{\frac
p2}ds)^{\frac 2p} \\
&=&\left\| \mathop{\displaystyle\int\!\!\!\int}_{_{\Bbb{R}_{+}^2}}\frac
1y\int_{s-2y}^{s+2y}|h(x,y,t)|^2dtdx\frac{dy}{y^2}\right\| _{^{L^{\frac
p2}(L^\infty (\Bbb{R})\otimes \mathcal{M})}}
\end{eqnarray*}
Notice that, for every positive $a$ with $\left\| a\right\| _{^{L^{(\frac
p2)^{{\prime }}}(L^\infty (\Bbb{R})\otimes \mathcal{M})}}\leq 1,$ by (\ref
{MAXC}) and (\ref{junge}) we have 
\begin{eqnarray*}
&&\tau \int_{\Bbb{R}}\mathop{\displaystyle\int\!\!\!\int}_{\Bbb{R}%
_{+}^2}\frac 1y\int_{s-2y}^{s+2y}|h(x,y,t)|^2dtdx\frac{dy}{y^2}a(s)ds \\
&=&\tau \int_{\Bbb{R}}\mathop{\displaystyle\int\!\!\!\int}_{\Bbb{R}%
_{+}^2}|h(x,y,t)|^2\frac 1y\int_{t-2y}^{t+2y}a(s)dsdx\frac{dy}{y^2}dt \\
&\leq &8\tau \int_{\Bbb{R}}\sum_{n=-\infty }^{+\infty
}\int_{2^{n-2}}^{2^{n-1}}\int_{\Bbb{R}}|h(x,y,t)|^2dx\frac{dy}{y^2}\frac
1{2^{n+1}}\int_{t-2^{n}}^{t+2^{n}}\!a(s)dsdt \\
&\leq &8\left\| \mathop{\displaystyle\int\!\!\!\int}_{\Bbb{R}%
_{+}^2}|h(x,y,t)|^2dx\frac{dy}{y^2}\right\| _{^{L^{\frac p2}(L^\infty (\Bbb{R%
})\otimes \mathcal{M})}}\left\| \sup_n|\frac
1{2^{n+1}}\int_{t-2^{n}}^{t+2^{n}}\!a(s)ds|\right\| _{^{L^{(\frac p2)^{{%
\prime }}}(L^\infty (\Bbb{R})\otimes \mathcal{M})}} \\
&\leq &c_p\left\| h\right\| _{L^p(L^\infty (\Bbb{R)}\otimes \mathcal{M}%
,L_c^2(\widetilde{\Bbb{R}_{+}^2}))}^2\leq c_p\left\| l_{\varphi} \right\| ^2
\end{eqnarray*}
Therefore by taking the supremum over all $a$ as above, we obtain 
\[
\left\| {\varphi} \right\| _{\mathcal{H}_c^p}^2\leq c_p\left\| l_{\varphi}
\right\| ^2 
\]
Combining this with (\ref{lbmo}) we get 
\[
\left\| {\varphi} \right\| _{\mathcal{H}_c^p}\leq c_p\left\| {\varphi}
\right\| _{^{\mathrm{BMO}_c^p}}. 
\]
And then $\left\| {\varphi} \right\| _{\mathcal{H}_c^p}\backsimeq \left\| {%
\varphi} \right\| _{^{\mathrm{BMO}_c^p}}$ for every ${\varphi} \in \mathcal{H%
}_c^p(\Bbb{R},\mathcal{M}).$

To prove $\mathrm{BMO}_c^p(\Bbb{R},\mathcal{M})$ and $\mathcal{H}_c^p(\Bbb{R}%
,\mathcal{M})$ are the same space, it remains to show that the family of $S_{%
\mathcal{M}}$-simple functions is dense in $\mathrm{BMO}_c^p(\Bbb{R},%
\mathcal{M}).$ From the proof of Theorem 4.4 we can see that for every ${%
\varphi} \in \mathrm{BMO}_c^p(\Bbb{R},\mathcal{M}),$ there exists a $h\in
L^{^\infty }(L^{^\infty }(\Bbb{R)}\otimes \mathcal{M},L_c^2)$ such that ${%
\varphi} =\Psi (h)$ and $\left\| \Psi (h)\right\| _{\mathrm{BMO}_c^p}\leq
c\left\| h\right\| _{L^p(L^{^\infty }(\Bbb{R)}\otimes \mathcal{M},L_c^2)}.$
Recall that the family of ''nice'' $h$'s(i.e. $h(x,y,t)=%
\sum_{i=1}^nm_if_i(t)\chi _{A_i}$ with $m_i\in S_{\mathcal{M}},A_i\in 
\widetilde{\Gamma },|A_i|<\infty $ and with scalar valued simple functions $%
f_i$) is dense in $L^p(L^\infty (\Bbb{R)}\otimes \mathcal{M},L_c^2).$
Choose ''nice'' $h_n\rightarrow h$ in $L^p(L^{^\infty }(\Bbb{R)}\otimes 
\mathcal{M},L_c^2).$ Let ${\varphi} _n=\Psi (h_n).$ Then ${\varphi}
_n\rightarrow {\varphi} $ in $\mathrm{BMO}_c^p(\Bbb{R},\mathcal{M}).$ Since
the ${\varphi} _n$'s$\,$ are continuous functions with compact support, we
can approximate them by simple functions in $\mathrm{BMO}_c^p(\Bbb{R},%
\mathcal{M}).$ This shows the density of simple functions in $\mathrm{BMO}%
_c^p(\Bbb{R},\mathcal{M})$ and thus completes the proof of the theorem.\qed

\medskip 
\noindent\textbf{Remark. }$ $By the same idea used in the
proof above, we can get the analogue of the classical duality result for the
tent spaces: $(T_c^p)^{*}=T_c^q$ $(1<p<\infty )$ with equivalent norms, where $%
T_c^p$ is defined as (\ref{tcp}).

\begin{theorem}
(i)$\quad \Psi $ extends to a bounded map from $L^\infty (L^\infty (\Bbb{R)}%
\otimes \mathcal{M},L_c^2(\widetilde{\Gamma }))\ $ into $\mathrm{BMO}_c(\Bbb{%
R},\mathcal{M})$ and 
\begin{equation}
\left\| \Psi (h)\right\| _{\mathrm{BMO}_c}\leq c\left\| h\right\|
_{L^{^\infty }(L^{^\infty }(\Bbb{R)}\otimes \mathcal{M},L_c^2)}  \label{psi}
\end{equation}

(ii)$\quad \Psi $ extends to a bounded map from $L^p(L^\infty (\Bbb{R)}%
\otimes \mathcal{M},L_c^2(\widetilde{\Gamma }))\ $into $\mathcal{H}_c^p(\Bbb{%
R},\mathcal{M})$ $(1<p<\infty )$and 
\begin{equation}
\left\| \Psi (h)\right\| _{\mathcal{H}_c^p}\leq c_p\left\| h\right\|
_{L^p(L^{^\infty }(\Bbb{R)}\otimes \mathcal{M},L_c^2)}.  \label{psip}
\end{equation}

(iii)$\quad $The statements (i) and (ii) also hold with column spaces
replaced by row spaces.
\end{theorem}

\noindent\textbf{Proof.} (\ref{psi}) is Lemma 2.2. The part of (\ref{psip})
concerning $p>2$ follows from Lemma 4.3 and Theorem 4.7. For $1<p<2,\,$ by
the duality between $\,\mathcal{H}_c^p$ and $\mathrm{BMO}_c^q,$ and Theorem
4.7, we have 
\begin{eqnarray}
\left\| \Psi (h)\right\| _{\mathcal{H}_c^p} &\leq &c\sup_{\left\| f\right\|
_{\mathrm{BMO}_c^q}\leq 1}\left| \tau \int_{\Bbb{R}}\Psi
(h)(s)f^{*}(s)ds\right|  \nonumber \\
&\le &\sup_{\left\| f\right\| _{\mathcal{H}_c^q}\leq 1}\left| \tau \int_{%
\Bbb{R}}\int_{\Bbb{R}}\mathop{\displaystyle\int\!\!\!\int}_\Gamma
h(x,y,t)\nabla P_y(x+t-s)dxdydtf^{*}(s)ds\right|  \nonumber \\
&=&\sup_{\left\| f\right\| _{\mathcal{H}_c^q}\leq c}\left| \tau \int_{\Bbb{R}%
}\mathop{\displaystyle\int\!\!\!\int}_\Gamma h(x,y,t)\nabla
f^{*}(x+t,y)dxdydt\right|  \nonumber \\
&\leq &c\left\| h\right\| _{L^p(L^\infty (\Bbb{R)\otimes \mathcal{M}},L_c^2%
\Bbb{)}}.  \label{hlp}
\end{eqnarray}
When $p=2,\,$similarly but taking Supremum over $\left\| f\right\| _{%
\mathcal{H}_c^2}\leq 1$ in the formula above, we have$\left\| \Psi
(h)\right\| _{\mathcal{H}_c^2}\leq \left\| h\right\| _{L^2(L^{^\infty }(\Bbb{%
R)}\otimes \mathcal{M},L_c^2)}.$

\begin{corollary}
$(\mathcal{H}_c^p(\Bbb{R},\mathcal{M}))^{*}=\mathcal{H}_c^q(\Bbb{R},\mathcal{%
M})$ with equivalent norms for all $1<p<\infty .\medskip $
\end{corollary}

\chapter{Reduction of BMO to dyadic BMO}

\setcounter{theorem}{0} \setcounter{equation}{0} Our approach in Chapter 3
towards the maximal inequality is to reduce it to the corresponding maximal
inequality for dyadic martingales. In this chapter, we pursue this idea. We
will see that \textrm{BMO} spaces can be characterized as intersections of
dyadic \textrm{BMO}. This result has many consequences. It will be used in
the next chapter for interpolation too.

\section{\textrm{BMO} is the intersection of two dyadic \textrm{BMO}}

Consider an increasing family of $\sigma $-algebras $\mathcal{F=}\{\mathcal{F%
}_n\}_{n\in \Bbb{Z}}\,$on $\Bbb{R}.$ Assume that each $\mathcal{F}_n$ is
generated by a sequence of atoms $\{F_n^k\}_{k\in \Bbb{Z}}.$ We are going to
introduce the \textrm{BMO}$^q$ spaces for martingales with respect to $%
\mathcal{F=}\{\mathcal{F}_n\}_{n\in \Bbb{Z}}.$ Let $2<q\leq \infty $ and ${%
\varphi }\in L^q(\mathcal{M},L_c^2(\Bbb{R},\frac{dt}{1+t^2})).$ Define 
\[
{\varphi }_{\mathcal{F}_n}^{\#}(t)=\frac 1{|F_n^k|}\int_{F_n^k\ni t}|{%
\varphi }(x)-{\varphi }_{F_n^k}|^2dx
\]
For ${\varphi }\in L^q(\mathcal{M},L_c^2(\Bbb{R},\frac{dt}{1+t^2}))($resp. $%
L^q(\mathcal{M},L_r^2(\Bbb{R},\frac{dt}{1+t^2}))),$ let 
\[
\left\| {\varphi }\right\| _{\mathrm{BMO}_c^{q,\mathcal{F}}}=\left\| \sup_n|{%
\varphi }_{\mathcal{F}_n}^{\#}|\right\| _{\frac q2}^{\frac 12}\quad \text{%
and \quad }\left\| {\varphi }\right\| _{\mathrm{BMO}_r^{q,\mathcal{F}%
}}=\left\| {\varphi }^{*}\right\| _{\mathrm{BMO}_c^{q,\mathcal{F}}}.
\]
And set
\begin{eqnarray*}
\mathrm{BMO}_c^{q,\mathcal{F}}(L^\infty (\Bbb{R)}\otimes \mathcal{M}) &=&\{{{%
\varphi }\in L^q(\mathcal{M},}L{_c^2(\Bbb{R},\frac{dt}{1+t^2})),~}\left\| {%
\varphi }\right\| _{\mathrm{BMO}_c^{q,\mathcal{F}}}<\infty \}, \\
\mathrm{BMO}_r^{q,\mathcal{F}}(L^\infty (\Bbb{R)}\otimes \mathcal{M}) &=&\{{%
\varphi }\in L^q(\mathcal{M},L_r^2(\Bbb{R},\frac{dt}{1+t^2})){,~}\left\| {%
\varphi }\right\| _{\mathrm{BMO}_r^{q,\mathcal{F}}}<\infty \}.
\end{eqnarray*}
Define $\mathrm{BMO}_{cr}^{q,\mathcal{F}}$\textrm{\ }to be the intersection
of $\mathrm{BMO}_c^{q,\mathcal{F}}$ and $\mathrm{BMO}_r^{q,\mathcal{F}}$
with the intersection norm $\max \{\left\| {\varphi }\right\| _{\mathrm{BMO}%
_c^{q,\mathcal{F}}},\left\| {\varphi }\right\| _{\mathrm{BMO}_r^{q,\mathcal{F%
}}}\}.$ These \textrm{BMO}$^q$ spaces were already studied in \cite{JX} for
general non-commutative martingales.

In the following, we will consider the spaces $\mathrm{BMO}_c^{q,\mathcal{D}%
}(L^\infty (\Bbb{R)}\otimes \mathcal{M}),$ $\mathrm{BMO}_c^{q,\mathcal{D}%
^{\prime }}(L^\infty (\Bbb{R)}\otimes \mathcal{M}),$ $\mathrm{BMO}_r^{q,%
\mathcal{D}}(L^\infty (\Bbb{R)}\otimes \mathcal{M}),$ $\mathrm{BMO}_r^{q,%
\mathcal{D}^{\prime }}(L^\infty (\Bbb{R)}\otimes \mathcal{M})$  etc. with
respect to the families $\mathcal{D},\mathcal{D}^{^{\prime }}$ of dyadic $%
\sigma $-algebras defined in Chapter 3.   

\begin{theorem}
Let $2<q\leq \infty .$ With equivalent norms, 
\begin{eqnarray*}
\mathrm{BMO}_c^q(\Bbb{R},\mathcal{M}) &=&\mathrm{BMO}_c^{q,\mathcal{D}%
}(L^\infty (\Bbb{R)}\otimes \mathcal{M})\cap \mathrm{BMO}_c^{q,\mathcal{D}%
^{\prime }}(L^\infty (\Bbb{R)}\otimes \mathcal{M});\mathrm{\ } \\
\mathrm{BMO}_r^q(\Bbb{R},\mathcal{M}) &=&\mathrm{BMO}_r^{q,\mathcal{D}%
}(L^\infty (\Bbb{R)}\otimes \mathcal{M})\cap \mathrm{BMO}_r^{q,\mathcal{D}%
^{\prime }}(L^\infty (\Bbb{R)}\otimes \mathcal{M}); \\
\mathrm{BMO}_{cr}^q(\Bbb{R},\mathcal{M}) &=&\mathrm{BMO}_{cr}^{q,\mathcal{D}%
}(L^\infty (\Bbb{R)}\otimes \mathcal{M})\cap \mathrm{BMO}_{cr}^{q,\mathcal{D}%
^{\prime }}(L^\infty (\Bbb{R)}\otimes \mathcal{M}).
\end{eqnarray*}
\end{theorem}

\noindent\textbf{Proof.} From Proposition 3.1, $\forall t\in \Bbb{R},h\in 
\Bbb{R}^{+}\times \Bbb{R}^{+},$ there exist $k_{t,h},N_h\in \Bbb{Z}$ such
that $I_{h,t}:=(t-h_1,t+h_2]$ is contained in $D_{N_h}^{k_{t,h}}$ or $%
D_{N_h}^{\prime k_{t,h}}$ and 
\[
|D_{N_h}^{k_{t,h}}|=|D_{N_h}^{\prime k_{t,h}}|\leq 6(h_1+h_2).
\]
If $I_{h,t}\subset D_{N_h}^{k_{t,h}},$ then 
\begin{eqnarray*}
{\varphi }_h^{\#}(t) &=&\frac 1{h_1+h_2}\int_{t-h_1}^{t+h_2}|{\varphi }(x)-{%
\varphi }_{I_{h,t}}|^2dx \\
&\leq &\frac 4{h_1+h_2}\int_{t-h_1}^{t+h_2}|{\varphi }(x)-{\varphi }%
_{D_{N_h}^{k_{t,h}}}|^2dx \\
&\leq &\frac{24}{|D_{N_h}^{k_{t,h}}|}\int_{D_{N_h}^{k_{t,h}}}|{\varphi }(x)-{%
\varphi }_{D_{N_h}^{k_{t,h}}}|^2dx \\
&\leq &24{\varphi }_{\mathcal{D}_{N_h}}^{\#}(t).
\end{eqnarray*}
Similarly, if $I_{h,t}\subset D_{N_h}^{^{\prime }k_{t,h}},$ then 
\[
{\varphi }_h^{\#}(t)\leq 24{\varphi }_{\mathcal{D}_{N_h}^{^{\prime
}}}^{\#}(t).
\]
Thus 
\begin{eqnarray*}
\left\| {\varphi }\right\| _{\mathrm{BMO}_c^q} &=&\left\| \sup_{h\in \Bbb{R}%
^{+}\times \Bbb{R}^{+}}|{\varphi }_h^{\#}|\right\| _{\frac q2}^{\frac 12} \\
&\leq &\sqrt{24}\left\| \sup_n|({\varphi }_{\mathcal{D}_n}^{\#}+{\varphi }_{%
\mathcal{D}_n^{^{\prime }}}^{\#})|\right\| _{\frac q2}^{\frac 12} \\
&\leq &4\sqrt{3}\mathrm{\max }(\left\| {\varphi }\right\| _{\mathrm{BMO}%
_c^{q,\mathcal{D}}},\left\| {\varphi }\right\| _{\mathrm{BMO}_c^{q,\mathcal{D%
}^{^{\prime }}}}).
\end{eqnarray*}
It is trivial that $\mathrm{\max }(\left\| {\varphi }\right\| _{\mathrm{BMO}%
_c^{q,\mathcal{D}}},\left\| {\varphi }\right\| _{\mathrm{BMO}_c^{q,\mathcal{D%
}^{^{\prime }}}})\leq \left\| {\varphi }\right\| _{\mathrm{BMO}_c^q}.$
Therefore 
\[
\mathrm{BMO}_c^q(\Bbb{R},\mathcal{M})=\mathrm{BMO}_c^{q,\mathcal{D}%
}(L^\infty (\Bbb{R)}\otimes \mathcal{M)}\cap \mathrm{BMO}_c^{q,\mathcal{D}%
^{^{\prime }}}(L^\infty (\Bbb{R)}\otimes \mathcal{M})
\]
with equivalent norms. The two other equalities in the theorem are immediate
consequences of this.\qed

\section{The equivalence of $\mathcal{H}_{cr}^p(\Bbb{R},\mathcal{M}\Bbb{)}$
and $L^p(L^\infty (\Bbb{R)}\otimes \mathcal{M}\Bbb{)(}1\Bbb{<}p\Bbb{<\infty )%
}$}

We denote the non-commutative martingale Hardy spaces defined in \cite{[PX2]}
and \cite{JX} with respect to $\mathcal{D}$ and $\mathcal{D}^{^{\prime }}\,$%
by $\mathcal{H}_c^{p,\mathcal{D}}(L^\infty (\Bbb{R)}\otimes \mathcal{M}),%
\mathcal{H}_c^{p,\mathcal{D}^{\prime }}(L^\infty (\Bbb{R)}\otimes \mathcal{M}%
)\ $etc.$(1\leq p<\infty )$. Note that 
\[
\mathcal{H}_c^2(\Bbb{R},\mathcal{M})=\mathcal{H}_c^{2,\mathcal{D}}(L^\infty (%
\Bbb{R)}\otimes \mathcal{M})=\mathcal{H}_c^{2,\mathcal{D}^{\prime
}}(L^\infty (\Bbb{R)}\otimes \mathcal{M})=L^2(L^\infty (\Bbb{R)}\otimes 
\mathcal{M}).
\]
By Theorems 4.4, 5.1 and the duality equality $(\mathcal{H}_c^{p,\mathcal{D}%
}(L^\infty (\Bbb{R)}\otimes \mathcal{M}))^{*}=\mathrm{BMO}_c^{q,\mathcal{D}%
}(L^\infty (\Bbb{R)}\otimes \mathcal{M})$  proved in \cite{JX}, \cite{JX} we
get the following result.\smallskip 

\begin{corollary}
$\mathrm{BMO}_{cr}^q(\Bbb{R},\mathcal{M})=L^q(L^\infty (\Bbb{R})\otimes 
\mathcal{M})$ with equivalent norms for $2<q<\infty .$
\end{corollary}

\noindent\textbf{Proof.} From the inequalities (4.5) and (4.7) of \cite{JX}
we have 
\begin{eqnarray*}
&&\mathrm{BMO}_c^{q,\mathcal{D}}(L^\infty (\Bbb{R)}\otimes \mathcal{M})\cap 
\mathrm{BMO}_r^{q,\mathcal{D}}(L^\infty (\Bbb{R)}\otimes \mathcal{M}) \\
&=&L^q(L^\infty (\Bbb{R})\otimes \mathcal{M}) \\
&=&\mathrm{BMO}_c^{q,\mathcal{D}^{\prime }}(L^\infty (\Bbb{R)}\otimes 
\mathcal{M})\cap \mathrm{BMO}_r^{q,\mathcal{D}^{\prime }}(L^\infty (\Bbb{R)}%
\otimes \mathcal{M})\quad 
\end{eqnarray*}
with equivalent norms. $\,$Therefore, by Theorem 5.1 
\begin{eqnarray*}
&&\mathrm{BMO}_{cr}^q(\Bbb{R},\mathcal{M}) \\
&=&\mathrm{BMO}_c^q(\Bbb{R},\mathcal{M})\cap \mathrm{BMO}_r^q(\Bbb{R},%
\mathcal{M}) \\
&=&\mathrm{BMO}_c^{q,\mathcal{D}}(L^\infty (\Bbb{R)}\otimes \mathcal{M})\cap 
\mathrm{BMO}_r^{q,\mathcal{D}}(L^\infty (\Bbb{R)}\otimes \mathcal{M}) \\
&&\qquad \cap \mathrm{BMO}_c^{q,\mathcal{D}^{\prime }}(L^\infty (\Bbb{R)}%
\otimes \mathcal{M})\cap \mathrm{BMO}_r^{q,\mathcal{D}^{\prime }}(L^\infty (%
\Bbb{R)}\otimes \mathcal{M}) \\
&=&L^q(L^\infty (\Bbb{R})\otimes \mathcal{M}).\qed
\end{eqnarray*}

\begin{corollary}
If $1\leq p<2,$ then 
\begin{eqnarray*}
\mathcal{H}_c^p(\Bbb{R},\mathcal{M}) &=&\mathcal{H}_c^{p,\mathcal{D}%
}(L^\infty (\Bbb{R)}\otimes \mathcal{M})+\mathcal{H}_c^{p,\mathcal{D}%
^{\prime }}(L^\infty (\Bbb{R)}\otimes \mathcal{M}), \\
\mathcal{H}_r^p(\Bbb{R},\mathcal{M}) &=&\mathcal{H}_r^{p,\mathcal{D}%
}(L^\infty (\Bbb{R)}\otimes \mathcal{M})+\mathcal{H}_r^{p,\mathcal{D}%
^{\prime }}(L^\infty (\Bbb{R)}\otimes \mathcal{M}), \\
\mathcal{H}_{cr}^p(\Bbb{R},\mathcal{M}) &=&\mathcal{H}_{cr}^{p,\mathcal{D}%
}(L^\infty (\Bbb{R)}\otimes \mathcal{M})+\mathcal{H}_{cr}^{p,\mathcal{D}%
^{\prime }}(L^\infty (\Bbb{R)}\otimes \mathcal{M}).
\end{eqnarray*}
If $p\geq 2,$ then 
\begin{eqnarray*}
\mathcal{H}_c^p(\Bbb{R},\mathcal{M}) &=&\mathcal{H}_c^{p,\mathcal{D}%
}(L^\infty (\Bbb{R)}\otimes \mathcal{M})\cap \mathcal{H}_c^{p,\mathcal{D}%
^{\prime }}(L^\infty (\Bbb{R)}\otimes \mathcal{M}), \\
\mathcal{H}_r^p(\Bbb{R},\mathcal{M}) &=&\mathcal{H}_r^{p,\mathcal{D}%
}(L^\infty (\Bbb{R)}\otimes \mathcal{M})\cap \mathcal{H}_r^{p,\mathcal{D}%
^{\prime }}(L^\infty (\Bbb{R)}\otimes \mathcal{M}), \\
\mathcal{H}_{cr}^p(\Bbb{R},\mathcal{M}) &=&\mathcal{H}_{cr}^{p,\mathcal{D}%
}(L^\infty (\Bbb{R)}\otimes \mathcal{M})\cap \mathcal{H}_{cr}^{p,\mathcal{D}%
^{\prime }}(L^\infty (\Bbb{R)}\otimes \mathcal{M}).
\end{eqnarray*}
\end{corollary}

\begin{corollary}
$\mathcal{H}_{cr}^p(\Bbb{R},\mathcal{M})=L^p(L^\infty (\Bbb{R})\otimes 
\mathcal{M})$ with equivalent norms for all $1<p<\infty .$
\end{corollary}

\noindent\textbf{Proof.} Recall the result 
\[
\mathcal{H}_{cr}^{p,\mathcal{D}}(L^\infty (\Bbb{R)}\otimes \mathcal{M})=L^p(%
\Bbb{R},\mathcal{M})=\mathcal{H}_{cr}^{p,\mathcal{D}^{\prime }}(L^\infty (%
\Bbb{R)}\otimes \mathcal{M})
\]
proved in \cite{[PX2]} and \cite{JX}. By Corollary 5.3, for $1<p<2,$ we have 
\begin{eqnarray*}
\mathcal{H}_{cr}^p(\Bbb{R},\mathcal{M}) &=&\mathcal{H}_c^p(\Bbb{R},\mathcal{M%
})+\mathcal{H}_r^p(\Bbb{R},\mathcal{M}) \\
&=&\mathcal{H}_c^{p,\mathcal{D}}(L^\infty (\Bbb{R)}\otimes \mathcal{M})+%
\mathcal{H}_c^{p,\mathcal{D}^{\prime }}(L^\infty (\Bbb{R)}\otimes \mathcal{M}%
) \\
&&\qquad +\mathcal{H}_r^{p,\mathcal{D}}(L^\infty (\Bbb{R)}\otimes \mathcal{M}%
)+\mathcal{H}_r^{p,\mathcal{D}^{\prime }}(L^\infty (\Bbb{R)}\otimes \mathcal{%
M}) \\
&=&\mathcal{H}_{cr}^{p,\mathcal{D}}(L^\infty (\Bbb{R)}\otimes \mathcal{M})+%
\mathcal{H}_{cr}^{p,\mathcal{D}^{\prime }}(L^\infty (\Bbb{R)}\otimes 
\mathcal{M}) \\
&=&L^p(L^\infty (\Bbb{R)}\otimes \mathcal{M})
\end{eqnarray*}
and, for $2\leq p<\infty ,$ 
\begin{eqnarray*}
\mathcal{H}_{cr}^p(\Bbb{R},\mathcal{M}) &=&\mathcal{H}_c^p(\Bbb{R},\mathcal{M%
})\cap \mathcal{H}_c^p(\Bbb{R},\mathcal{M}) \\
&=&\mathcal{H}_c^{p,\mathcal{D}}(L^\infty (\Bbb{R)}\otimes \mathcal{M})\cap 
\mathcal{H}_c^{p,\mathcal{D}^{\prime }}(L^\infty (\Bbb{R)}\otimes \mathcal{M}%
) \\
&&\qquad \qquad \qquad \cap \mathcal{H}_r^{p,\mathcal{D}}(L^\infty (\Bbb{R)}%
\otimes \mathcal{M})\cap \mathcal{H}_r^{p,\mathcal{D}^{\prime }}(L^\infty (%
\Bbb{R)}\otimes \mathcal{M}) \\
&=&\mathcal{H}_{cr}^{p,\mathcal{D}}(L^\infty (\Bbb{R)}\otimes \mathcal{M}%
)\cap \mathcal{H}_{cr}^{p,\mathcal{D}^{\prime }}(L^\infty (\Bbb{R)}\otimes 
\mathcal{M}) \\
&=&L^p(L^\infty (\Bbb{R)}\otimes \mathcal{M}).\qed
\end{eqnarray*}

\medskip \noindent
\textbf{Remark.} In \cite{JMX1} and \cite{JMX}, M. Junge, C. Le Merdy and Q.
Xu have studied the Littlewood-Paley theory for semigroups on
non-commutative $L^p$-spaces. Among many results, they proved in particular,
that for many nice semigroups, the corresponding non-commutative Hardy
spaces defined by the Littlewood-Paley $g$-function coincide with the
underlying non-commutative $L^p$-spaces ($1<p<\infty $). In their viewpoint,
the semigroup in the context of our paper is the Poisson semigroup
tensorized by the identity of $L^p(\mathcal{M}).$ This semigroup satisfies
all assumptions of \cite{JMX}. Thus if we define our Hardy spaces $\mathcal{H%
}_{cr}^p(\Bbb{R},\mathcal{M})$ by the $g$-function $G_c(f)$ and $G_r(f)$
(which is the same as that defined by $S_c(f)$ and $S_r(f)$ in virtue of
Theorem 4.6), then Corollary 5.4 is a particular case of a general result
from \cite{JMX}. We should emphasize that the method in \cite{JMX} is
completely different from ours. It is based on the $H^\infty $ functional
calculus. It seems that the method in \cite{JMX} does not permit to deal
with the Lusin square functions $S_c(f)$ and $S_r(f).\medskip $

\chapter{Interpolation}

\setcounter{theorem}{0} \setcounter{equation}{0} In this chapter, we
consider the interpolation for non-commutative Hardy spaces and \textrm{BMO. 
}The main results in this chapter are function space analogues of those in 
\cite{MM} for non-commutative martingales. On the other hand, they are also
the extensions to the present non-commutative setting of the scalar results
in \cite{JJ}. Recall that the non-commutative $L^p$ spaces associated with a
semifinite von Neumann algebra form an interpolation scale with respect to
both the complex and real interpolation methods. And, as the column (resp.
row) subspaces of $L^p(\mathcal{M}\otimes B(L^2(\Omega )))\,,$ the spaces $%
L^p(L^{^\infty }(\Bbb{R)}\otimes \mathcal{M},L_c^2(\widetilde{\Gamma }))$
form an interpolation scale also.

\section{The complex interpolation}

We first consider the complex interpolation. 

Let $\mathrm{BMO}_c^{\mathcal{D}}(L^\infty (\Bbb{R)}\otimes \mathcal{M})$
and $\mathcal{H}_c^{p,\mathcal{D}}(L^\infty (\Bbb{R)}\otimes \mathcal{M}%
)$ (resp. $\mathrm{BMO}_c^{\mathcal{D'}}(L^\infty (\Bbb{R)}\otimes \mathcal{M})$
and $\mathcal{H}_c^{p,\mathcal{D'}}(L^\infty (\Bbb{R)}\otimes \mathcal{M}
))$ $(1\leq p<\infty )$ be the non-commutative martingale BMO spaces and Hardy
spaces defined in \cite{JX} with respect to the usual dyadic
filtration $\mathcal{D}$ (resp. the dyadic filtration $\mathcal{D'}$) 
described in Chapter 3. 

\begin{lemma}
For $\,1<p<\infty ,$ we have 
\begin{eqnarray}
(\mathrm{BMO}_c^{\mathcal{D}}(L^\infty (\Bbb{R)}\otimes \mathcal{M}),%
\mathcal{H}_c^{1,\mathcal{D}}(L^\infty (\Bbb{R)}\otimes \mathcal{M}))_{\frac
1p} =&\mathcal{H}_c^{p,\mathcal{D}}(L^\infty (\Bbb{R)}\otimes \mathcal{M}),
\label{dic} \\
(\mathrm{BMO}_r^{\mathcal{D}}(L^\infty (\Bbb{R)}\otimes \mathcal{M}),%
\mathcal{H}_r^{1,\mathcal{D}}(L^\infty (\Bbb{R)}\otimes \mathcal{M}))_{\frac
1p} =&\mathcal{H}_r^{p,\mathcal{D}}(L^\infty (\Bbb{R)}\otimes \mathcal{M}),
\label{dir} \\
(X,Y)_{\frac 1p} =&L^p(L^\infty (\Bbb{R)}\otimes \mathcal{M}).  \label{dicr}
\end{eqnarray}
where $X=\mathrm{BMO}_{cr}^{\mathcal{D}}(L^\infty (\Bbb{R)}\otimes \mathcal{M%
})$ or $L^\infty (L^\infty (\Bbb{R)}\otimes \mathcal{M})$ and $Y=\mathcal{H}%
_{cr}^{1,\mathcal{D}}(L^\infty (\Bbb{R)}\otimes \mathcal{M})\,$or $%
L^1(L^\infty (\Bbb{R)}\otimes \mathcal{M}).$ Moreover, the same results hold for 
$\mathrm{BMO}_c^{\mathcal{D'}}(L^\infty (\Bbb{R)}\otimes \mathcal{M})$
and $\mathcal{H}_c^{p,\mathcal{D'}}(L^\infty (\Bbb{R)}\otimes \mathcal{M}).$
\end{lemma}

\noindent\textbf{Proof. } For each $k\in \Bbb{N}$ and 
each projection $p$ of $\mathcal {M}$ with $\tau (p)<\infty,$ 
denote by $\mathcal{H}_c^{q,\mathcal{D}
}(L^\infty (-2^k,2^k)\otimes p\mathcal{M}p)$ the
subspace of $\mathcal{H}_c^{q,\mathcal{D}}(L^\infty (\Bbb{R)}\otimes \mathcal{M}%
)$ consisting of elements supported on $(-2^k,2^k)$ and with values in
$p\mathcal{M}p$. 
By dualizing Theorem 3.1 of \cite{MM} 
we get, for $1<r\leq q<\infty$,
\begin{eqnarray*}
&&\left(\mathcal{H}_c^{1,\mathcal{D}}(L^\infty (-2^k,2^k)\otimes p\mathcal{M}p),
\mathcal{H}_c^{\frac r{r-1},\mathcal{D}}(L^\infty (-2^k,2^k)\otimes 
p\mathcal{M}p)\right)_{\frac rq} \\
&=&\mathcal{H}_c^{\frac q{q-1},\mathcal{D}}(L^\infty (-2^k,2^k)\otimes
p\mathcal{M}p). 
\end{eqnarray*}
Note that the union of all these $\mathcal{H}_c^{r,\mathcal{D}
}(L^\infty (-2^k,2^k)\otimes p\mathcal{M}p)$ is
dense in $\mathcal{H}_c^{r,\mathcal{D}}(L^\infty (\Bbb{R)}\otimes \mathcal{M}%
)$. By approximation we get
\begin{eqnarray}
(\mathcal{H}_c^{1,\mathcal{D}}(L^\infty (\Bbb{R)}\otimes \mathcal{M}),
\mathcal{H}_c^{\frac r{r-1},\mathcal{D}}(L^\infty (\Bbb{R)}\otimes \mathcal{M}))_{\frac
rq} =\mathcal{H}_c^{\frac q{q-1},\mathcal{D}}(L^\infty (\Bbb{R)}\otimes \mathcal{M})
\label{id} 
\end{eqnarray}
Dualizing (\ref{id}) we have
\begin{equation}
(\mathrm{BMO}_c^{\mathcal{D}}(L^\infty (\Bbb{R)}\otimes \mathcal{M}),%
\mathcal{H}_c^{r,\mathcal{D}}(L^\infty (\Bbb{R)}\otimes \mathcal{M}%
))_{\frac rq}=\mathcal{H}_c^{q,\mathcal{D}}(L^\infty (\Bbb{R)}\otimes 
\mathcal{M}).  \label{id3}
\end{equation}
Combining (\ref{id}) and (\ref{id3}) we get (\ref{dic}) by
Wolff's interpolation theorem (see \cite{w}). The equalities (\ref{dir}),(\ref{dicr}) and the 
arguments for the dyadic filtration $\mathcal{D'}$ can be proved similarly.

\begin{theorem}
Let $1<p<\infty .$ Then with equivalent norms, 
\begin{eqnarray}
(\mathrm{BMO}_c(\Bbb{R},\mathcal{M}),\mathcal{H}_c^1(\Bbb{R},\mathcal{M}%
))_{\frac 1p} &=&\mathcal{H}_c^p(\Bbb{R},\mathcal{M}),  \label{HBMO} \\
(\mathrm{BMO}_r(\Bbb{R},\mathcal{M}),\mathcal{H}_r^1(\Bbb{R},\mathcal{M}%
))_{\frac 1p} &=&\mathcal{H}_r^p(\Bbb{R},\mathcal{M}),  \label{hbmor} \\
(X,Y)_{\frac 1p} &=&L^p(L^\infty (\Bbb{R)}\otimes \mathcal{M}).
\label{hbmocr}
\end{eqnarray}
where $X=\mathrm{BMO}_{cr}(\Bbb{R},\mathcal{M})$ or $L^\infty (L^\infty (%
\Bbb{R)}\otimes \mathcal{M})$ and $Y=\mathcal{H}_{cr}^1(\Bbb{R},\mathcal{M}%
)\,$or $L^1(L^\infty (\Bbb{R)}\otimes \mathcal{M}).$
\end{theorem}

\noindent\textbf{Proof. } Note that 
\[
\mathcal{H}_c^2(\Bbb{R},\mathcal{M})=\mathcal{H}_c^{2,\mathcal{D}}(\Bbb{R},%
\mathcal{M})=\mathcal{H}_c^{2,\mathcal{D}^{\prime }}(\Bbb{R},\mathcal{M}).
\]
$\,$Let $2<q<\infty .$ By Theorem 5.1 and Lemma 6.1 we have 
\begin{eqnarray*}
&&(\mathrm{BMO}_c(\Bbb{R},\mathcal{M}),\mathcal{H}_c^2(\Bbb{R},\mathcal{M}%
))_{\frac 2q} \\
&=&(\mathrm{BMO}_c^{\mathcal{D}}(L^\infty (\Bbb{R)}\otimes \mathcal{M})\cap 
\mathrm{BMO}_c^{\mathcal{D}^{\prime }}(L^\infty (\Bbb{R)}\otimes \mathcal{M}%
),\mathcal{H}_c^2(\Bbb{R},\mathcal{M}))_{\frac 2q} \\
&\subseteq &(\mathrm{BMO}_c^{\mathcal{D}}(L^\infty (\Bbb{R)}\otimes \mathcal{%
M}),\mathcal{H}_c^2(\Bbb{R},\mathcal{M}))_{\frac 2q}\cap (\mathrm{BMO}_c^{%
\mathcal{D}^{\prime }}(L^\infty (\Bbb{R)}\otimes \mathcal{M}),\mathcal{H}%
_c^2(\Bbb{R},\mathcal{M}))_{\frac 2q} \\
&\subseteq &\mathcal{H}_c^{q,\mathcal{D}}(L^\infty (\Bbb{R)}\otimes \mathcal{%
M})\cap \mathcal{H}_c^{q,\mathcal{D}^{^{\prime }}}(L^\infty (\Bbb{R)}\otimes 
\mathcal{M}) \\
&=&\mathcal{H}_c^q(\Bbb{R},\mathcal{M}).
\end{eqnarray*}
Then by duality 
\begin{equation}
(\mathcal{H}_c^1(\Bbb{R},\mathcal{M}),\mathcal{H}_c^2(\Bbb{R},\mathcal{M}%
))_{\frac 2q}\supseteq \mathcal{H}_c^{q^{\prime }}(\Bbb{R},\mathcal{M}).
\label{p'q'}
\end{equation}
The converse of (\ref{p'q'}) can be easily proved since the map $\Phi $
defined by $\Phi (f)=\nabla f(x+t,y)\chi _\Gamma (x,y)$ is isometric from$\,%
\mathcal{H}_c^{q^{\prime }}(\Bbb{R},\mathcal{M})$ to $L^{q^{\prime
}}(L^\infty (\Bbb{R})\otimes \mathcal{M},L_c^2(\widetilde{\Gamma }))$ for $%
q\geq 1.$ Thus we have 
\begin{equation}
(\mathcal{H}_c^1(\Bbb{R},\mathcal{M}),\mathcal{H}_c^2(\Bbb{R},\mathcal{M}%
))_{\frac 2q}=\mathcal{H}_c^{q^{\prime }}(\Bbb{R},\mathcal{M}).  \label{bh}
\end{equation}
Dualizing this equality once more, we get 
\begin{equation}
(\mathrm{BMO}_c(\Bbb{R},\mathcal{M}),\mathcal{H}_c^2(\Bbb{R},\mathcal{M}%
))_{\frac 2q}=\mathcal{H}_c^q(\Bbb{R},\mathcal{M}).  \label{hh}
\end{equation}
Note that by Proposition 2.1 and Theorem 4.8, $\mathcal{H}_c^q$ is
complemented in $L^q(L^{^\infty }(\Bbb{R)}\otimes \mathcal{M},L_c^2(%
\widetilde{\Gamma }))(1<q<\infty )$ via the embedding $\Phi .$ Hence, from
the interpolation result (\ref{interp}) we have 
\begin{equation}
(\mathcal{H}_c^q(\Bbb{R},\mathcal{M}),\mathcal{H}_c^{q^{\prime }}(\Bbb{R},%
\mathcal{M}))_{\frac 12}=\mathcal{H}_c^2(\Bbb{R},\mathcal{M})  \label{hh2}
\end{equation}
Combining (\ref{bh}), (\ref{hh}) and (\ref{hh2}) we get (\ref{HBMO}) by
Wolff's interpolation theorem (see \cite{w}). (\ref{hbmor}) can be proved
similarly. For (\ref{hbmocr}), by Lemma 6.1 and Theorem 5.1, 
\begin{eqnarray*}
&&(\mathrm{BMO}_{cr}(\Bbb{R},\mathcal{M}),L^1(L^\infty (\Bbb{R})\otimes 
\mathcal{M)})_{\frac 1p} \\
&=&(\mathrm{BMO}_{cr}^{\mathcal{D}}(L^\infty (\Bbb{R)}\otimes \mathcal{M}%
)\cap \mathrm{BMO}_{cr}^{\mathcal{D}^{\prime }}(L^\infty (\Bbb{R)}\otimes 
\mathcal{M}),L^1(L^\infty (\Bbb{R})\otimes \mathcal{M)})_{\frac 1p} \\
&\subseteq &(\mathrm{BMO}_{cr}^{\mathcal{D}}(L^\infty (\Bbb{R)}\otimes 
\mathcal{M}),L^1(L^\infty (\Bbb{R})\otimes \mathcal{M}))_{\frac 1p} \\
&&\qquad \cap (\mathrm{BMO}_{cr}^{\mathcal{D}^{\prime }}(L^\infty (\Bbb{R)}%
\otimes \mathcal{M}),L^1(L^\infty (\Bbb{R})\otimes \mathcal{M)})_{\frac 1p}
\\
&=&L^p(L^\infty (\Bbb{R})\otimes \mathcal{M)}
\end{eqnarray*}
On the other hand, since $\mathrm{BMO}_{cr}(\Bbb{R},\mathcal{M})\supset
L^\infty (L^\infty (\Bbb{R})\otimes \mathcal{M)},$ 
\begin{eqnarray*}
&&(\mathrm{BMO}_{cr}(\Bbb{R},\mathcal{M}),L^1(L^\infty (\Bbb{R})\otimes 
\mathcal{M)})_{\frac 1p} \\
&\supseteq &(L^\infty (L^\infty (\Bbb{R})\otimes \mathcal{M)},L^1(L^\infty (%
\Bbb{R})\otimes \mathcal{M)})_{\frac 1p} \\
&=&L^p(L^\infty (\Bbb{R})\otimes \mathcal{M)}.
\end{eqnarray*}
Therefore, 
\[
(\mathrm{BMO}_{cr}(\Bbb{R},\mathcal{M}),L^1(L^\infty (\Bbb{R})\otimes 
\mathcal{M)})_{\frac 1p}=L^p(L^\infty (\Bbb{R})\otimes \mathcal{M)}.
\]
By duality we have 
\[
(L^\infty (L^\infty (\Bbb{R})\otimes \mathcal{M}),\mathcal{H}_{cr}^1(\Bbb{R},%
\mathcal{M}))_{\frac 1p}=L^p(L^\infty (\Bbb{R})\otimes \mathcal{M)}.
\]
Finally, 
\begin{eqnarray*}
(L^\infty (L^\infty (\Bbb{R})\otimes \mathcal{M}),\mathcal{H}_{cr}^1(\Bbb{R},%
\mathcal{M}))_{\frac 1p} &\subseteq &(\mathrm{BMO}_{cr}(\Bbb{R},\mathcal{M}),%
\mathcal{H}_{cr}^1(\Bbb{R},\mathcal{M}))_{\frac 1p} \\
&\subseteq &(\mathrm{BMO}_{cr}(\Bbb{R},\mathcal{M}),L^1(L^\infty (\Bbb{R}%
)\otimes \mathcal{M}))_{\frac 1p}.
\end{eqnarray*}
Hence 
\[
(\mathrm{BMO}_{cr}(\Bbb{R},\mathcal{M}),\mathcal{H}_{cr}^1(\Bbb{R},\mathcal{M%
}))_{\frac 1p}=L^p(L^\infty (\Bbb{R})\otimes \mathcal{M}).
\]
Thus we have obtained all equalities in the theorem.\qed
\medskip 

\noindent\textbf{Remark.} We know little about $(\mathrm{BMO}_c(\Bbb{R},%
\mathcal{M}),L^1(L^\infty (\Bbb{R})\otimes \mathcal{M})_{\frac 1p}$ even for 
$p=2.$

\section{The real interpolation}

The following theorem is devoted to the real interpolation.

\begin{theorem}
Let $1\leq p<\infty .$ Then with equivalent norms, 
\begin{equation}
(X,Y)_{\frac 1p,p}=L^p(L^\infty (\Bbb{R})\otimes \mathcal{M}).  \label{chbmo}
\end{equation}
where $X=\mathrm{BMO}_{cr}(\Bbb{R},\mathcal{M})$ or $L^\infty (L^\infty (%
\Bbb{R})\otimes \mathcal{M})$ and $Y=\mathcal{H}_{cr}^1(\Bbb{R},\mathcal{M}%
)\,$or $L^1(L^\infty (\Bbb{R})\otimes \mathcal{M}).$
\end{theorem}

\noindent\textbf{Proof.} By Theorem 4.3 of \cite{MM} and Theorem 5.1 we
have(using the same argument as above for the complex method) 
\[
(\mathrm{BMO}_{cr}(\Bbb{R},\mathcal{M}),L^1(L^\infty (\Bbb{R})\otimes 
\mathcal{M}))_{\frac 1p,p}\subseteq L^p(L^\infty (\Bbb{R})\otimes \mathcal{M}%
). 
\]
On the other hand, for $1<p<\infty ,$ 
\begin{eqnarray*}
(\mathrm{BMO}_{cr}(\Bbb{R},\mathcal{M}),L^1(L^\infty (\Bbb{R})\otimes 
\mathcal{M}))_{\frac 1p,p} &\supseteq &(L^\infty (L^\infty (\Bbb{R})\otimes 
\mathcal{M}),L^1(L^\infty (\Bbb{R})\otimes \mathcal{M}))_{\frac 1p,p} \\
&=&L^p(L^\infty (\Bbb{R})\otimes \mathcal{M}).
\end{eqnarray*}
Therefore 
\[
(\mathrm{BMO}_{cr}(\Bbb{R},\mathcal{M}),L^1(L^\infty (\Bbb{R})\otimes 
\mathcal{M}))_{\frac 1p,p}=L^p(L^\infty (\Bbb{R})\otimes \mathcal{M}),\quad
1<p<\infty . 
\]
By duality we have 
\[
(L^\infty (L^\infty (\Bbb{R})\otimes \mathcal{M}),\mathcal{H}_{cr}^1(\Bbb{R},%
\mathcal{M}))_{\frac 1p,p}=L^p(L^\infty (\Bbb{R})\otimes \mathcal{M}),\quad
1<p<\infty . 
\]
Noting again that 
\begin{eqnarray*}
(L^\infty (L^\infty (\Bbb{R})\otimes \mathcal{M}),\mathcal{H}_{cr}^1(\Bbb{R},%
\mathcal{M}))_{\frac 1p,p} &\subseteq &(\mathrm{BMO}_{cr}(\Bbb{R},\mathcal{M}%
),\mathcal{H}_{cr}^1(\Bbb{R},\mathcal{M}))_{\frac 1p,p} \\
&\subseteq &(\mathrm{BMO}_{cr}(\Bbb{R},\mathcal{M}),L^1(L^\infty (\Bbb{R}%
)\otimes \mathcal{M}))_{\frac 1p,p},
\end{eqnarray*}
we conclude 
\[
\mathrm{BMO}_{cr}(\Bbb{R},\mathcal{M}),\mathcal{H}_{cr}^1(\Bbb{R},\mathcal{M}%
))_{\frac 1p,p}=L^p(L^\infty (\Bbb{R})\otimes \mathcal{M})),\quad 1<p<\infty
.\qed
\]

\medskip

\noindent\textbf{Remark.}Very recently, Junge and Musat got a
John-Nirenberg theorem for BMO spaces of non-commutative martingales 
(see \cite{JM}). By using Proposition 3.1 and the duadic trick of this article, they
got a John-Nirenberg theorem for non-commutative BMO spaces discussed here,
which can also be proved as a consequence of the interpolation results
established in this chapter. Unlike the classical case, the John-Nirenberg
theorem for non-commutative BMO spaces will no longer be the equivalence of 
\[
\sup_{I\subset \Bbb{R}}\left\| \left( \frac 1{|I|}\int_I|{\varphi }-{\varphi 
}_I|^pd\mu \right) ^{\frac 1p}\right\| _{\mathcal{M}} 
\]
for different $p,1\leq p\leq \infty $. In fact, if $\mathcal{M}=M_n$ the
algebra of $n$ by $n$ matrices, it can be proved that the best constant $%
c_n$ such that 
\[
\sup_{I\subset \Bbb{R}}\left\| \frac 1{|I|}\int_I|{\varphi }-{\varphi }%
_I|^2d\mu \right\| _{M_n}^{\frac 12}\leq c_n\sup_{I\subset \Bbb{R}}\left\|
\frac 1{|I|}\int_I|{\varphi }-{\varphi }_I|d\mu \right\| _{M_n},\ 
\]
holds for ${\varphi \in \rm{BMO}}_c{(}\Bbb{R}{,M}_n)$ will be at least $c\log n$
as $n\rightarrow \infty .$ And the corresponding constant for $M_n$ valued martingales
could be $cn^{\frac12}$ if no additional assumption on the related filtration. 
What remains true is the equivalence of 
\[
\sup_{I\subset \Bbb{R}}\sup_{\tau |a|^p\leq 1,}|I|^{-\frac 1p}\left\|
(f-f_I)a\chi _I\right\| _{L^p(\Bbb{R},\mathcal{M})}+\sup_{I\subset \Bbb{R}%
}\sup_{\tau |a|^p\leq 1,}|I|^{-\frac 1p}\left\| a\chi _I(f-f_I)\right\|
_{L^p(\Bbb{R},\mathcal{M})} 
\]
for different $p,2\leq p<\infty $ (see Theorem 1.2 of \cite{JM}) and
the equivalence of 
\[
\sup_{\mbox{cube }I\subset \Bbb{R}}\sup_{\tau |a|^p\leq 1,}\{|I|^{-\frac
1p}\left\| (f-f_I)a\chi _I\right\| _{\mathcal{H}_c^p(\Bbb{R},\mathcal{M})}\} 
\]
for different $p$, $2\leq p<\infty .$ See \cite{JM}, \cite{m2} for more information on
this.

\section{Fourier multipliers}

We close this chapter by a result on Fourier multipliers. Recall that $H^1(%
\Bbb{R})$ denotes the classical Hardy space on $\Bbb{R}.$ We will also need $%
H^1(\Bbb{R},H),$ the $H^1$ on $\Bbb{R}$ with values in a Hilbert space $H.$
Recall that we say a bounded map $M:H^1(\Bbb{R})\Bbb{\rightarrow }H^1(\Bbb{R}%
)$ is a Fourier multiplier if there exists a function $m\in L^\infty (\Bbb{R}%
)$ such that 
\[
\widehat{Mf}=m\widehat{f},\quad \forall f\in H^1(\Bbb{R}) 
\]
where $\widehat{f}$ is the Fourier transform of $f.$

\begin{theorem}
Let $M$ be a Fourier multiplier of the classical Hardy space $H^1(\Bbb{R}).$
Then $M$ extends in a natural way to a bounded map on $\mathrm{BMO}_c(\Bbb{R}%
,\mathcal{M})$ and $\mathcal{H}_c^p(\Bbb{R},\mathcal{M})$ for all $1\leq
p<\infty $ and 
\begin{eqnarray}
\left\| M:\mathrm{BMO}_c(\Bbb{R},\mathcal{M})\rightarrow \mathrm{BMO}_c(\Bbb{%
R},\mathcal{M})\right\|  &\leq c\left\| M:H^1(\Bbb{R})\rightarrow H^1(\Bbb{R%
})\right\| ,  \label{m} \\
\left\| M:\mathcal{H}_c^p(\Bbb{R},\mathcal{M})\rightarrow \mathcal{H}_c^p(%
\Bbb{R},\mathcal{M})\right\|  &\leq c\left\| M:H^1(\Bbb{R})\rightarrow H^1(%
\Bbb{R})\right\| .  \label{Mh}
\end{eqnarray}
Similar assertions also hold for $\mathrm{BMO}_r(\Bbb{R},\mathcal{M}),%
\mathrm{BMO}_{cr}(\Bbb{R},\mathcal{M}),\mathcal{H}_c^p(\Bbb{R},\mathcal{M})$
and $\mathcal{H}_{cr}^p(\Bbb{R},\mathcal{M})$.
\end{theorem}

\noindent\textbf{Proof. } Assume $\left\| M:H^1(\Bbb{R})\rightarrow H^1(\Bbb{%
R})\right\| =1.$ Let $H$ be the Hilbert space on which $\mathcal{M\,}$acts.
We start by showing the (well known) fact that $M$ is bounded on $H^1(\Bbb{R}%
,H).$ Denote by $R$ the Hilbert transform. Recall that $\left\| f\right\|
_{H^1(\Bbb{R},H)}\backsimeq \left\| f\right\| _{L^1(\Bbb{R},H)}+\left\|
Rf\right\| _{L^1(\Bbb{R},H)}$ for every $f\in $ $H^1(\Bbb{R},H).$ Denote by $%
\{e_{\lambda} \}_{{\lambda }\in {\Lambda }}$ the orthogonal normalized basis
of $H.$ Then $f=(f_{\lambda} )_{{\lambda }\in {\Lambda }}$ with $f_{\lambda}
=\langle e_{\lambda} ,f\rangle e_{\lambda} .$ Note that if $f\in $ $H^1(\Bbb{R}%
,H)$ then at most countably many $f_{\lambda} $'s are non zero. Let ${%
\varepsilon }=({\varepsilon }_n)_{n\in \Bbb{N}}$ be a sequence of
independent random variables on some probability space $(\Omega ,P)$ such
that $P({\varepsilon }_n=1)=P({\varepsilon }_n=-1)=\frac 12,$ $\forall n\in 
\Bbb{N}.\,$Notice that $MR=RM.$ Let $f\in $ $H^1(\Bbb{R},H).$ Let $\{{%
\lambda }_n:n\in \Bbb{N}\}$ be an enumeration of the ${\lambda }$'s such
that $f_{\lambda} \not{=}0$. Then by Khintchine's inequality, 
\begin{eqnarray*}
\left\| Mf\right\| _{H^1(\Bbb{R},H)} &\backsimeq &\int_{\Bbb{R}}((\sum_{n\in 
\Bbb{N}}|Mf_{{\lambda }_n}|^2)^{\frac 12}+(\sum_{n\in \Bbb{N}}|RMf_{{\lambda 
}_n}|^2)^{\frac 12})dt \\
&\backsimeq &\int_{\Bbb{R}}\int_\Omega |\sum_{n\in \Bbb{N}}{\varepsilon }%
_nMf_{{\lambda }_n}|dP({\varepsilon })dt+\int_{\Bbb{R}}\int_\Omega
|\sum_{n\in \Bbb{N}}{\varepsilon }_nMRf_{{\lambda }_n}|dP({\varepsilon })dt
\\
&\backsimeq &\int_\Omega \left\| M(\sum_{n\in \Bbb{N}}{\varepsilon }_nf_{{%
\lambda }_n})\right\| _{H^1(\Bbb{R},H)}dP({\varepsilon }) \\
&\leq &c\int_\Omega \left\| \sum_{n\in \Bbb{N}}{\varepsilon }_nf_{{\lambda }%
_n}\right\| _{H^1(\Bbb{R},H)}dP({\varepsilon }) \\
&\leq &c\left\| f\right\| _{H^1(\Bbb{R},H)}
\end{eqnarray*}
Therefore, as announced 
\[
\left\| M:H^1(\Bbb{R},H)\rightarrow H^1(\Bbb{R},H)\right\| \leq c_1. 
\]
Then by transposition 
\[
\left\| M:\mathrm{BMO}(\Bbb{R},H)\rightarrow \mathrm{BMO}(\Bbb{R},H)\right\|
\leq c_2; 
\]
whence, in virtue of (\ref{bmoh}), 
\[
\left\| M:\mathrm{BMO}_c(\Bbb{R},\mathcal{M})\rightarrow \mathrm{BMO}_c(\Bbb{%
R},\mathcal{M})\right\| \leq c_2. 
\]
Thus by duality 
\[
\left\| M:\mathcal{H}_c^1(\Bbb{R},\mathcal{M})\rightarrow \mathcal{H}_c^1(%
\Bbb{R},\mathcal{M})\right\| \leq c_3. 
\]
Then by Theorem 6.1 we have 
\[
\left\| M:\mathcal{H}_c^p(\Bbb{R},\mathcal{M})\rightarrow \mathcal{H}_c^p(%
\Bbb{R},\mathcal{M})\right\| \leq c_4. 
\]
Hence we have obtained the assertion concerning the column spaces. The other
assertions are immediate consequences of this one.\qed
\medskip

\medskip \noindent
\textbf{Remark.} As we have mentioned in the remark at the end of Chapter 2,
all the results of Chapter 2 can be generalized to the case of the functions
defined on $\Bbb{R}^n$ or $\Bbb{T}^n$, where $\Bbb{T}$ is the unit circle, and 
the relevant constants will be still absolute constants. By
Proposition 2.5 of \cite{m}, all the results carried out in other chapters
of this paper can also be generalized to the case of functions defined on 
$\Bbb{R}^n$ or $\Bbb{T}^n$ by the same approach as in Chapter
3, 5, 6. Unfortunately, the relevant constants there will depend on 
$n$ because the constant in Proposition 2.5 of \cite{m} depends on $n$.
This could be corrected if we could find a direct proof of the non-commutative
Hardy-Littlewood maximal inequality.

\newpage

\noindent\textbf{Acknowledgments.} Most of the work carried out in this paper was
done under the direction of Quanhua Xu. The author is greatly indebted to Q.
Xu for having suggested to him the subject of this paper, for many helpful
discussion and very careful reading of this paper. He has paid a lot of
attention to all parts of this paper and made many essential modifications
to its original version. The author is also very grateful to his current
advisor Gilles Pisier for his very careful reading of the manuscript and
very useful comments which led to many corrections and improvements. Finally, the author 
would like to thank
the referee for his very careful reading of the manuscript and very useful comments.

\medskip
{\small Math. Dept., Texas A\&M Univ., College Station, TX, 77843,
U. S. A. , and Math. Dept., WuHan Univ., 430072, 
P.R. China;  
email: tmei@math.tamu.edu
}
\end{document}